\documentclass[11pt]{article}
\usepackage{amsmath,amssymb,amsfonts,amsthm,authblk}
\usepackage{amsrefs}
\usepackage{braket}
\usepackage{moreverb,rotating,graphics,bm}
\usepackage[T1]{fontenc}
\usepackage{epsfig}
\usepackage{dsfont}
\usepackage{subcaption}
\usepackage{ulem}
\usepackage{cancel}
\usepackage{xcolor}
\usepackage{enumitem}
\usepackage{mathrsfs}
\usepackage{multirow}
\usepackage{algorithm}
\usepackage{algpseudocode}
\usepackage{booktabs}
\usepackage{hyperref}
\usepackage{tcolorbox}
\hypersetup{colorlinks,breaklinks}
\hypersetup{
	colorlinks=true,
	allcolors={black},
	urlcolor={black}
}

\tcbuselibrary{breakable, listings, skins}
\allowdisplaybreaks

\newtheorem{theorem}{Theorem}

\newtheorem{remark}{Remark}
\newtheorem{lemma}{Lemma}

\newtheorem{definition}{Definition}
\newtheorem{example}{Example}

\algrenewcommand\algorithmicrequire{\textbf{Input:}}
\algrenewcommand\algorithmicensure{\textbf{Output:}}

\def\tr{{\rm Tr \,}}

\def\1{{\mathds{1}}}

\newcommand{\revise}[1]{{\color{black}#1}}

\graphicspath{{./figures/}}

\DefineSimpleKey{bib}{primaryclass}{}
\DefineSimpleKey{bib}{archiveprefix}{}
\newcommand{\squarebracket}[1]{[#1]}
\newcommand{\arxivurl}[1]{\href{https://arxiv.org/abs/#1}{\texttt{arXiv: #1}}}
\renewcommand{\MR}[1]{~MR\href{https://mathscinet.ams.org/mathscinet/article?mr=#1}{#1}}

\BibSpec{arXiv}{%
  +{}{\PrintAuthors}{author}
  +{,}{ \textit}{title}
  +{,}{ \arxivurl}{eprint}
  +{}{ \squarebracket}{primaryclass}
  +{,}{ }{date}
  +{.}{}{transition}
}

\BibSpec{inproceedings}{%
  +{}{\PrintAuthors}{author}
  +{,}{ \textit}{title}
  +{,}{ }{booktitle}
  +{}{ }{volume}
  +{,}{ }{publisher}
  +{,}{ }{date}
  +{,}{ pp.~}{pages}
  +{.}{}{transition}
}

\newcommand{\tabincell}[2]{\begin{tabular}{@{}#1@{}}#2\end{tabular}}
\DeclareMathOperator*{\argmin}{arg\,min}

\newcommand{\rmb}{{\rm b}}

\newcommand{\rmo}{{\rm o}}
\newcommand{\upD}{\mathrm{D}}

\newcommand{\upN}{\mathrm{N}}
\newcommand{\upP}{\mathrm{P}}

\newcommand{\upR}{\revise{\mathrm{R}}}

\newcommand{\upT}{\mathrm{T}}
\newcommand{\dd}{\,\mathrm{d}}
\newcommand{\diag}{\mathrm{diag}}

\newcommand{\dist}{\mathrm{dist}}

\newcommand{\eps}{\varepsilon}

\newcommand{\Proj}{\mathrm{Proj}}
\newcommand{\Span}{\mathrm{span}}
\newcommand{\st}{\mathrm{s.~t.}}
\newcommand{\trace}{\mathrm{Tr}}
\newcommand{\two}{\mathrm{II}}

\newcommand{\calB}{\mathcal{B}}

\newcommand{\calE}{\mathcal{E}}

\newcommand{\calI}{\mathcal{I}}

\newcommand{\calL}{\mathcal{L}}
\newcommand{\calM}{\mathcal{M}}
\newcommand{\calN}{\mathcal{N}}
\newcommand{\calO}{\mathcal{O}}

\newcommand{\calS}{\mathcal{S}}

\newcommand{\calU}{\mathcal{U}}
\newcommand{\calV}{\mathcal{V}}

\newcommand{\calX}{\mathcal{X}}
\newcommand{\calY}{\mathcal{Y}}

\newcommand{\N}{\mathbb{N}}
\newcommand{\R}{\mathbb{R}}

\newcommand{\bdelta}{\boldsymbol{\delta}}

\newcommand{\veca}{\boldsymbol{a}}

\newcommand{\vecc}{\boldsymbol{c}}
\newcommand{\vecd}{\boldsymbol{d}}
\newcommand{\vece}{\boldsymbol{e}}

\newcommand{\vech}{\boldsymbol{h}}
\newcommand{\vecm}{\boldsymbol{m}}

\newcommand{\vecr}{\boldsymbol{r}}

\newcommand{\vecu}{\boldsymbol{u}}
\newcommand{\vecv}{\boldsymbol{v}}
\newcommand{\vecw}{\boldsymbol{w}}
\newcommand{\vecx}{\boldsymbol{x}}
\newcommand{\vecy}{\boldsymbol{y}}
\newcommand{\vecz}{\boldsymbol{z}}

\newcommand{\PT}{{\rm PT}}
\newcommand{\VT}{{\rm VT}}

\newcommand{\lrbrace}[1]{\left\{#1\right\}}
\newcommand{\lrbracket}[1]{\left(#1\right)}
\newcommand{\lrsquare}[1]{\left[#1\right]}

\newcommand{\abs}[1]{\left|#1\right|}
\newcommand{\inner}[1]{\left\langle#1\right\rangle}

\newcommand{\norm}[1]{\left\Vert#1\right\Vert}
\newcommand{\snorm}[1]{\Vert#1\Vert}

\newcommand{\Exp}{\mathrm{Exp}}
\newcommand{\Id}{\mathrm{Id}}
\newcommand{\Gr}{\mathrm{Gr}}
\newcommand{\St}{\mathrm{St}}

\newcommand{\Ver}{\mathrm{Ver}}
\newcommand{\Hor}{\mathrm{Hor}}
\newcommand{\grad}[1]{\mathrm{grad}_{#1}\,}
\newcommand{\hess}[1]{\mathrm{Hess}_{#1}\,}

\newcommand{\RHF}{\mathrm{RHF}}

\newcommand{\asym}{\mathrm{asym}}
\newcommand{\sym}{\mathrm{sym}}

\definecolor{lightblue}{rgb}{0.957,0.963,0.975}
\definecolor{midblue}{rgb}{0.937,0.943,0.965}
\definecolor{deepblue}{rgb}{0.325,0.427,0.569}
\definecolor{blocktitleblue}{rgb}{0.225,0.427,0.669}
\definecolor{lightred}{rgb}{0.996,0.969,0.969}
\definecolor{midred}{rgb}{0.976,0.949,0.949}
\definecolor{deepred}{rgb}{0.686,0.133,0.098}
\definecolor{deepgreen}{rgb}{0,0.5,0}
\definecolor{halfgray}{gray}{0.55}

\definecolor{lightpurple}{rgb}{0.978,0.978,1.0}
\definecolor{deeppurple}{rgb}{0.353,0.275,0.478}

\title{Constrained Dynamics for Searching Saddle Points on \revise{Embedded Riemannian Submanifolds of Euclidean Space}}

\author{Yukuan Hu}
\author{Laura Grazioli}
\affil{CERMICS, \'Ecole des Ponts - Institut Polytechnique de Paris and Inria, \linebreak 6-8 avenue Blaise Pascal, Cit\'e Descartes, 77455 Marne-la-Vall\'ee, France}

\begin{document}

\maketitle

\begin{abstract}
	Finding constrained saddle points on \revise{embedded Riemannian submanifolds of Euclidean space} is significant for analyzing energy landscapes arising in physics and chemistry. Existing works \revise{exploit explicit} global\revise{/local} regular level-set representations \revise{of manifolds, which may be unavailable or computationally inconvenient for manifolds represented through, e.g., projectors, factorizations, or rank constraints.} 
In this paper, we develop a constrained saddle dynamic\revise{\sout{s}} \revise{based on embedded-submanifold geometric primitives, completely avoiding the use of explicit representations}. \revise{In particular, o}ur dynamic\revise{\sout{s}} is formulated compactly on the Grassmann bundle of the tangent bundle. By analyzing the Grassmann bundle geometry, we 
rigorously establish the local linear stability of the dynamic\revise{\sout{s}} and the local linear convergence of the resulting algorithms. Remarkably, our analysis provides \revise{the first iterate convergence result for discretized algorithms to saddle points of prescribed indices in embedded-submanifold settings}. Moreover, by \revise{virtue of the Grassmann bundle formulation}, we remove unnecessary nondegeneracy assumptions on the eigenvalues of the Riemannian Hessian that are present in existing works. We also point out that locating saddle points can be more ill-condition\revise{ed} than finding local minimizers, and requires using nonredundant parametrizations. Finally, numerical experiments on linear eigenvalue problems and electronic excited-state calculations showcase the effectiveness of the proposed algorithms and corroborate the established local theory.
\end{abstract}

\section{Introduction}

\par Finding the saddle points (SPs) of potential energy functionals is a fundamental task in various scientific and engineering applications, particularly those involving energy landscape analysis. For example, the transition states between two (meta)stable states, which are crucial in physics \cite{goldstein1969viscous}, chemistry \cite{truhlar1996current}, and biology \cite{shakhnovich1996conserved}, can be identified as index-1 SPs \cite{ambrosetti1973dual}. Higher-index SPs are useful for collective, multi-mode, or concerted transitions and play a key role in the construction of solution landscape \cite{heidrich1986saddle,yin2021searching}. Here, an unconstrained index-$k$ SP is a critical point where the Euclidean Hessian has exactly $k$ negative eigenvalues. Additionally, \revise{embedded Riemannian submanifolds of Euclidean space} can naturally arise as constraint sets due to the incorporation of physical laws, such as in the Thomson problem \cite{mehta2016kinetic,thomson1904xxiv} and Bose-Einstein condensation \cite{altmann2025riemannian,bao2013mathematical,bose1924plancks,einstein1925quantentheorie}. In such cases, index-$k$ constrained SPs can be analogously defined by using the Riemannian gradient and Hessian. 

\par In comparison with finding local or global minimizers, locating SPs exhibits two distinct and major challenges: (1) there are always unknown descent directions at SPs, rendering off-the-shelf optimization methods unstable; (2) in general, it is impossible to construct a global merit function which is variationally minimized (or maximized) at SPs, posing significant difficulties for designing globally convergent numerical methods \cite{levitt2017convergence}. 

\par There have been numerous algorithmic developments concerning saddle search in \revise{\sout{the}} unconstrained settings, together with rigorous theoretical analyses. On the contrary, the exploration in \revise{\sout{the}} constrained settings \revise{requires explicit global/local level-set representations of submanifolds} (see Eq. \eqref{eqn:special Riemannian manifold} later for \revise{a} definition). This limitation rules out applications such as electronic excited-state calculations \cite{edelman1998geometry,grazioli2026critical2,payne1992iterative,vidal2024geometric}, where quotient structures \revise{and rank constraints arise in quantum-mechanical models, and explicit representations are unavailable or computationally inconvenient.} 
In some cases, omitting the quotient structures allows applying the existing methods directly. But unlike searching for local or global minimizers, this simplification can be detrimental in our context; see Example \ref{exm:condition number at index 1 Grassmann} and Section \ref{subsec:linear eigenvalue problem} later for an example involving the Stiefel and Grassmann manifolds. 

\par In this work, we develop a constrained saddle dynamic\revise{\sout{s}} \revise{based on embedded-submanifold geometric primitives and avoid the use of explicit level-set representations completely. We establish} the theoretical properties of the continuous dynamic\revise{\sout{s}} and the resulting discretized algorithms. 
We further demonstrate the effectiveness of \revise{the developed} algorithms through electronic excited-state calculations on standard benchmark molecular systems.

\subsection{Literature review}\label{subsec:literature review}

\par In the following, we review the existing works for finding SPs in both unconstrained and constrained settings. Most of them primarily focus on developing \textit{locally} convergent methods, with only few exceptions \cite{lelievre2024using,su2025improvedhighindexsaddledynamics}. Our discussions are confined to cases where objective values, first-order derivatives, and (approximate) Hessian-vector products are available, excluding direct applications of Newton-type methods \cite{baker1986algorithm,banerjee1985search,cerjan1981finding,simons1983walking}. 

\medskip

\par\noindent\textit{Unconstrained settings.} In these cases, the numerical methods for finding the index-1 SPs can be mainly categorized into two classes: single-ended and double-ended methods. The double-ended methods (also known as the path-finding or chain-of-states methods) \cite{henkelman2002methods,koslover2007comparison,liu2024stability} are fed with two candidates of local minimizers and target at finding the minimum energy path (which passes through an index-1 SP under certain conditions by the mountain pass theorem \cite{ambrosetti1973dual}). The single-ended ones (also known as the surface walking or eigenvector-following methods) \cite{olsen2004comparison} start with a single initial point without {\it a priori} knowledge about the final state. In this work, we focus on the class of single-ended methods\footnote{In fact, the double-ended methods are unsuitable for locating higher-index SPs by their nature and cannot be easily generalized.}, covering the activation-relaxation technique ({\it nouveau}) \cite{barkema1996event,barkema2001activation,cances2009some,crippen1971minimization,doye1997surveying,machado2011optimized,malek2000dynamics}, gentlest ascent dynamics \cite{e2011gentlest,quapp2014locating}, (shrinking) dimer methods \cite{henkelman1999dimer,kastner2008superlinearly,poddey2008dynamical,zhang2012shrinking,zhang2016optimization}, among others \cite{li2001minimax,miron2001step}. All the mentioned three methods can be set in the flow
\begin{equation}
	\begin{aligned}
		\frac{\dd\vecx}{\dd t}(t)&=-\Proj_{\vecv(t)^\perp}\big(\grad{}f(\vecx(t))\big)+\inner{\vecv(t),\grad{}f(\vecx(t))}\vecv(t)\\
		&=-R_{\vecv(t)}\big(\grad{}f(\vecx(t))\big),
	\end{aligned}
	\label{eqn:x-dynamics unconstrained}
\end{equation}
where $f:\calE\to\R$ is the objective (or energy) functional, $\inner{\bullet,\bullet}:\calE\times\calE\to\R$ is the inner product of the ambient Euclidean space, $\vecv(t)\in\calE$ is an additional \revise{normalized} direction variable, $\Proj_{\vecv(t)^\perp}$ denotes the orthogonal projection operator onto $\Span\{\vecv(t)\}^\perp$, defined as
\begin{equation*}
	\Proj_{\vecv(t)^\perp}(\vecu):=\vecu-\inner{\vecv(t),\vecu}\vecv(t),\quad\forall~\vecu\in\calE,
\end{equation*}
and $R_{\vecv(t)}$ represents the (Householder) reflection operator defined using $\vecv(t)$:
$$R_{\vecv(t)}(\vecu):=\vecu-2\inner{\vecv(t),\vecu}\vecv(t)=\Proj_{\vecv(t)^\perp}(\vecu)-\inner{\vecv(t),\vecu}\vecv(t),\quad\forall~\vecu\in\calE.$$
Briefly speaking, instead of following the gradient flow, which leads the trajectory to a local or global minimizer, the dynamic\revise{\sout{s}} \eqref{eqn:x-dynamics unconstrained} increases the objective value by climbing up along $\pm\vecv(t)$ (whose sign is determined by the angle between $\vecv(t)$ and $\grad{}f(\vecx(t))$), while decreases the value in all the directions perpendicular to $\vecv(t)$. An effective candidate for $\vecv(t)$ is the normalized eigenvector corresponding to the lowest eigenvalue of $\hess{}f(\vecx(t))$. In a neighborhood of a nondegenerate index-1 unconstrained SP, the dynamic\revise{\sout{s}} \eqref{eqn:x-dynamics unconstrained} can then be viewed as the gradient flow of a local strongly convex merit function \cite{gao2015iterative,gould2016dimer}, thereby providing a local stabilization for the index-1 SP. This underlies the activation-relaxation technique \revise{(\textit{nouveau})}. 
Instead of solving \revise{for the lowest eigenvector} directly, the gentlest ascent dynamics include\revise{\sout{s}} a direction dynamic\revise{\sout{s}}, 
$$\frac{\dd\vecv}{\dd t}(t)=-\Proj_{\vecv(t)^\perp}\big(\hess{}f(\vecx(t))[\vecv(t)]\big),$$
which follows from the Euler-Lagrange equation of \revise{the associated Rayleigh-Ritz minimization}. The (shrinking) dimer methods further approximate the Hessian-vector product through a central finite-difference scheme\revise{, in which the grid size} 
is sometimes driven to $0^+$ as $t\to+\infty$ \cite{zhang2012shrinking}. 

\par In recent years, there have been extensions for locating index-$k$ unconstrained SPs \cite{chu2025generalized,e2011gentlest,gao2015iterative,luo2025accelerated,luo2022convergence,quapp2014locating,yin2021searching,yin2019high}. For instance, the gentlest ascent dynamics can be generalized by incorporating $k$ direction dynamics \cite{yin2019high} \revise{($k<\dim(\calE)$)}, namely,
$$\frac{\dd\vecx}{\dd t}(t)=-R_{V(t)}\big(\grad{}f(\vecx(t))\big)$$
and
\begin{align}
	&\frac{\dd\vecv_i}{\dd t}(t)=\label{eqn:v-dynamics unconstrained index-k}\\
    &-\Proj_{\vecv_i(t)^\perp}\big(\hess{}f(\vecx(t))[\vecv_i(t)]\big)+2\sum_{j=1}^{i-1}\Proj_{\vecv_j(t)}\big(\hess{}f(\vecx(t))[\vecv_i(t)]\big),\nonumber\\
	&\qquad i=1,\ldots,k,\quad\text{with}\quad V(t):=(\vecv_1(t),\ldots,\vecv_k(t))~\text{\revise{orthonormalized}},\nonumber
\end{align}
where $R_{V(t)}$ is the reflection operator defined using $V(t)$:
$$R_{V(t)}(\vecu):=\vecu-2\sum_{j=1}^k\inner{\vecv_j(t),\vecu}\vecv_j(t),\quad\forall~\vecu\in\calE,$$
and $\Proj_{\vecv_j(t)}$ refers to the orthogonal projection operator onto $\Span\{\vecv_j(t)\}$ ($j=1,\ldots,k$). The second term on the right-hand side of Eq. \eqref{eqn:v-dynamics unconstrained index-k} is introduced to maintain the orthonormality condition $\inner{\vecv_i(t),\vecv_j(t)}=\delta_{ij}$ ($i,j=1,\ldots,k$), by combining the Lagrangian formalism with operator splitting. Some improvements have also been made for acceleration and stabilization, including those leveraging second-order information \cite{cances2009some,quapp2014locating,yin2019high,zhang2016optimization}, based on local merit functions \cite{gao2015iterative,gould2016dimer,gu2025iterative}, and incorporating additional inertial terms \cite{luo2025accelerated}. Some works have established the linear stability of the dynamic\revise{\sout{s}} at unconstrained SPs \cite{e2011gentlest,levitt2017convergence,luo2025accelerated,yin2019high,zhang2012shrinking} and the local convergence of the discretized algorithms \cite{gao2015iterative,gould2016dimer,levitt2017convergence,luo2025accelerated,luo2022convergence,zhang2012shrinking} under the nondegeneracy assumption. The error estimates for different discretization schemes can be found in \cite{luo2024semi,zhang2022error}. Recently, there have been some attempts dealing with stochastic and degenerate settings \cite{cui2025efficient,jiang2025nullspacepreservinghighindexsaddledynamics,shi2025stochastic}. A package has been designed for solution landscape exploration and construction based on the dynamic\revise{\sout{s}} \eqref{eqn:v-dynamics unconstrained index-k} \cite{liu2025saddlescape}.

\medskip

\par\noindent\textit{Constrained settings.} The exploration in this context remains rather limited \cite{li2015gentlest,liu2023constrained,yin2022constrained,zhang2012constrained,zhang2023discretization}. All these works \revise{exploit explicit global/local level-set representations of Riemannian submanifolds, e.g.}, 
\begin{equation}
	\calM:=\lrbrace{\vecx\in\R^n\mid\vecc(\vecx)=0}\quad\text{with}\quad \vecc(\vecx):=\big(c_1(\vecx),\ldots,c_q(\vecx)\big)^\top\in\R^q,
	\label{eqn:special Riemannian manifold}
\end{equation}
where $q<n$ and the functions $c_i$'s are smooth and regular, in that ${\rm rank}(\grad{}\vecc(\vecx))=q$ for any $\vecx\in\calM$ (note that $\grad{}\vecc(\vecx)\in\R^{q\times n}$). In this case, a constrained saddle dynamic\revise{\sout{s}}, targeting index-$k$ constrained SPs \revise{($k<\dim(\calM)$)}, can be derived using the Lagrangian function with operator splitting as follows \cite{yin2022constrained}:
\begin{equation}
    \frac{\dd\vecx}{\dd t}(t):=-R_{\vecx(t),V(t)}\big(\grad{\calM}f(\vecx(t))\big)
    \label{eqn:x-dynamics constrained existing}
\end{equation}
and
\begin{align}
	&~\frac{\dd\vecv_i}{\dd t}(t):=-\Proj_{\vecx(t),\vecv_i(t)^\perp}\big(\hess{\calM}f(\vecx(t))[\vecv_i(t)]\big)\label{eqn:v-dynamics constrained existing}\\
    &\quad+2\sum_{j=1}^{i-1}\Proj_{\vecx(t),\vecv_j(t)}\big(\hess{\calM}f(\vecx(t))[\vecv_i(t)]\big)+N\bigg(\vecx(t),\frac{\dd\vecx}{\dd t}(t),\vecv_i(t)\bigg),\nonumber\\
	&\qquad i=1,\ldots,k,\quad\text{with}~V(t):=(\vecv_1(t),\ldots,\vecv_k(t)),\nonumber
\end{align}
where the operators $\Proj_{\vecx(t),\vecv_i(t)^\perp}$, $\Proj_{\vecx(t),\vecv_j(t)}$, and $R_{\vecx(t),V(t)}$ are similarly defined on the tangent space $\upT_{\vecx(t)}\calM$, and
$$N\bigg(\vecx,\frac{\dd\vecx}{\dd t},\vecv_i\bigg):=-\grad{}\vecc(\vecx)^\top\lrbracket{\grad{}\vecc(\vecx)\cdot\grad{}\vecc(\vecx)^\top}^{-1}\lrbracket{\hess{}\vecc(\vecx)\lrsquare{\frac{\dd\vecx}{\dd t}}}\vecv_i.$$
Note that the invertibility of $\grad{}\vecc(\vecx)\cdot\grad{}\vecc(\vecx)^\top$ follows from the full rank assumption. We also recall that 
$$\hess{}\vecc(\vecx)\lrsquare{\vecu}=(\hess{}c_1(\vecx)[\vecu],\ldots,\hess{}c_q(\vecx)[\vecu])^\top\in\R^{q\times n},\quad\forall~\vecu\in\R^n.$$
Similar dynamics have been derived in \cite{li2015gentlest,liu2023constrained,zhang2012constrained}. In \cite{yin2022constrained}, the authors establish the linear stability of dynamic\revise{\sout{s}} at index-$k$ constrained SPs. That being said, we shall remark that \revise{it is difficult to apply the dynamics \eqref{eqn:x-dynamics constrained existing}-\eqref{eqn:v-dynamics constrained existing} directly to embedded submanifolds that lack easy-to-use explicit global/local representations, such as} 
the Grassmann manifolds, fixed-rank manifolds, and more general cases which find applications in, e.g., electronic excited-state calculations \cite{edelman1998geometry,grazioli2026critical2,payne1992iterative,vidal2024geometric}. Moreover, the local convergence \revise{to constrained SPs with prescribed indices has} not been investigated for the discretized algorithms, due to the complication of manifold settings. 

\subsection{Contributions}

\par In this article, we develop a constrained saddle dynamic\revise{\sout{s}} (CSD) for finding index-$k$ constrained SPs on \revise{any embedded Riemannian submanifold of Euclidean space}. Instead of tracking the lowest $k$-dimensional invariant subspace of the Riemannian Hessian with $k$ separate direction variables in  $\upT_{\vecx(t)}\calM$, we adopt an orthogonal projector $P(t):\upT_{\vecx(t)}\calM\to\upT_{\vecx(t)}\calM$ as a single variable, respecting the inherent quotient structure \revise{due to rotational invariance (see Section \ref{sec:algorithmic developments} later for details)}. The position-projector pair $(\vecx(t),P(t))$ is treated compactly using the CSD formulated on the {\it Grassmann bundle} of $k$-planes in the tangent bundle $\upT\calM$, denoted by $\Gr_k(\upT\calM)$ (see Eq. \eqref{eqn:Grassmann bundle} later for definition). By studying the geometry of $\Gr_k(\upT\calM)$ \revise{(Lemmas \ref{lem:bundle structures as embedded submanifolds}, \ref{lem:tangent space of Grassmann bundle}, and \ref{lem:general retraction on Grassmann bundle})}, we reveal that the time derivative of $P(t)$ necessarily contains a term defined by the second fundamental form of the manifold, accounting for the varying tangent spaces along the trajectory. This term vanishes in \revise{unconstrained} settings (cf. Eq. \eqref{eqn:v-dynamics unconstrained index-k}) and, after horizontal lifts \revise{(Lemma \ref{lem:decomposition of tangent space of Stiefel bundle})}, reduces to the third term on the right-hand side of Eq. \eqref{eqn:v-dynamics constrained existing} for the submanifolds \eqref{eqn:special Riemannian manifold}. Notably, %
\revise{our derivations of CSD rely only on embedded-submanifold geometric primitives and are completely free of explicit global/local level-set representations.}

\par In theory, we establish the global well-definedness of the CSD (Theorem \ref{thm:global well-definedness}), its linear stability at index-$k$ constrained SPs (Theorem \ref{thm:linear stability of CSD}), and the {\it first} local linear convergence results for the discretized algorithm in the manifold-constrained settings (Theorem \ref{thm:local convergence of discretized algorithm}). Moreover, compared with the existing linear stability results, our analysis removes unnecessary nondegeneracy assumptions on eigenvalues by \revise{leveraging the Grassmann bundle formulation} 
(Remarks \ref{rem:weaker assumptions} and \ref{rem:removal of unnecessary assumptions}). We also demonstrate through an example (Example \ref{exm:condition number at index 1 Grassmann}) 
that finding constrained SPs can (1) be worse conditioned and (2) require choosing nonredundant parametrizations, in contrast to locating global or local minimizers.  

\par Finally, we demonstrate the effectiveness of the developed algorithm on linear eigenvalue problems (Section \ref{subsec:linear eigenvalue problem}) and electronic excited-state calculations for standard benchmark molecules (Section \ref{subsec:electronic excited-state calcs}). We also corroborate numerically the influence of problem data and the importance of removing parametrization redundancies when searching for SPs.

\medskip

\par\noindent {\it Organization.} This paper is organized as follows. 
In Section \ref{sec:algorithmic developments}, we first investigate the geometry of the Grassmann bundle $\Gr_k(\upT\calM)$, upon which the CSD is built. In Section \ref{sec:theoretical analysis}, we establish the theoretical properties of the CSD as well as its discretized version. In Section \ref{sec:numerical experiments} \revise{and Appendix \ref{appsec:numerical comparisons between retractions}}, we report the numerical results on linear eigenvalue problems and electronic excited-state calculations. The conclusions are drawn in Section \ref{sec:conclusions}. \revise{Some fundamental concepts of embedded Riemannian submanifolds are recalled in Appendix \ref{appsec:preliminaries on Riemannian manifolds} for self-containedness. }

\medskip

\par\noindent{\it Notations.} Throughout this paper, scalars, vectors, and matrices are usually denoted by lowercase, bold lowercase, and uppercase letters, respectively. The sets or spaces are presented by calligraphic letters. In particular, we write the spaces of all $k\times k$ real symmetric and \revise{antisymmetric} matrices as $\R_{\sym}^{k\times k}$ and $\R_{\asym}^{k\times k}$, respectively. The Stiefel manifold of $k$-frames and the Grassmann manifold of $k$-planes in a vector space $\calV$ are denoted by $\St_k(\calV)$ and $\Gr_k(\calV)$, respectively. The orthogonal group of degree $k$ is given by $\calO(k)$. \revise{For a Euclidean space $\calE$, $\calS(\calE)$ denotes the Euclidean space of self-adjoint linear operators on $\calE$. The notation $\calL(\calX,\calY)$ refers to the vector space of linear maps from $\calX$ to $\calY$.} The notation ``${\rm cl}$'' means taking the closure of a set. The notations ``$\inner{\bullet,\bullet}$'' and ``$\norm{\bullet}$'' calculate the inner product and norm of vectors in the ambient space. The notation ``$[\bullet,\bullet]$'' represents the commutator of two matrices, defined as $[A,B]:=AB-BA$. The identity mapping over a vector space $\calV$ is denoted by $\Id_{\calV}$. We write the orthogonal projection operator onto a vector space $\calV$ as $\Proj_{\calV}$; if $\calV=\Span\{\vecv\}$ or $\Span\{\vecv\}^\perp$ for some vector $\vecv$, we simply write $\Proj_{\vecv}$ or $\Proj_{\vecv^\perp}$. \revise{For $V:=(\vecv_1,\ldots,\vecv_k)$, $\vecv_i$ refers to the $i$-th component of $V$.} The reflection operator defined by $V=(\vecv_1,\ldots,\vecv_k)$ is denoted by $R_V$; if $k=1$, we simply write $R_{\vecv}$. 

\par For a Riemannian manifold $\calM$ with $\vecx\in\calM$, the tangent space to $\calM$ at $\vecx$ is denoted by $\upT_{\vecx}\calM$, the normal space to $\calM$ at $\vecx$ by $\upN_{\vecx}\calM$, the tangent bundle of $\calM$ by $\upT\calM$, any retraction on $\calM$ by $\upR$, the exponential mapping in particular by $\Exp$, the Riemannian distance by $\dist_{\calM}$, \revise{any vector transport by $\VT$, the parallel transport in particular by $\PT$,} and the second fundamental form at $\vecx$ by $\two_{\vecx}$. These notations are sometimes equipped with additional superscripts to indicate manifolds. If $\calM$ is a quotient manifold and $\overline{\calM}$ is its total space, with $\pi$ the associated quotient map, the tangent space $\upT_{\vecx}\overline{\calM}$ to $\overline{\calM}$ at $\vecx$ can be decomposed into vertical and horizontal subspaces, denoted by $\Ver_{\vecx}^\pi\,\overline{\calM}$ and $\Hor_{\vecx}^\pi\,\overline{\calM}$, respectively. A general Riemannian metric on $\upT_{\vecx}\calM$ is denoted by $\inner{\bullet,\bullet}_{\vecx}$. The orthogonal projection and reflection operators defined on $\upT_{\vecx}\calM$ are described by an additional subscript $\vecx$, e.g., $\Proj_{\vecx,\vecv}$ and $R_{\vecx,\vecv}$. 
For a smooth function $f$, we write its differential, Euclidean gradient, and Euclidean Hessian as $\upD f$, $\grad{}f$, and $\hess{}f$, respectively. If it is defined on a Riemannian manifold $\calM$, its Riemannian gradient and Hessian are denoted by $\grad{\calM}f$ and $\hess{\calM}f$, respectively. 

\par The notation ``$\oplus$'' stands for the direct sum of two vector spaces, and ``$\cong$'' for the diffeomorphism between two vector spaces. When describing algorithms, we use superscripts within brackets to refer to the iteration numbers. \revise{For notational ease, we sometimes use ``h.c.'' to denote the Hermitian conjugate of the part in front.}

\section{Algorithmic developments}\label{sec:algorithmic developments}

\par In this section, we develop a constrained saddle dynamic\revise{\sout{s}} (CSD) for locating index-$k$ constrained SPs, \revise{applicable to any embedded Riemannian submanifold of Euclidean space}. Existing derivations \cite{li2015gentlest,liu2023constrained,yin2022constrained} \revise{exploit explicit level-set representations of submanifolds} (see Eq. \eqref{eqn:special Riemannian manifold}) and introduce $k$ different direction dynamics to track the lowest $k$-dimensional invariant subspace of the Riemannian Hessian. 
However, it is not difficult to observe the quotient structure under the hood: only the subspace is pursued and it is invariant under the choice of orthonormal basis. Therefore, we instead adopt a single dynamic\revise{\sout{s}} of an orthogonal projector, which, together with the position dynamic\revise{\sout{s}}, will amount to the CSD on the Grassmann bundle $\Gr_k(\upT\calM)$ (see Section \ref{subsec:development of CSD}). To this end, we first \revise{introduce $\Gr_k(\upT\calM)$ and the closely related Stiefel bundle of tangent bundle $\St_k(\upT\calM)$}, and investigate \revise{their} geometries (see Section \ref{subsec:geometries of Stiefel and Grassmann bundles}). 
\revise{We refer readers to Appendix \ref{appsec:preliminaries on Riemannian manifolds} for some fundamentals on embedded Riemannian submanifolds if necessary.} For \revise{simplicity}, we assume \revise{in this section} that $\calE=\R^n$ \revise{with $n\in\N$}. 
Nevertheless, our arguments apply to \revise{cases of general Euclidean spaces}. 

\subsection{$\St_k(\upT\calM)$ and $\Gr_k(\upT\calM)$, and their geometries}\label{subsec:geometries of Stiefel and Grassmann bundles}

\revise{The Stiefel bundle of $k$-frames and the Grassmann bundle of $k$-planes in the tangent bundle $\upT\calM$ are respectively defined by
\begin{align}
	\St_k(\upT\calM)&:=\lrbrace{(\vecx,V)\mid\vecx\in\calM,~V\in\St_k(\upT_{\vecx}\calM)},\label{eqn:Stiefel bundle}\\
	\Gr_k(\upT\calM)&:=\lrbrace{(\vecx,P)\mid\vecx\in\calM,~P\in\Gr_k(\upT_{\vecx}\calM)},\label{eqn:Grassmann bundle}
\end{align}
where $\St_k(\upT_{\vecx}\calM)$ and $\Gr_k(\upT_{\vecx}\calM)$ stand for the Stiefel manifold of ordered $k$-tuples of orthonormal vectors in $\upT_{\vecx}\calM$ and the Grassmann manifold of $k$-dimensional linear subspaces of $\upT_{\vecx}\calM$, respectively. In some contexts, $\St_k(\upT_{\vecx}\calM)$ and $\Gr_k(\upT_{\vecx}\calM)$ are called fibers over $\vecx\in\calM$. When $k=1$, $\Gr_1(\upT\calM)=:\mathbb{P}(\upT\calM)$ is called the projective bundle of $\upT\calM$. In some references (e.g., \cite{saunders1989geometry}), $\Gr_k(\upT\calM)$ is also called contact bundle.} 
\revise{Fiberwise, $\St_k(\upT\calM)$ and $\Gr_k(\upT\calM)$ possess the same geometry as the classical manifolds $\St_k(\R^d)$ and $\Gr_k(\R^d)$, respectively, where $d:=\dim(\calM)$. But unlike the classical ones, these fibers vary with the base point, which produces additional subtleties (vide infra).}

\par \revise{The two bundle structures are complicated in their definitions, particularly because the second components move with base points. This may necessitate the use of intrinsic local coordinate charts when deriving their geometries. Nevertheless, we show that both $\St_k(\upT\calM)$ and $\Gr_k(\upT\calM)$ are embedded Riemannian submanifolds of some ambient Euclidean spaces. This result enables us later to study them in an extrinsic way.

\begin{lemma}[$\St_k(\upT\calM)$ and $\Gr_k(\upT\calM)$ as embedded submanifolds]\label{lem:bundle structures as embedded submanifolds}
    Let $\calM$ be an embedded Riemannian submanifold of $\R^n$ with $\dim(\calM)=d$ and let $1\le k<d$. Then the following two statements hold.
    \begin{enumerate}[label=\textup{(\roman*)}]
        \item $\St_k(\upT\calM)$ is an embedded Riemannian submanifold of $\R^n\times\R^{n\times k}$ with 
        $$\dim(\St_k(\upT\calM))=d+kd-\frac{k(k+1)}{2}.$$

        \item $\Gr_k(\upT\calM)$ is an embedded Riemannian submanifold of $\R^n\times\R_{\rm sym}^{n\times n}$ with
        $$\dim(\Gr_k(\upT\calM))=d+k(d-k).$$
    \end{enumerate}
\end{lemma}

\begin{proof}
    Let $\bar\vecx\in\calM$. Choose an open neighborhood $\calU\subseteq\calM$ of $\vecx$ and a smooth orthonormal frame, $\vece_1(\vecx)$, $\ldots$, $\vece_d(\vecx)$, of $\upT\calM$ over $\calU$. Let $E(\vecx):=\big(\vece_1(\vecx),\ldots,\vece_d(\vecx)\big)\in\R^{n\times d}$. Then it holds that $E(\vecx)^\top E(\vecx)=I_d$ and $E(\vecx)E(\vecx)^\top=\Proj_{\upT_{\vecx}\calM}$.
    
    \smallskip
    
    \par\noindent (i) Define 
    $$F_1:\calU\times\St_k(\R^d)\to\calM\times\R^{n\times k},\quad F_1(\vecx,Y):=\big(\vecx,E(\vecx)Y\big).$$
    Then we have 
    $$F_1\big(\calU\times\St_k(\R^d)\big)=\St_k(\upT\calM)\cap(\calU\times\R^{n\times k}).$$
    The forward inclusion is straightforward. Conversely, for any $\vecx\in\calU$ and $V\in\St_k(\upT_{\vecx}\calM)$, let $Y:=E(\vecx)^\top V\in\R^{d\times k}$, we have $F_1(\vecx,Y)=(\vecx,E(\vecx)E(\vecx)^\top V)=(\vecx,V)$ and $Y^\top Y=V^\top E(\vecx)E(\vecx)^\top V=V^\top V=I_k$. Therefore, the equality holds. 
    
    The map $F_1$ is injective. Indeed, for any $(\vecx_1,V_1)$, $(\vecx_2,V_2)\in\St_k(\upT\calM)\cap(\calU\times\R^{n\times k})$, let $(\vecx_1,Y_1)$, $(\vecx_2,Y_2)\in\calU\times\St_k(\R^d)$ be their pre-images, respectively. Given $\vecx_1=\vecx_2=\vecx\in\calM$, it holds that $E(\vecx)Y_1=V_1=V_2=E(\vecx)Y_2$, which implies $Y_1=Y_2$ after left-multiplying $E(\vecx)^\top$. The inverse of $F_1$ is continuous, as it is precisely $F_1^{-1}(\vecx,V)=(\vecx,E(\vecx)^\top V)\in\calU\times\St_k(\R^d)$ for any $\vecx\in\calU$ and $V\in\St_k(\upT_{\vecx}\calM)$. Finally, the differential of $F_1$ at $(\vecx,Y)\in\calU\times\St_k(\R^d)$ is 
    $$\upD F_1(\vecx,Y)[\delta\vecx,\delta Y]=\big(\delta\vecx,\upD E(\vecx)[\delta\vecx]Y+E(\vecx)\delta Y\big),~\forall~(\delta\vecx,\delta Y)\in\upT_{(\vecx,Y)}\big(\calU\times\St_k(\R^d)\big).$$
    Then $\delta\vecx=0$ and $\upD E(\vecx)[\delta\vecx]Y+E(\vecx)\delta Y=0$ imply $E(\vecx)\delta Y=0$, which further yields $\delta Y=0$ after left-multiplying $E(\vecx)^\top$. Thus, $F_1$ is an immersion. As a result, $F_1$ smoothly embeds $\calU\times\St_k(\R^d)\cong\St_k(\upT\calM)\cap(\calU\times\R^{n\times k})$ into $\calM\times\R^{n\times k}$. Due to the arbitrariness of $\bar\vecx$, $\St_k(\upT\calM)$ is an embedded submanifold of $\calM\times\R^{n\times k}$. Its dimension is 
    $$\dim(\St_k(\upT\calM))=\dim(\calU\times\St_k(\R^d))=d+kd-\frac{k(k+1)}{2}.$$
    
    \smallskip

    \par\noindent (ii) Define 
    $$F_2:\calU\times\Gr_k(\R^d)\to\calM\times\R_{\sym}^{n\times n},\quad F_2(\vecx,W):=\big(\vecx,E(\vecx)WE(\vecx)^\top\big).$$
    Then we have
    $$F_2\big(\calU\times\Gr_k(\R^d)\big)=\Gr_k(\upT\calM)\cap(\calU\times\R_{\sym}^{n\times n}).$$
    Indeed, for the forward inclusion, pick any $(\vecx,W)\in\calU\times\Gr_k(\R^d)$, let $P:=E(\vecx)WE(\vecx)^\top$, we can verify that $P^\top=P$, that
    $$P^2=E(\vecx)WE(\vecx)^\top E(\vecx)WE(\vecx)^\top=E(\vecx)W^2E(\vecx)^\top=E(\vecx)WE(\vecx)^\top=P,$$
    that ${\rm Ran}(P)\subseteq\upT_{\vecx}\calM$ due to
    $$\Proj_{\upT_{\vecx}\calM}\circ P=E(\vecx)E(\vecx)^\top E(\vecx)WE(\vecx)^\top=E(\vecx)WE(\vecx)^\top=P,$$
    and that ${\rm rank}(P)={\rm rank}(W)=k$, as desired. Conversely, for any $\vecx\in\calU$ and $P\in\Gr_k(\upT_{\vecx}\calM)$, let $W:=E(\vecx)^\top PE(\vecx)$, we have
    $$F_2(\vecx,W)=\big(\vecx,E(\vecx)E(\vecx)^\top PE(\vecx)E(\vecx)^\top\big)=\big(\vecx,\Proj_{\upT_{\vecx}\calM}\circ P\circ\Proj_{\upT_{\vecx}\calM}\big)=(\vecx,P).$$

    Similar to the proof of (i), we can show the injectivity of $F_2$. The inverse of $F_2$ is precisely $F_2^{-1}(\vecx,P)=\big(\vecx,E(\vecx)^\top PE(\vecx)\big)\in\calU\times\Gr_k(\R^d)$ for any $\vecx\in\calU$ and $P\in\Gr_k(\upT_{\vecx}\calM)$, which is continuous. Finally, $F_2$ is also an immersion, as its differential at $(\vecx,W)\in\calU\times\Gr_k(\R^d)$ is
    \begin{multline*}
        \upD F_2(\vecx,W)[\delta\vecx,\delta W]=\big(\delta\vecx,\upD E(\vecx)[\delta\vecx]WE(\vecx)^\top\\+E(\vecx)\delta WE(\vecx)^\top+E(\vecx)W\upD E(\vecx)[\delta\vecx]^\top\big),\quad\forall~(\delta\vecx,\delta W)\in\upT_{(\vecx,W)}\big(\calU\times\Gr_k(\R^d)\big).
    \end{multline*}
    Therefore, $F_2$ smoothly embeds $\calU\times\Gr_k(\R^d)\cong\Gr_k(\upT\calM)\cap\big(\calU\times\R_{\sym}^{n\times n}\big)$ into $\calM\times\R_{\sym}^{n\times n}$. Again due to the arbitrariness of $\bar\vecx$, $\Gr_k(\upT\calM)$ is an embedded submanifold of $\calM\times\R_{\sym}^{n\times n}$. Its dimension is 
    $$\dim(\Gr_k(\upT\calM))=\dim(\calU\times\Gr_k(\R^d))=d+k(d-k).$$
\end{proof}}

\par \revise{Lemma \ref{lem:bundle structures as embedded submanifolds} allows us to investigate the geometries of $\St_k(\upT\calM)$ and $\Gr_k(\upT\calM)$ safely in an extrinsic way.} We start with the characterization of the tangent space to $\St_k(\upT\calM)$. 

\begin{lemma}[Tangent space to $\St_k(\upT\calM)$]
	Let $\calM$ be \revise{an embedded} Riemannian submanifold of \revise{$\R^n$} with $\dim(\calM)=d$ \revise{and let $1\le k<d$}. For any $(\vecx,V)\in\St_k(\upT\calM)$, the tangent space 
	\begin{equation}
		\upT_{(\vecx,V)}\St_k(\upT\calM)=\lrbrace{(\bdelta,\Gamma)\left|\begin{array}{l}
				\bdelta\in\upT_{\vecx}\calM,~A\in\R_{\asym}^{k\times k},~B\in\R^{(d-k)\times k}\\[0.1cm]
				\Gamma=VA+V_\perp B+\big(\two_{\vecx}(\bdelta,\vecv_1),\ldots,\two_{\vecx}(\bdelta,\vecv_k)\big)
			\end{array}\right.},
		\label{eqn:tangent space of Stiefel bundle}
	\end{equation}
	where $V_{\perp}\in\St_{d-k}(\upT_{\vecx}\calM)$ satisfies $V^\top V_\perp=0$. The characterization \eqref{eqn:tangent space of Stiefel bundle} is independent from the choice of $V_\perp$. 
\end{lemma}
\begin{proof}
	For any $(\bdelta,\Gamma)\in\upT_{(\vecx,V)}\St_k(\upT\calM)$, take a smooth curve $c:=(c_1,c_2):\R\supseteq\calI\to\St_k(\upT\calM)$ on $\St_k(\upT\calM)$ such that $c(0)=(\vecx,V)$ and $c'(0)=(\bdelta,\Gamma)$. It is obvious that $c_1$ is a smooth curve on $\calM$ passing through $\vecx$ at the origin, which implies $\bdelta\in\upT_{\vecx}\calM$. For the second part, note that $c_2(s)^\top c_2(s)=I_k$ and $c_2(s)_j=\Proj_{\upT_{c_1(s)}\calM}(c_2(s)_j)$, where $c_2(s)_j$ refers to the $j$-th column of $c_2(s)$ ($j=1,\ldots,k$). It therefore holds by \revise{Lemma \ref{lem:bundle structures as embedded submanifolds} and} the product rule that $\Gamma^\top V+V^\top\Gamma=0$ and 
	\begin{align}
		\Gamma_j=\left.\frac{\dd c_2(s)_j}{\dd s}\right|_{s=0}&=\left.\frac{\dd}{\dd s}\Proj_{\upT_{c_1(s)}\calM}\right|_{s=0}(\vecv_j)+\Proj_{\upT_{\vecx}\calM}(\Gamma_j)\label{eqn:Stiefel bundle tangent space 1}\\
		&=\Proj_{\upT_{\vecx}\calM}(\Gamma_j)+\two_{\vecx}(\bdelta,\vecv_j),\nonumber
	\end{align}
	for $j=1,\ldots,k$. Since the columns of $V$ and $V_\perp$ constitute an orthonormal basis of $\upT_{\vecx}\calM$, we have 
	$$\big(\Proj_{\upT_{\vecx}\calM}(\Gamma_1),\ldots,\Proj_{\upT_{\vecx}\calM}(\Gamma_k)\big)=VA+V_\perp B$$
	for some $A\in\R^{k\times k}$ and $B\in\R^{(d-k)\times k}$. Plugging this and Eq. \eqref{eqn:Stiefel bundle tangent space 1} into $\Gamma^\top V+V^\top\Gamma=0$ yields $A\in\R_{\asym}^{k\times k}$ because $V^\top V_{\perp}=0$ and $\two_{\vecx}(\bdelta,\vecv_j)\in\upN_{\vecx}\calM$ (see \revise{Appendix \ref{appsec:preliminaries on Riemannian manifolds}}). The proof is complete by noticing the arbitrariness of $(\bdelta,\Gamma)$.
\end{proof}

\par By the definitions in Eqs. \eqref{eqn:Stiefel bundle} and \eqref{eqn:Grassmann bundle}, there is a natural projection $\pi:\St_k(\upT\calM)\to\Gr_k(\upT\calM)$, defined as 
$$\St_k(\upT\calM)\owns(\vecx,V)\mapsto\pi(\vecx,V):=(\vecx,VV^\top)\in\Gr_k(\upT\calM),$$
which is smooth and surjective, and is a submersion; in fact, we have the quotient structure
\begin{equation}
	\Gr_k(\upT\calM)\cong\St_k(\upT\calM)/\calO(k).
	\label{eqn:Grassmann bundle quotient}
\end{equation}
Therefore, the tangent space to $\Gr_k(\upT\calM)$ at $(\vecx,P)$ can be readily obtained by computing the range of $\upD\pi(\vecx,V)$, where $V\in\St_k(\upT_{\vecx}\calM)$ satisfies $P=VV^\top$. This is given in the following lemma without proof. 

\begin{lemma}[Tangent space to $\Gr_k(\upT\calM)$]\label{lem:tangent space of Grassmann bundle}
	Let $\calM$ be \revise{an embedded Riemannian submanifold of $\R^n$} with $\dim(\calM)=d$ \revise{and let $1\le k<d$}. For any $(\vecx,P)\in\Gr_k(\upT\calM)$, the tangent space 
	\begin{align}
		&\upT_{(\vecx,P)}\Gr_k(\upT\calM)=\lrbrace{(\bdelta,\Delta)\mid\bdelta\in\upT_{\vecx}\calM,~\Delta=\Proj_{\upT_{\vecx}\calM}\circ\Delta\circ\Proj_{\upT_{\vecx}\calM}+\widehat\two_{\vecx}(\bdelta,P)}\label{eqn:tangent space of Grassmann bundle}\\
		&\qquad=\lrbrace{(\bdelta,\Delta)\mid\bdelta\in\upT_{\vecx}\calM,~\Delta=V_\perp BV^\top+VB^\top V_\perp^\top+\widehat\two_{\vecx}(\bdelta,P),~B\in\R^{(d-k)\times k}},\nonumber
	\end{align}
	where $V:=(\vecv_1,\ldots,\vecv_k)\in\St_k(\upT_{\vecx}\calM)$ and $V_\perp:=(\vecv_{k+1},\ldots,\vecv_d)\in\St_{d-k}(\upT_{\vecx}\calM)$ satisfy $P=VV^\top$ and $V^\top V_\perp=0$, the operator $\widehat\two_{\vecx}$ is defined by
	\begin{equation}
		\widehat\two_{\vecx}(\vecu,Q):=\sum_{\ell=1}^k\two_{\vecx}(\vecu,\vecu_\ell)\vecu_\ell^\top+\sum_{\ell=1}^k\vecu_\ell\two_{\vecx}(\vecu,\vecu_\ell)^\top,
		\label{eqn:extended second fundamental form}
	\end{equation}
	for any $\vecu\in\upT_{\vecx}\calM$ and $Q=UU^\top\in\Gr_k(\upT_{\vecx}\calM)$ with $U:=(\vecu_1,\ldots,\vecu_k)\in\St_k(\upT_{\vecx}\calM)$. The characterization \eqref{eqn:tangent space of Grassmann bundle} is independent from the choices of $V$ and $V_\perp$. 
	
	\par A natural basis for $\upT_{(\vecx,V)}\Gr_k(\upT\calM)$ is $\{(\bdelta_q,\Delta_q)\}_{q=1}^d\cup\{(\bdelta_{ij},\Delta_{ij})\}_{i=1,\ldots,k,~j=k+1,\ldots,d}$, where
	\begin{equation}
		(\bdelta_q,\Delta_q):=\lrbracket{\vecv_q,\widehat\two_{\vecx}(\vecv_q,P)},\quad q=1,\ldots,d,
		\label{eqn:horizontal part of tangent space of Grassmann bundle}
	\end{equation}
	and 
	\begin{equation}
		(\bdelta_{ij},\Delta_{ij}):=\lrbracket{0,\frac{1}{\sqrt{2}}\big(\vecv_i\vecv_j^\top+\vecv_j\vecv_i^\top\big)},\quad i=1,\ldots,k,~j=k+1,\ldots,d.
		\label{eqn:vertical part of tangent space of Grassmann bundle}
	\end{equation}
\end{lemma}

\par By virtue of the quotient structure \eqref{eqn:Grassmann bundle quotient}, each tangent space to $\St_k(\upT\calM)$ can be decomposed into vertical and horizontal parts induced by $\pi$.

\begin{lemma}[Decomposition of the tangent space to $\St_k(\upT\calM)$]\label{lem:decomposition of tangent space of Stiefel bundle}
	Let $\calM$ be \revise{an embedded Riemannian submanifold of $\R^n$} with $\dim(\calM)=d$ \revise{and let $1\le k<d$}. For any $(\vecx,V)\in\St_k(\upT\calM)$, the tangent space to $\St_k(\upT\calM)$ at $(\vecx,V)$ admits a direct sum decomposition as 
	$$\upT_{(\vecx,V)}\St_k(\upT\calM)=\Ver_{(\vecx,V)}^\pi\St_k(\upT\calM)\oplus\Hor_{(\vecx,V)}^\pi\St_k(\upT\calM),$$
	where
	\begin{align*}
		\Ver_{(\vecx,V)}^\pi\St_k(\upT\calM)&:=\lrbrace{(0,VA)\mid A\in\R_{\asym}^{k\times k}},\\
		\Hor_{(\vecx,V)}^\pi\St_k(\upT\calM)&:=\lrbrace{\lrbracket{\bdelta,V_\perp B+\big(\two_{\vecx}(\bdelta,\vecv_i)\big)_{i=1}^k}\mid\bdelta\in\upT_{\vecx}\calM,~B\in\R^{(d-k)\times k}}
	\end{align*}
	represent the vertical and horizontal subspaces of $\upT_{(\vecx,V)}\St_k(\upT\calM)$ induced by $\pi$, respectively, and $V_\perp\in\St_{d-k}(\upT_{\vecx}\calM)$ satisfies $V^\top V_\perp=0$. The characterization of the horizontal subspace is independent from the choice of $V_\perp$. Moreover, 
	$$\Hor_{(\vecx,V)}^\pi\St_k(\upT\calM)\cong\upT_{(\vecx,VV^\top)}\Gr_k(\upT\calM).$$
	For any $(\bdelta,\Delta)\in\upT_{(\vecx,VV^\top)}\Gr_k(\upT\calM)$, its horizontal lift is $(\bdelta,\Delta V)\in\Hor_{(\vecx,V)}^\pi\St_k(\upT\calM)$. 
\end{lemma}
\begin{proof}
	By definition, $\Ver_{(\vecx,V)}^\pi\St_k(\upT\calM)=\ker(\upD\pi(\vecx,V))$. Since 
	$$\upD\pi(\vecx,V)[(\bdelta,\Gamma)]=(\bdelta,V\Gamma^\top+\Gamma V^\top),\quad\forall~(\bdelta,\Gamma)\in\upT_{(\vecx,V)}\St_k(\upT\calM),$$
	and recall the characterization \eqref{eqn:tangent space of Stiefel bundle}, we have for $(\bdelta,\Gamma)\in\Ver_{(\vecx,V)}^\pi\St_k(\upT\calM)$ that $\bdelta=0$ and 
	$$V_\perp BV^\top+VB^\top V_\perp^\top=0\Leftrightarrow (V,~V_\perp)\begin{pmatrix}
		& B^\top\\
		B & 
	\end{pmatrix}\begin{pmatrix}
		V^\top \\ V_\perp^\top
	\end{pmatrix}=0.$$
	Note that the terms involving the second fundamental form vanish due to $\bdelta=0$. Since the columns of $V$ and $V_\perp$ form an orthonormal basis of $\upT_{\vecx}\calM$, the above equation implies that $B=0$. Therefore, $\Gamma=VA$ for some $A\in\R_{\asym}^{k\times k}$. This verifies the expression for the vertical subspace. Again by the characterization \eqref{eqn:tangent space of Stiefel bundle}, the horizontal subspace is clear since it is the orthogonal complement of the vertical subspace. The proof is complete. 
\end{proof}

\par It is known that a natural Riemannian metric for the tangent bundle is the Sasaki metric \cite{musso1988riemannian,sasaki1958differential}. Basically, the Sasaki metric uses the Riemannian connection of $\calM$ to split the tangent space to $\upT\calM$ at any point into horizontal and vertical subspaces, then defines the inner product using the original metric on each piece. Below, we define a Sasaki-type metric for $\Gr_k(\upT\calM)$. 

\begin{definition}[Sasaki-type metric for $\Gr_k(\upT\calM)$]\label{def:sasaki metric}
	Let $\calM$ be \revise{an embedded Riemannian submanifold of $\R^n$} with $\dim(\calM)=d$ \revise{and let $1\le k<d$}. For any $(\vecx,P)\in\Gr_k(\upT\calM)$, the Sasaki-type metric on $\upT_{(\vecx,P)}\Gr_k(\upT\calM)$ is defined as
	\begin{align*}
		\inner{(\bdelta,\Delta),(\tilde\bdelta,\tilde\Delta)}_{(\vecx,P)}:=\inner{\bdelta,\tilde\bdelta}+\inner{\Delta_{\vecx},\tilde\Delta_{\vecx}},\quad\forall~(\bdelta,\Delta),~(\tilde\bdelta,\tilde\Delta)\in\upT_{(\vecx,P)}\Gr_k(\upT\calM),
	\end{align*}
	where 
	\begin{equation}
		\Delta_{\vecx}:=\Proj_{\upT_{\vecx}\calM}\circ\Delta\circ\Proj_{\upT_{\vecx}\calM},~~\tilde\Delta_{\vecx}:=\Proj_{\upT_{\vecx}\calM}\circ\tilde\Delta\circ\Proj_{\upT_{\vecx}\calM}.
		\label{eqn:projected direction}
	\end{equation}
	Moreover, we denote the norm induced by the Sasaki-type metric with $\snorm{\bullet}_{(\vecx,P)}$. Then for any $(\bdelta,\Delta)\in\upT_{(\vecx,P)}\Gr_k(\upT\calM)$, 
	$$\norm{(\bdelta,\Delta)}_{(\vecx,P)}=\lrbracket{\snorm{\bdelta}^2+\snorm{\Delta_{\vecx}}^2}^{\frac12}.$$
\end{definition}

To preserve the feasibility condition in discretized algorithms, we need to define a retraction on $\Gr_k(\upT\calM)$. One possible choice based on the retraction on $\calM$ is given in the following lemma.

\revise{
\begin{lemma}[\revise{Retraction on $\Gr_k(\upT\calM)$}]\label{lem:general retraction on Grassmann bundle}
	Let $\calM$ be an embedded Riemannian submanifold of $\R^n$ with $\dim(\calM)=d$ and let $1\le k<d$. Suppose that
    \begin{enumerate}[label=\textup{(\roman*)}]
        \item $\upR^{\calM}:\upT\calM\to\calM$ is a retraction on $\calM$;

        \item for any $(\vecx,P)\in\Gr_k(\upT\calM)$ with $V\in\St_k(\upT_{\vecx}\calM)$ such that $P=VV^\top$, $\widehat\upR_{P,V}^{\Gr_k(\upT_{\vecx}\calM)}:\upT_P\big(\Gr_k(\upT_{\vecx}\calM)\big)\to\R^{n\times k}$ is smooth and fulfills 
        \begin{equation}
            \begin{aligned}
                &~\widehat\upR_{P,V}^{\Gr_k(\upT_{\vecx}\calM)}(0)=V,\\
                &~{\rm Ran}\Big(\widehat\upR_{P,V}^{\Gr_k(\upT_{\vecx}\calM)}(\Delta)\Big)\subseteq\upT_{\vecx}\calM,~~{\rm rank}\Big(\widehat\upR_{P,V}^{\Gr_k(\upT_{\vecx}\calM)}(\Delta)\Big)=k,\\
                &~\widehat\upR_{P,VQ}^{\Gr_k(\upT_{\vecx}\calM)}(\Delta)=\widehat\upR_{P,V}^{\Gr_k(\upT_{\vecx}\calM)}(\Delta)Q,\quad\forall~Q\in\calO(k),\\
                &~\upD\widehat\upR_{P,V}^{\Gr_k(\upT_{\vecx}\calM)}(0)[\Delta]V^\top+V\upD\widehat\upR_{P,V}^{\Gr_k(\upT_{\vecx}\calM)}(0)[\Delta]^\top=\Delta,
            \end{aligned}
            \label{eqn:condition on hat R_P}
        \end{equation}
        for any $\Delta\in\upT_P\big(\Gr_k(\upT_{\vecx}\calM)\big)$ sufficiently small;

        \item $\VT^{\calM}$ is a vector transport associated with $\upR^{\calM}$ on $\calM$ such that, for any $(\vecx,\bdelta)\in\upT\calM$, 
        \begin{equation}
            \upD\big(\iota\mapsto\VT_{(\vecx,\iota)}^{\calM}\big)(0)[\bdelta](\vecv)=\two_{\vecx}(\bdelta,\vecv)\in\upN_{\vecx}\calM
            \label{eqn:condition on vector transport}
        \end{equation}
        holds for any $\vecv\in\upT_{\vecx}\calM$, and, for any $\vecx\in\calM$, there exists a neighborhood $\calN_1$ of the origin in $\upT_{\vecx}\calM$ such that $\VT_{(\vecx,\bdelta)}^{\calM}:\upT_{\vecx}\calM\to\upT_{\upR_{\vecx}^{\calM}(\bdelta)}\calM$ is injective for any $\bdelta\in\calN_1$. 
    \end{enumerate}
    Then $\upR^{\Gr_k(\upT\calM)}:\upT(\Gr_k(\upT\calM))\to\Gr_k(\upT\calM)$ defined as
	\begin{equation}
	    \upR_{(\vecx,P)}^{\Gr_k(\upT\calM)}(\bdelta,\Delta):=\big(\upR_{\vecx}^{\calM}(\bdelta),\tilde P_{\vecx,P}(\bdelta,\Delta)\big),
        \label{eqn:retraction over Grassmann bundle}
	\end{equation}
    for any $\big((\vecx,P),(\bdelta,\Delta)\big)\in\upT\big(\Gr_k(\upT\calM)\big)$, with
    \begin{align}
        \tilde P_{\vecx,P}(\bdelta,\Delta):=&~\tilde V_{\vecx,P,V}(\bdelta,\Delta)\tilde V_{\vecx,P,V}(\bdelta,\Delta)^\top,\label{eqn:tilde P}\\
        \tilde V_{\vecx,P,V}(\bdelta,\Delta):=&~Y_{\vecx,P,V}(\bdelta,\Delta)\big(Y_{\vecx,P,V}(\bdelta,\Delta)^\top Y_{\vecx,P,V}(\bdelta,\Delta)\big)^{-1/2},\label{eqn:tilde V}\\
        Y_{\vecx,P,V}(\bdelta,\Delta):=&~\VT_{(\vecx,\bdelta)}^{\calM}\big(Z_{\vecx,P,V}(\Delta)\big),\label{eqn:Y}\\
        Z_{\vecx,P,V}(\Delta):=&~\widehat\upR_{P,V}^{\Gr_k(\upT_{\vecx}\calM)}(\Delta_{\vecx}),\label{eqn:Z}
    \end{align}
    is a retraction on $\Gr_k(\upT\calM)$, where $\VT_{(\vecx,\bdelta)}^{\calM}$ acts columnwisely in Eq. \eqref{eqn:Y} and $\Delta_{\vecx}$ is defined in Eq. \eqref{eqn:projected direction}. Its representative on $\St_k(\upT\calM)$ is
    \begin{equation}
        \upR_{(\vecx,V)}^{\Gr_k(\upT\calM)}(\bdelta,\Delta V):=\big(\upR_{\vecx}^{\calM}(\bdelta),\tilde V_{\vecx,P,V}(\bdelta,\Delta)\big).
        \label{eqn:retraction over Stiefel bundle}
    \end{equation}
\end{lemma}
\begin{proof}
    For the ease of notation and unless stated, we do not specify the reference point $(\vecx,P)$ or $(\vecx,P,V)$ in the subscripts of $\tilde P$, $\tilde V$, $Y$, and $Z$ in this proof. Given $(\vecx,P)\in\Gr_k(\upT\calM)$, we first prove that $\upR^{\Gr_k(\upT\calM)}_{(\vecx,P)}$in Eq. \eqref{eqn:retraction over Grassmann bundle} is well-defined in a neighborhood of the origin. Let $\calN_2\subseteq\upT_{\vecx}\calM$ be a neighborhood of the origin on which $\upR_{\vecx}^{\calM}$ is well-defined, and let $\calN_3\subseteq\upT_P\Gr_k(\upT_{\vecx}\calM)$ be the neighborhood of the origin (which is different from the previous one) such that Eq. \eqref{eqn:condition on hat R_P} holds for any $\Delta\in\calN_2$. 
    Suppose that $\bdelta\in\calN_1\cap\calN_2$ and $\Delta\in\calN_3$. The $\vecx$-part is obvious; hence, we focus on the $P$-part afterwards. From Eq. \eqref{eqn:Z} and the injectivity of $\VT_{(\vecx,\bdelta)}^{\calM}$, we have that $Y(\bdelta,\Delta)$ is of full column rank. Indeed, let $\veca\in\ker(Y(\bdelta,\Delta))$. Then by the linearity of $\VT_{(\vecx,\bdelta)}^{\calM}$, 
    $$0=Y(\bdelta,\Delta)\veca=\big(\VT_{(\vecx,\bdelta)}^{\calM}(\vecz_j(\Delta))\big)_{j=1}^k\veca=\VT_{(\vecx,\bdelta)}^{\calM}\big(Z(\Delta)\veca\big).$$
    Because of the injectivity of $\VT_{(\vecx,\bdelta)}^{\calM}$ on $\upT_{\vecx}\calM$, the above equality implies $\veca\in\ker(Z(\Delta))=\{0\}$, as desired. Consequently, $\tilde V(\bdelta,\Delta)$ is well-defined by Eq. \eqref{eqn:tilde V}, which in turn yields the well-definedness of $\tilde P(\bdelta,\Delta)$. The independence of $\tilde P(\bdelta,\Delta)$ on the choice of $V$ can be derived from the equivariance \eqref{eqn:condition on hat R_P}: for any $Q\in\calO(k)$, 
    \begin{align*}
        Z_{\vecx,P,VQ}(\Delta)&=\widehat\upR_{P,VQ}^{\Gr_k(\upT_{\vecx}\calM)}(\Delta_{\vecx})=\widehat\upR_{P,V}^{\Gr_k(\upT_{\vecx}\calM)}(\Delta_{\vecx})Q=Z_{\vecx,P,V}(\Delta)Q,\\
        Y_{\vecx,P,VQ}(\bdelta,\Delta)&=\VT_{(\vecx,\bdelta)}^{\calM}\big(Z_{\vecx,P,VQ}(\Delta)\big)=Y_{\vecx,P,V}(\bdelta,\Delta)Q,\\
        \tilde V_{\vecx,P,VQ}(\bdelta,\Delta)&=Y_{\vecx,P,V}(\bdelta,\Delta)Q\big(Q^\top Y_{\vecx,P,V}(\bdelta,\Delta)^\top Y_{\vecx,P,V}(\bdelta,\Delta)Q\big)^{-1/2}\\
        &=Y_{\vecx,P,V}(\bdelta,\Delta)QQ^\top\big(Y_{\vecx,P,V}(\bdelta,\Delta)^\top Y_{\vecx,P,V}(\bdelta,\Delta)\big)^{-1/2}Q\\
        &=Y_{\vecx,P,V}(\bdelta,\Delta)\big(Y_{\vecx,P,V}(\bdelta,\Delta)^\top Y_{\vecx,P,V}(\bdelta,\Delta)\big)^{-1/2}Q.
    \end{align*}
    From these arguments, it is also not difficult to verify that $\upR^{\Gr_k(\upT\calM)}$ is smooth and its codomain is $\Gr_k(\upT\calM)$. 
    
    
    \par Next we prove that $\upR^{\Gr_k(\upT\calM)}$ is indeed a retraction on $\Gr_k(\upT\calM)$. It is clear that $\upR_{(\vecx,P)}^{\Gr_k(\upT\calM)}(0,0)=(\vecx,P)$. Therefore, we focus on its differential at $(0,0)$ below. The first component is obvious. Regarding the second component, we have by Eqs. \eqref{eqn:tilde P}-\eqref{eqn:Z} and the chain rule that
    \begin{align}
        \upD\tilde P(0,0)[\bdelta,\Delta]&=\upD\tilde V(0,0)[\bdelta,\Delta]V^\top+V\upD\tilde V(0,0)[\bdelta,\Delta]^\top,\label{eqn:D tilde P}\\
        \upD\tilde V(0,0)[\bdelta,\Delta]&=\upD Y(0,0)[\bdelta,\Delta]+V\Theta(\bdelta,\Delta),\label{eqn:D tilde V}\\
        \upD Y(0,0)[\bdelta,\Delta]&=\lrbracket{\upD\big(\iota\mapsto\VT_{(\vecx,\iota)}^{\calM}\big)(0)[\bdelta](\vecv_j)}_{j=1}^k+\upD Z(0)[\Delta],\label{eqn:DY}\\
        \Delta_{\vecx}&=\upD Z(0)[\Delta]V^\top+V\upD Z(0)[\Delta]^\top,\label{eqn:D Z}
    \end{align}
    where 
    $$\Theta(\bdelta,\Delta):=\upD\lrbracket{(\iota,\Xi)\mapsto \big(Y(\iota,\Xi)^\top Y(\iota,\Xi)\big)^{-1/2}}(0,0)[\bdelta,\Delta],$$
    and Eq. \eqref{eqn:D Z} is due to condition \eqref{eqn:condition on hat R_P}. Note that by Eq. \eqref{eqn:DY},
    \begin{align*}
        &~\upD\lrbracket{(\iota,\Xi)\mapsto Y(\iota,\Xi)^\top Y(\iota,\Xi)}(0,0)[\bdelta,\Delta]\\
        =&~\upD Y(0,0)[\bdelta,\Delta]^\top V+V^\top\upD Y(0,0)[\bdelta,\Delta]\\
        =&~V^\top\lrbracket{\upD\big(\iota\mapsto\VT_{(\vecx,\iota)}^{\calM}\big)(0)[\bdelta](\vecv_j)}_{j=1}^k+V^\top\upD Z(0)[\Delta]+{\rm h.c.}=0,
    \end{align*}
    where the last equality uses Eq. \eqref{eqn:condition on vector transport} and the fact that $\upD Z(0)[\Delta]=VA+V_{\perp}B$ with $V_{\perp}\in\St_{d-k}(\upT_{\vecx}\calM)$ the orthonormal basis for the orthogonal complement of $\Span(V)$ in $\upT_{\vecx}\calM$, $A\in\R_{\asym}^{k\times k}$, and $B\in\R^{(d-k)\times k}$ (cf. \cite{bendokat2024grassmann}). Therefore, for a sufficiently small $t>0$, 
    $$Y(t\bdelta,t\Delta)^\top Y(t\bdelta,t\Delta)=I_k+O(t^2),$$
    which implies 
    $$\big(Y(t\bdelta,t\Delta)^\top Y(t\bdelta,t\Delta)\big)^{-1/2}=I_k+O(t^2)\quad\Rightarrow\quad\Theta(\bdelta,\Delta)=0.$$
    Combining this together with Eqs. \eqref{eqn:condition on vector transport} and \eqref{eqn:D tilde P}-\eqref{eqn:D Z}, we have 
    \begin{align*}
        \upD\tilde P(0,0)[\bdelta,\Delta]&=\upD Y(0,0)[\bdelta,\Delta]V^\top+{\rm h.c.}\\
        &=\lrbracket{\upD\big(\iota\mapsto\VT_{(\vecx,\iota)}^{\calM}\big)(0)[\bdelta](\vecv_j)}_{j=1}^kV^\top+\upD Z(0)[\Delta]V^\top+{\rm h.c.}\\
        &=\widehat\two_{\vecx}(\bdelta,P)+\Delta_{\vecx}=\Delta.
    \end{align*}
    The last equality follows from Lemma \ref{lem:tangent space of Grassmann bundle}. The proof is complete.
\end{proof}}

\revise{\begin{remark}[\revise{Examples of retractions on $\Gr_k(\upT\calM)$}]\label{rem:example retractions}
    The assumption (i) in Lemma \ref{lem:general retraction on Grassmann bundle} is not restrictive. In the following, we pay special attention to the last two assumptions (ii)-(iii) and provide concrete examples. 

    \par For the mapping in (ii), one natural choice is any $\calO(k)$-equivariant frame lift of a Grassmann retraction. Taking the exponential mapping for instance, we have from \cite{bendokat2024grassmann}
    $$Z(\Delta)=VU\cos(\Sigma)U^\top+W\sin(\Sigma)U^\top,\quad\forall~\Delta\in\upT_P\Gr_k(\upT_{\vecx}\calM),$$
    where $W\in\St_k(\upT_{\vecx}\calM)$, $\Sigma=\diag(\sigma_1,\ldots,\sigma_k)\in\R_{+}^{k\times k}$, and $U\in\calO(k)$ satisfy $\Delta V=W\Sigma U^\top$. Another (simple) strategy is to take 
    $$\widehat\upR_{P,V}^{\Gr_k(\upT_{\vecx}\calM)}(\Delta)=(P+\Delta)V\in\R^{n\times k},\quad\forall~(P,\Delta)\in\upT\big(\Gr_k(\upT_{\vecx}\calM)\big).$$
    Note that the range of $(P+\Delta)V$ is still a subspace of $\upT_{\vecx}\calM$, but $(P+\Delta)V$ is not necessarily orthonormalized. Nevertheless, this can be amended by the polar decomposition in Eq. \eqref{eqn:tilde V} later. 

    \par For the mapping in (iii), we provide the following two choices.
    \begin{itemize}
        \item Option I: parallel transport. In this case, $\VT_{(\vecx,\bdelta)}^{\calM}=\PT_{1\leftarrow0}^{c_1}$, where $c_1(t)\in\calM$ is the smooth curve defined by $c_1(t):=\upR_{\vecx}^{\calM}(t\bdelta)$. The invertibility of parallel transport implies the injectivity. Regarding its differential, we first note $\Proj_{\upT_{c_1(t)}\calM}\big(\VT_{(\vecx,t\bdelta)}^{\calM}(\vecv)\big)=\VT_{(\vecx,t\bdelta)}^{\calM}(\vecv)$. By the chain rule, 
        $$\upD\big(\iota\mapsto\VT_{(\vecx,\iota)}^{\calM}\big)(0)[\bdelta](\vecv)=\two_{\vecx}(\bdelta,\vecv)+\Proj_{\upT_{\vecx}\calM}\lrbracket{\upD\big(\iota\mapsto\VT_{(\vecx,\iota)}^{\calM}\big)(0)[\bdelta](\vecv)}.$$
        Condition \eqref{eqn:condition on vector transport} then holds after noticing that the second term on the right-hand side vanishes due to the property of parallel transport. 
    
        \par We remark that, with the vector transport chosen as the parallel transport and $Z(\Delta)\in\St_k(\upT_{\vecx}\calM)$, the definition of $\upR^{\Gr_k(\upT\calM)}$ can be simplified. Indeed, due to the fact that the parallel transport is an isometry between tangent spaces, it must hold that $Y(\bdelta,\Delta)^\top Y(\bdelta,\Delta)=I_k$. 
        In addition, if $\calM$ is Grassmannian and $\upR^{\calM}$ is chosen to be the exponential mapping on $\calM$, $\VT_{(\vecx,\bdelta)}^{\calM}$ can be explicitly constructed; see, e.g., \cite{bendokat2024grassmann}.

        \item Option II: metric projection-based vector transport. In this case, $\VT_{(\vecx,\bdelta)}^{\calM}=\Proj_{\upT_{c_1(1)}\calM}|_{\upT_{\vecx}\calM}$, where $c_1(t)$ is the same curve as the one defined in option I. Since $\ker(\VT_{(\vecx,\bdelta)}^{\calM})=\upT_{\vecx}\calM\cap\upN_{c(1)}\calM$, $\VT_{(\vecx,\bdelta)}^{\calM}$ is injective if and only if the largest principal angle between $\upT_{\vecx}\calM$ and $\upT_{c(1)}\calM$ is less than $\pi/2$. In \cite{boissonnat2019reach}, the authors provide an upper bound for this angle using the Riemannian distance between $\vecx=c(0)$ and $c(1)$, which implies that $\VT_{(\vecx,\bdelta)}^{\calM}$ in this option is injective as long as $\snorm{\bdelta}$ is sufficiently small. Regarding the differential, by Eq. \eqref{eqn:differential of orthogonal projector}, we have
        $$\upD\big(\iota\mapsto\VT_{(\vecx,\iota)}^{\calM}\big)(0)[\bdelta](\vecv)=\upP_{\vecx,\bdelta}(\vecv)=\two_{\vecx}(\bdelta,\vecv),$$
        which verifies condition \eqref{eqn:condition on vector transport}.
    \end{itemize}

    \par Provided that the assumptions in Lemma \ref{lem:general retraction on Grassmann bundle} are fulfilled, different definitions of retractions on $\Gr_k(\upT\calM)$ agree up to the first order and may differ locally from the second order. This discrepancy does not affect the framework of convergence analysis in Section \ref{sec:theoretical analysis}, which exploits merely first-order properties. Nevertheless, there can be distinctions in the involved constants and numerical performance (see Appendix \ref{appsec:numerical comparisons between retractions} for an illustration).
\end{remark}}

\subsection{Constrained saddle dynamics on $\Gr_k(\upT\calM)$}\label{subsec:development of CSD}

\par Equipped with the above geometrical tools of $\Gr_k(\upT\calM)$, we are ready to develop the constrained saddle dynamic\revise{\sout{s}} (CSD). Since $(\vecx(t),P(t))\in\Gr_k(\upT\calM)$ for any $t$, we require from Eq. \eqref{eqn:tangent space of Grassmann bundle} that
\begin{equation}
	\frac{\dd\vecx}{\dd t}(t)\in\upT_{\vecx(t)}\calM,~~\frac{\dd P}{\dd t}(t)=\Proj_{\upT_{\vecx(t)}\calM}\circ\frac{\dd P}{\dd t}(t)\circ\Proj_{\upT_{\vecx(t)}\calM}+\widehat\two_{\vecx}\lrbracket{\frac{\dd\vecx}{\dd t}(t),P(t)}
	\label{eqn:requirements on dx and dv}
\end{equation}
hold for all the time, where $\widehat\two_{\vecx}$ is defined in Eq. \eqref{eqn:extended second fundamental form}. To achieve this, we set
\begin{equation}
	\begin{aligned}
		\frac{\dd\vecx}{\dd t}(t)&:=-R_{\vecx(t),P(t)}\big(\grad{\calM}f(\vecx(t))\big),\\
		\frac{\dd P}{\dd t}(t)&:=-\Proj_{\upT_{P(t)}\Gr_k(\upT_{\vecx(t)}\calM)}\circ\hess{\calM}f(\vecx(t))+\widehat\two_{\vecx(t)}\lrbracket{\frac{\dd\vecx}{\dd t}(t),P(t)},
	\end{aligned}
	\label{eqn:deterministic constrained saddle dynamics}
\end{equation}
where $R_{\vecx(t),P(t)}$ is the reflection operator defined by $P(t)$ on $\upT_{\vecx(t)}\calM$:
\begin{equation}
	R_{\vecx(t),P(t)}(\vecu):=\vecu-2P(t)\vecu,\quad\forall~\vecu\in\upT_{\vecx(t)}\calM,
	\label{eqn:reflection operator in CSD}
\end{equation}
and $\Proj_{\upT_{P(t)}\Gr_k(\upT_{\vecx(t)}\calM)}$ refers to the orthogonal projection operator onto the tangent space to $\Gr_k(\upT_{\vecx(t)}\calM)$ at $P(t)$, which has a closed-form expression \revise{by, e.g., \cite{bendokat2024grassmann}},
\begin{equation}
	\Proj_{\upT_{P(t)}\Gr_k(\upT_{\vecx(t)}\calM)}(M):=\lrsquare{P(t),\lrsquare{P(t),M}},
	\label{eqn:Grassmann manifold projection onto tangent space}
\end{equation}
for any symmetric bilinear form $M$ on $\upT_{\vecx(t)}\calM$. The $\vecx$-part is the same as in Eq. \eqref{eqn:x-dynamics constrained existing} except that we adopt the orthogonal projector. The first ingredient of the $P$-part is responsible for tracking the lowest $k$-dimensional invariant subspace of $\hess{\calM}f(\vecx(t))$ in $\upT_{\vecx(t)}\calM$, following the Euler-Lagrange equation of the Rayleigh quotient minimization, 
$$\min_{\tilde P}~\trace(\tilde P\circ\hess{\calM}f(\vecx(t))),\quad\st~\tilde P\in\Gr_k(\upT_{\vecx(t)}\calM),$$
while the second ingredient accounts for the changes in the tangent spaces, as explained in \revise{Appendix \ref{appsec:preliminaries on Riemannian manifolds}}. The dynamic\revise{\sout{s}} \eqref{eqn:deterministic constrained saddle dynamics} can be horizontally lifted to \revise{a dynamic\revise{\sout{s}}} on $\St_k(\upT\calM)$ as shown in Lemma \ref{lem:decomposition of tangent space of Stiefel bundle}: 
\begin{equation}
	\begin{aligned}
		\frac{\dd\vecx}{\dd t}(t)&:=-R_{\vecx(t),V(t)}\big(\grad{\calM}f(\vecx(t))\big),\\
		\frac{\dd\vecv_i}{\dd t}(t)&:=-\Proj_{\vecx(t),V(t)^\perp}\big(\hess{\calM}f(\vecx(t))[\vecv_i(t)]\big)+\two_{\vecx(t)}\lrbracket{\frac{\dd\vecx}{\dd t}(t),\vecv_i(t)},\\
		&\qquad i=1,\ldots,k,\quad\text{with}\quad V(t):=(\vecv_1(t),\ldots,\vecv_k(t)).
	\end{aligned}
	\label{eqn:deterministic constrained saddle dynamics lifted}
\end{equation}

\par If the submanifolds \eqref{eqn:special Riemannian manifold} are considered, the second fundamental form in Eq. \eqref{eqn:deterministic constrained saddle dynamics lifted} recovers the third term of the right-hand side of Eq. \eqref{eqn:v-dynamics constrained existing}. However, thanks to the explicit characterization of the tangent space to $\Gr_k(\upT\calM)$ in Lemma \ref{lem:tangent space of Grassmann bundle}, our methodology \revise{does not rely on explicit global/local level-set representations of $\calM$} since we only require \revise{embedded-submanifold geometric primitives (e.g., $\upT_{\vecx}\calM$, $\Proj_{\upT_{\vecx}\calM}$, and $\two_{\vecx}$)}. 
Some examples are listed as follows.

\begin{example}[Second fundamental form in special cases]\mbox{}
	\begin{itemize} 
		\item Riemannian manifolds induced by global defining functions: $\calE=\R^n$, $\calM$ is defined in Eq. \eqref{eqn:special Riemannian manifold} \revise{\cite{li2015gentlest,liu2023constrained,yin2022constrained,zhang2012constrained,zhang2023discretization}}. For any $\vecx\in\calM$, we have
		\begin{align*}
			\upT_{\vecx}\calM&=\lrbrace{\vecv\in\R^n\mid\grad{}\vecc(\vecx)\vecv=0},\\
			\Proj_{\upT_{\vecx}\calM}(\vecv)&=\lrbracket{I_n-\grad{}\vecc(\vecx)^\top\big(\grad{}\vecc(\vecx)\cdot\grad{}\vecc(\vecx)^\top\big)^{-1}\grad{}\vecc(\vecx)}\vecv
		\end{align*}
		for any $\vecv\in\calE$. Therefore, by Eq. \eqref{eqn:differential of orthogonal projector}, for any $\vecu$, $\vecv\in\upT_{\vecx}\calM$, 
		\begin{align*}
			\two_{\vecx}(\vecu,\vecv)&=-\big(\hess{}\vecc(\vecx)[\vecu]\big)^\top\lrbracket{\grad{}\vecc(\vecx)\cdot\grad{}\vecc(\vecx)^\top}^{-1}\grad{}\vecc(\vecx)\vecv\\
			&\quad-\grad{}\vecc(\vecx)^\top\upD\lrbracket{\vecy\mapsto\big(\grad{}\vecc(\vecy)\cdot\grad{}\vecc(\vecy)^\top\big)^{-1}}(\vecx)[\vecu]\cdot\grad{}\vecc(\vecx)\vecv\\
			&\quad-\grad{}\vecc(\vecx)^\top\lrbracket{\grad{}\vecc(\vecx)\cdot\grad{}\vecc(\vecx)^\top}^{-1}\lrbracket{\hess{}\vecc(\vecx)[\vecu]}\vecv\\
			&=-\grad{}\vecc(\vecx)^\top\lrbracket{\grad{}\vecc(\vecx)\cdot\grad{}\vecc(\vecx)^\top}^{-1}\lrbracket{\hess{}\vecc(\vecx)[\vecu]}\vecv,
		\end{align*}
		where the last equality is due to $\grad{}\vecc(\vecx)\vecv=0$. Notably, this recovers the third term on the right-hand side of Eq. \eqref{eqn:v-dynamics constrained existing} if we let $\vecu=\frac{\dd\vecx}{\dd t}(t)$. Some special cases falling into this class and the corresponding second fundamental forms are listed as follows:
		\begin{itemize}
			\item Flat space: $\calM=\calE$, $\two_{\vecx}\equiv0$.
			
			\item Sphere: $\calE=\R^{d+1}$, $\calM=\mathbb{S}^d$, 
			$$\two_{\vecx}(\vecu,\vecv)=-\inner{\vecu,\vecv}\vecx,\quad\forall~\vecx\in\calM,~\vecu,\vecv\in\upT_{\vecx}\calM.$$
			
			\item Stiefel manifold: $\calE=\R^{n\times p}$, $\calM=\St_p(\R^n)$, 
			$$\two_X(U,V)=-\frac12X\big(U^\top V+V^\top U\big),\quad\forall~X\in\St_p(\R^n),~U,V\in\upT_X\calM.$$
		\end{itemize}
		
		\item Grassmann manifold: $\calE=\R_{\sym}^{n\times n}$, $\calM=\Gr_p(\R^n)$. For any $P\in\calM$, we have from \cite{bendokat2024grassmann} that
		$$\upT_P\calM=\lrbrace{\Gamma\in\R_{\sym}^{n\times n}\mid\Gamma P+P\Gamma=\Gamma},\quad\Proj_{\upT_P\calM}(\Gamma)=[P,[P,\Gamma]],\quad\forall~\Gamma\in\calE.$$
		Therefore, by Eq. \eqref{eqn:differential of orthogonal projector}, for any $\Gamma$, $\Delta\in\upT_P\calM$, 
		$$\two_P(\Delta,\Gamma)=[\Delta,[P,\Gamma]]+[P,[\Delta,\Gamma]]=[\Delta,[P,\Gamma]].$$
		The last equality is due to
		\begin{align*}
			[P,[\Delta,\Gamma]]&=[P,[\Delta P+P\Delta,\Gamma P+P\Gamma]]\\
			&=[P,\Delta P\Gamma+P\Delta\Gamma P-\Gamma P\Delta-P\Gamma\Delta P]\\
			&=P\Delta\Gamma P-P\Gamma\Delta P-P\Delta\Gamma P+P\Gamma\Delta P=0,
		\end{align*}
		where we have used $\Gamma P+P\Gamma=\Gamma$, $\Delta P+P\Delta=\Delta$, $P\Gamma P=0$, and $P\Delta P=0$.
		
		\item \revise{Manifold of fixed-rank matrices}: $\calE=\R^{m\times n}$, $\calM=\R_r^{m\times n}$. For any $X\in\calM$, let $X=U\Sigma V^\top$ be its singular value decomposition, where $U\in\St_r(\R^m)$, $\Sigma=\diag(\sigma_1,\ldots,\sigma_r)\in\R_{++}^{r\times r}$, and $V\in\St_r(\R^n)$. We have from \cite[Section 7.5]{boumal2023introduction} that 
		\begin{align*}
			&\upT_X\calM=\lrbrace{\left.\begin{pmatrix}
					U & U_\perp
				\end{pmatrix}\begin{pmatrix}
					A & B \\ C & 0
				\end{pmatrix}\begin{pmatrix}
					V^\top \\ V_\perp^\top
				\end{pmatrix}~\right|~A\in\R^{r\times r},~B\in\R^{r\times (n-r)},~C\in\R^{(m-r)\times r}},\\
			&\Proj_{\upT_X\calM}(Z)=Z-(I_m-UU^\top)Z(I_n-VV^\top),\quad\forall~Z\in\calE.
		\end{align*}
		Note that for any smooth curve $c:\R\supseteq\calI\to\calM$ satisfying $X(0)=X$ and $X'(0)=\Psi\in\upT_X\calM$, there exist smooth mappings $U(t)\in\St_r(\R^m)$, $\Sigma(t)\in\R^{r\times r}$, and $V(t)\in\St_r(\R^n)$ such that $X(t)=U(t)\Sigma(t)V(t)^\top$ and $U(0)=U$, $\Sigma(0)=\Sigma$, $V(0)=V$ hold in a small neighborhood of the origin. \revise{Note that $\Sigma(t)$ is invertible but need not be diagonal.} By the product rule, 
		$$U'(0)\Sigma V^\top +U\Sigma'(0)V^\top+U\Sigma V'(0)^\top=X'(0)=\Psi,$$
		which yields
		\begin{align*}
			\Sigma'(0)=U^\top\Psi V,~~U'(0)=(I_m-UU^\top)\Psi V\Sigma^{-1},~~V'(0)=(I_n-VV^\top)\Psi^\top U\Sigma^{-1}.
		\end{align*}
		Therefore, by Eq. \eqref{eqn:differential of orthogonal projector}, for any $\Phi$, $\Psi\in\upT_X\calM$, 
		\begin{align*}
			\two_X(\Psi,\Phi)&=\lrbracket{(I_m-UU^\top)\Psi V\Sigma^{-1}U^\top+U\Sigma^{-1}V^\top\Psi^\top(I_m-UU^\top)}\Phi(I_n-VV^\top)\\
			&\quad+(I_m-UU^\top)\Phi\lrbracket{(I_n-VV^\top)\Psi^\top U\Sigma^{-1}V^\top+V\Sigma^{-1}U^\top\Psi(I_n-VV^\top)}\\
			&=(I_m-UU^\top)\lrbracket{\Psi V\Sigma^{-1}U^\top\Phi+\Phi V\Sigma^{-1}U^\top\Psi}(I_n-VV^\top),
		\end{align*}
		where the last equality is due to the fact that $(I_m-UU^\top)\Phi(I_n-VV^\top)=0$. 
	\end{itemize}
	
\end{example}

\begin{remark}[\revise{Relation between different dynamics}]
    \par We shall point out that\revise{, when $\calM$ is defined by Eq. \eqref{eqn:special Riemannian manifold},} the CSD \eqref{eqn:deterministic constrained saddle dynamics lifted} \revise{differs from the existing dynamics in Eqs. \eqref{eqn:x-dynamics constrained existing} and \eqref{eqn:v-dynamics constrained existing} in off-diagonal tangent terms of $V$-dynamic\revise{\sout{s}} when $k>1$:} 
in the existing dynamics \eqref{eqn:v-dynamics constrained existing}, the component is
$$-\Proj_{\vecx(t),\vecv_i(t)^\perp}\big(\hess{\calM}f(\vecx(t))[\vecv_i(t)]\big)+2\sum_{j=1}^{i-1}\Proj_{\vecx(t),\vecv_j(t)}\big(\hess{\calM}f(\vecx(t))[\vecv_i(t)]\big),$$
while in ours, the component reads
$$-\Proj_{\vecx(t),\vecv_i(t)^\perp}\big(\hess{\calM}f(\vecx(t))[\vecv_i(t)]\big)+\sum_{j\ne i}\Proj_{\vecx(t),\vecv_j(t)}\big(\hess{\calM}f(\vecx(t))[\vecv_i(t)]\big).$$
The existing one is derived from the Lagrangian formalism combined with operator splitting \cite{yin2022constrained}, which is favorable due to its decoupled nature. 
\revise{Nonetheless, these two dynamics essentially correspond to the same trajectory in the orthogonal projector space (namely, dynamic\revise{\sout{s}} \eqref{eqn:deterministic constrained saddle dynamics}). This can be verified by the quotient mapping associated with the quotient structure \eqref{eqn:Grassmann bundle quotient}. With the introduction of Grassmann bundle, we will be able to remove unnecessary nondegeneracy assumptions on the spectrum of Riemannian Hessian in theoretical analysis; see Section \ref{sec:theoretical analysis}.} 
\end{remark}

\par It is straightforward to verify that Eq. \eqref{eqn:requirements on dx and dv} holds if $(\vecx(t),P(t))\in\Gr_k(\upT\calM)$. Indeed, we are able to show the global well-definedness of the dynamic\revise{\sout{s}} \eqref{eqn:deterministic constrained saddle dynamics} if $\calM$ is compact.

\begin{theorem}[Global well-definedness of the dynamic\revise{\sout{s}} \eqref{eqn:deterministic constrained saddle dynamics}]\label{thm:global well-definedness}
	Suppose that $\calM$ is a compact \revise{embedded Riemannian submanifold of $\R^n$} with $\dim(\calM)=d$ \revise{and let $1\le k<d$}. If $(\vecx(0),P(0))\in\Gr_k(\upT\calM)$, then the trajectory generated by the dynamic\revise{\sout{s}} \eqref{eqn:deterministic constrained saddle dynamics} always stays on $\Gr_k(\upT\calM)$. 
\end{theorem}
\begin{proof}
	Since $\calM$ is compact in $\calE$, so is $\Gr_k(\upT\calM)$ in its ambient space. \revise{Indeed, 
    $$\Gr_k(\upT\calM)=\lrbrace{(\vecx,P)\mid\vecx\in\calM,~P\in\Gr_k(\R^n),~\Proj_{\upT_{\vecx}\calM}\circ P\circ\Proj_{\upT_{\vecx}\calM}=P},$$
    which is a bounded subset in $\calM\times\R_{\sym}^{n\times n}$. Its closedness is due to the fact that the map $(\vecx,P)\mapsto\Proj_{\upT_{\vecx}\calM}\circ P\circ\Proj_{\upT_{\vecx}\calM}-P$ is continuous.} The conclusion then follows directly from \cite[Corollary 9.17]{lee2012introduction}.
\end{proof}

\par We then discretize the dynamic\revise{\sout{s}} \eqref{eqn:deterministic constrained saddle dynamics} by directly applying the forward Euler scheme and using the retraction \eqref{eqn:retraction over Grassmann bundle}, which yields Algorithm \ref{alg:discretized deterministic constrained saddle dynamics}.
\begin{algorithm}[!t]
	\caption{Discretized CSD on $\Gr_k(\upT\calM)$.}
	\label{alg:discretized deterministic constrained saddle dynamics}
	\begin{algorithmic}[1]
		\Require{Initial point $(\vecx^{(0)},P^{(0)})\in\Gr_k(\upT\calM)$, step size $\eta>0$, maximum iteration number $maxit\in\N$, convergence tolerance $tol>0$.}
		\State{Set $t:=0$, $r_{\vecx}^{(0)}=r_P^{(0)}:=\infty$.}
		\While{$t<maxit$ \revise{and} $\max\{r_{\vecx}^{(t)},r_P^{(t)}\}>tol$}
		\State{Compute the $\vecx$-direction: 
			$$\vecd_{\vecx}^{(t)}:=-R_{\vecx^{(t)},P^{(t)}}\big(\grad{\calM}f(\vecx^{(t)})\big)\in\upT_{\vecx^{(t)}}\calM.$$}
		\State{Compute the $P$-direction:
			$$\vecd_P^{(t)}:=-\Proj_{\upT_{P^{(t)}}\Gr_k(\upT_{\vecx^{(t)}}\calM)}\circ\hess{\calM}f(\vecx^{(t)})\in\upT_{P^{(t)}}\Gr_k(\upT_{\vecx^{(t)}}\calM).$$
		}
		\State{Update $\vecx$ and $P$ using the retraction defined in Eq. \eqref{eqn:retraction over Grassmann bundle}:
			$$(\vecx^{(t+1)},P^{(t+1)}):=\upR^{\Gr_k(\upT\calM)}_{(\vecx^{(t)},P^{(t)})}\Big(\eta\vecd_{\vecx}^{(t)},\revise{\eta\big(\vecd_P^{(t)}+\widehat\two_{\vecx^{(t)}}(\vecd_{\vecx}^{(t)},P^{(t)})\big)}\Big)\in\Gr_k(\upT\calM).$$
		}\label{step:update v direction}
		\If{$t=0$}
		\State{Set $r_{\vecx}^{(t+1)}=r_P^{(t+1)}:=1$.}
		\Else
		\State{Update $r_{\vecx}^{(t+1)}:=\snorm{\vecd_{\vecx}^{(t)}}/\snorm{\vecd_{\vecx}^{(0)}}$ and $r_P^{(t+1)}:=\snorm{\vecd_P^{(t)}}/\snorm{\vecd_P^{(0)}}$.}
		\EndIf
		\State{Set $t:=t+1$.}
		\EndWhile
		\Ensure{$(\vecx^{(t)},P^{(t)})\in\Gr_k(\upT\calM)$.}
	\end{algorithmic}
\end{algorithm}
Note that in step \ref{step:update v direction}, we do not \revise{need to} compute the term involving the second fundamental form. This is credited to the construction of the retraction \eqref{eqn:retraction over Grassmann bundle}, in which the normal component makes no difference. In light of the quotient structure \eqref{eqn:Grassmann bundle quotient} \revise{and Lemmas \ref{lem:decomposition of tangent space of Stiefel bundle} and \ref{lem:general retraction on Grassmann bundle}}, we \revise{have that the horizontal lift of $\vecd_P^{(t)}$ in step \ref{step:update v direction} is exactly $\vecd_V^{(t)}:=(\vecd_{\vecv_1}^{(t)},\ldots,\vecd_{\vecv_k}^{(t)})$, where
$$\vecd_{\vecv_i}^{(t)}:=-\Proj_{\vecx^{(t)},(V^{(t)})^\perp}\big(\hess{\calM}f(\vecx^{(t)})[\vecv_i^{(t)}]\big),~~i=1,\ldots,k,$$
and thus} obtain a representative version running on $\St_k(\upT\calM)$; see Algorithm \ref{alg:discretized deterministic constrained saddle dynamics lifted}. In comparison with Algorithm \ref{alg:discretized deterministic constrained saddle dynamics}, Algorithm \ref{alg:discretized deterministic constrained saddle dynamics lifted} can be computationally \revise{cheaper}. 
\begin{algorithm}[!t]
	\caption{Discretized CSD using representatives on $\St_k(\upT\calM)$.}
	\label{alg:discretized deterministic constrained saddle dynamics lifted}
	\begin{algorithmic}[1]
		\Require{Initial point $(\vecx^{(0)},V^{(0)})\in\St_k(\upT\calM)$, step size $\eta>0$, maximum iteration number $maxit\in\N$, convergence tolerance $tol>0$.}
		\State{Set $t:=0$, $r_{\vecx}^{(0)}=r_V^{(0)}:=\infty$.}
		\While{$t<maxit$ \revise{and} $\max\{r_{\vecx}^{(t)},r_V^{(t)}\}>tol$}
		\State{Compute the $\vecx$-direction: 
			$$\vecd_{\vecx}^{(t)}:=-R_{\vecx^{(t)},V^{(t)}}\big(\grad{\calM}f(\vecx^{(t)})\big)\in\upT_{\vecx^{(t)}}\calM.$$}
		\State{Compute the $V$-direction:
			$$\vecd_{\vecv_i}^{(t)}:=-\Proj_{\vecx^{(t)},(V^{(t)})^\perp}\big(\hess{\calM}f(\vecx^{(t)})[\vecv_i^{(t)}]\big)\in\upT_{\vecx^{(t)}}\calM,~~i=1,\ldots,k,$$
			and set $\vecd_V^{(t)}:=(\vecd_{\vecv_1}^{(t)},\ldots,\vecd_{\vecv_k}^{(t)})\in\upT_{V^{(t)}}\St_k(\upT_{\vecx^{(t)}}\calM).$
		}
		\State{Update $\vecx$ and $V$ using the retraction defined in Eq. \eqref{eqn:retraction over Stiefel bundle}:
			$$(\vecx^{(t+1)},V^{(t+1)}):=\upR^{\Gr_k(\upT\calM)}_{(\vecx^{(t)},V^{(t)})}\Big(\eta\vecd_{\vecx}^{(t)},\revise{\eta\big(\vecd_{\vecv_i}^{(t)}+\two_{\vecx^{(t)}}(\vecd_{\vecx}^{(t)},\vecv_i^{(t)})\big)_{i=1}^k}\Big)\in\St_k(\upT\calM).$$
		}
		\If{$t=0$}
		\State{Set $r_{\vecx}^{(t+1)}=r_V^{(t+1)}:=1$.}
		\Else
		\State{Update $r_{\vecx}^{(t+1)}:=\snorm{\vecd_{\vecx}^{(t)}}/\snorm{\vecd_{\vecx}^{(0)}}$ and $r_V^{(t+1)}:=\snorm{\vecd_V^{(t)}}/\snorm{\vecd_V^{(0)}}$.}
		\EndIf
		\State{Set $t:=t+1$.}
		\EndWhile
		\Ensure{$(\vecx^{(t)},V^{(t)})\in\St_k(\upT\calM)$.}
	\end{algorithmic}
\end{algorithm}

\section{Theoretical analysis}\label{sec:theoretical analysis}

In this section, we analyze the linear stability of the dynamic\revise{\sout{s}} \eqref{eqn:deterministic constrained saddle dynamics} and the local convergence of Algorithm \ref{alg:discretized deterministic constrained saddle dynamics}. 

\subsection{Linear stability analysis of the dynamics}

\par Let $\vecw:=(\vecx,P)\in\Gr_k(\upT\calM)$ and $\vech:\Gr_k(\upT\calM)\to\upT(\Gr_k(\upT\calM))$ be the vector field defined as $\vech:=(h_1,h_2)$, where
\begin{equation}
	\begin{aligned}
		h_1(\vecw)&:=-R_{\vecw}\big(\grad{\calM}f(\vecx)\big),\\
		h_2(\vecw)&:=-\Proj_{\upT_{P}\Gr_k(\upT_{\vecx}\calM)}\circ\hess{\calM}f(\vecx)+\widehat\two_{\vecx}(h_1(\vecw),P).
	\end{aligned}
	\label{eqn:solution map in dynamics}
\end{equation}
Note that $R_{\vecw}$ should be identified as the one in Eq. \eqref{eqn:reflection operator in CSD}. \revise{It} follows that Eq. \eqref{eqn:deterministic constrained saddle dynamics} can be rewritten compactly as the following dynamic\revise{\sout{s}} on $\Gr_k(\upT\calM)$:
\begin{equation}
	\frac{\dd\vecw}{\dd t}(t)=\vech(\vecw(t)),
	\label{eqn:deterministic CSD compact}
\end{equation}
with the initial condition $\vecw(0)=(\vecx(0),P(0))\in\Gr_k(\upT\calM)$. 

\par In the following, we aim to establish the linear stability of the dynamics \eqref{eqn:deterministic CSD compact} using the geometrical tools for $\Gr_k(\upT\calM)$ in Section \ref{subsec:geometries of Stiefel and Grassmann bundles}. Before that, we first figure out the expression of the differential $\upD\vech(\vecw)$. 

\begin{lemma}[Differential of $\vech$]\label{lem:differential of solution map}
	Suppose that $f$ is $C^3$, $\calM$ is \revise{an embedded} Riemannian submanifold of a Euclidean space $\calE$ with $\dim(\calM)=d$, \revise{$1\le k<d$}, and $\vech:\Gr_k(\upT\calM)\to\upT(\Gr_k(\upT\calM))$ is the vector field on $\Gr_k(\upT\calM)$ defined in Eq. \eqref{eqn:solution map in dynamics}. Then for any $\vecw:=(\vecx,P)\in\Gr_k(\upT\calM)$ and $D:=(\bdelta,\Delta)\in\upT_{\vecw}(\Gr_k(\upT\calM))$, it holds that 
	\begin{align*}
		\upD h_1(\vecw)[D]&=-R_{\vecw}\big(\hess{\calM}f(\vecx)[\bdelta]\big)+\vecm_1(\vecw)[\bdelta]+\vecm_2(\vecw)[\Delta],\\
		\upD h_2(\vecw)[D]&=-[\Delta,[P,\hess{\calM}f(\vecx)]]-[P,[\Delta,\hess{\calM}f(\vecx)]]+\vecm_3(\vecw)[\bdelta]+\vecm_4(\vecw)[\Delta],
	\end{align*}
	where \revise{$\vecm_1(\vecw)\in\calL(\upT_{\vecx}\calM,\calE)$, $\vecm_2(\vecw)\in\calL(\calS(\calE),\calE)$, $\vecm_3(\vecw)\in\calL(\upT_{\vecx}\calM,\calS(\calE))$, $\vecm_4(\vecw)\in\calL(\calS(\calE),\calS(\calE))$ are defined as}
	\begin{align*}
		&\vecm_1(\vecw)[\bdelta]:=-\upD(\vecy\mapsto\revise{{\rm Op}}^R_{\vecy,P})(\vecx)[\bdelta]\big(\grad{\calM}f(\vecx)\big)-\revise{{\rm Op}}^R_{\vecw}\big(\two_{\vecx}(\bdelta,\grad{\calM}f(\vecx))\big),\\
		&\vecm_2(\vecw)[\Delta]:=-\upD(Q\mapsto\revise{\rm Op}^R_{\vecx,Q})(P)[\Delta]\big(\grad{\calM}f(\vecx)\big),\\
		&\vecm_3(\vecw)[\bdelta]:=-\upD(\vecy\mapsto \revise{\rm Op}^{\Proj}_{\vecy,P})(\vecx)[\bdelta]\circ\hess{\calM}f(\vecx)-\revise{\rm Op}^{\Proj}_{\vecw}\circ\upD\bar H(\vecx)[\bdelta]\\
		&\qquad\qquad\qquad+\upD(\vecy\mapsto \revise{\rm Op}^{\widehat\two}_{\vecy})(\vecx)[\bdelta](h_1(\vecw),P)+\revise{\rm Op}^{\widehat\two}_{\vecx}\big(\upD(\vecy\mapsto h_1(\vecy,P))(\vecx)[\bdelta],P\big),\\
		&\vecm_4(\vecw)[\Delta]:=\revise{\rm Op}^{\widehat\two}_{\vecx}\big(\upD(Q\mapsto h_1(\vecx,Q))(P)[\Delta],P\big)+\revise{\rm Op}^{\widehat\two}_{\vecx}(h_1(\vecw),\Delta),
	\end{align*}
	\revise{$(\vecx,P)\mapsto\revise{\rm Op}^{R}_{\vecx,P}$ and $(\vecx,P)\mapsto\revise{\rm Op}^{\Proj}_{\vecx,P}$ are respectively smooth extensions of $(\vecx,P)\mapsto R_{\vecx,P}$ and $(\vecx,P)\mapsto\Proj_{\upT_P\Gr_k(\upT_{\vecx}\calM)}$ in $\calE\times\calS(\calE)$, $\vecx\mapsto\revise{\rm Op}^{\widehat\two}_{\vecx}$ and $\bar H$ are respectively smooth extensions of $\vecx\mapsto\widehat\two_{\vecx}$ and $\hess{\calM}f$ in $\calE$}. 
\end{lemma}
\begin{proof}
	Let $c:\R\supseteq\calI\to\Gr_k(\upT\calM)$ be a smooth curve on $\Gr_k(\upT\calM)$ such that $c(0)=\vecw$ and $c'(0)=D$. We thus have \revise{from Lemma \ref{lem:bundle structures as embedded submanifolds}}
	$$\upD\vech(\vecw)[D]=(\vech\circ c)'(0)=\big((h_1\circ c)'(0),~(h_2\circ c)'(0)\big).$$
	By the product and chain rules,
	\begin{align*}
		&~(h_1\circ c)'(0)=\lim_{s\to0}\lrbracket{h_1(c(s))-h_1(c(0))}/s\nonumber\\
		=&~-\big(\upD(\vecy\mapsto\revise{\rm Op}^R_{\vecy,P})(\vecx)[\bdelta]+\upD(Q\mapsto\revise{\rm Op}^R_{\vecx,Q})(P)[\Delta]\big)\big(\grad{\calM}f(\vecx)\big)\nonumber\\
		&~-\revise{\rm Op}^R_{\vecw}\lrbracket{\upP_{\vecx,\bdelta}\big(\grad{}f(\vecx)\big)+\Proj_{\upT_{\vecx}\calM}\big(\hess{}f(\vecx)[\bdelta]\big)}\nonumber\\
		=&~-\big(\upD(\vecy\mapsto\revise{\rm Op}^R_{\vecy,P})(\vecx)[\bdelta]+\upD(Q\mapsto\revise{\rm Op}^R_{\vecx,Q})(P)[\Delta]\big)\big(\grad{\calM}f(\vecx)\big)\nonumber\\
		&~-\revise{\rm Op}^R_{\vecw}\lrbracket{\hess{\calM}f(\vecx)[\bdelta]+\two_{\vecx}\big(\bdelta,\grad{\calM}f(\vecx)\big)}\nonumber\\
		=&~-R_{\vecw}\big(\hess{\calM}f(\vecx)[\bdelta]\big)+\vecm_1(\vecw)[\bdelta]+\vecm_2(\vecw)[\Delta],
	\end{align*}
	where the second equality uses Eqs. \eqref{eqn:definition of Riemannian Hessian} and \eqref{eqn:exchange property of projection differential}, and
	\begin{align*}
		&~(h_2\circ c)'(0)=\lim_{s\to0}\lrbracket{h_2(c(s))-h_2(c(0))}/s\\
		=&~-\upD(\vecy\mapsto\revise{\rm Op}^{\Proj}_{\vecy,P})(\vecx)[\bdelta]\circ\hess{\calM}f(\vecx)-\revise{\rm Op}^{\Proj}_{\vecw}\circ\upD\bar H(\vecx)[\bdelta]\\
		&~+\upD(\vecy\mapsto\revise{\rm Op}^{\widehat\two}_{\vecy})(\vecx)[\bdelta](h_1(\vecw),P)+\revise{\rm Op}^{\widehat\two}_{\vecx}\big(\upD(\vecy\mapsto h_1(\vecy,P))(\vecx)[\bdelta],P\big)\\
		&~-\upD(Q\mapsto\revise{\rm Op}^{\Proj}_{\vecx,Q})(P)[\Delta]\circ\hess{\calM}f(\vecx)+\revise{\rm Op}^{\widehat\two}_{\vecx}(h_1(\vecw),\Delta)\\
		&~+\revise{\rm Op}^{\widehat\two}_{\vecx}\big(\upD(Q\mapsto h_1(\vecx,Q))(P)[\Delta],P\big)\\
		=&~-[\Delta,[P,\hess{\calM}f(\vecx)]]-[P,[\Delta,\hess{\calM}f(\vecx)]]+\vecm_3(\vecw)[\bdelta]+\vecm_4(\vecw)[\Delta],
	\end{align*}
	where the last equality follows from Eq. \eqref{eqn:Grassmann manifold projection onto tangent space}. The proof is complete. 
\end{proof}

\begin{theorem}[Linear stability of the dynamics \eqref{eqn:deterministic CSD compact}]\label{thm:linear stability of CSD}
	Suppose that $f$ is $C^3$, $\calM$ is \revise{an embedded} Riemannian submanifold of a Euclidean space $\calE$ with $\dim(\calM)=d$, \revise{$1\le k<d$}, $\vecw^\star:=(\vecx^\star,P^\star)\in\Gr_k(\upT\calM)$, and $\hess{\calM}f(\vecx^\star)$ is nondegenerate, whose eigenvalues are $\lambda_1^\star\le\cdots\le\lambda_k^\star<\lambda_{k+1}^\star\le\cdots\le\lambda_d^\star$. Then $\vecw^\star$ is a linearly steady state of the dynamic\revise{\sout{s}} \eqref{eqn:deterministic CSD compact} if and only if $\vecx^\star$ is an index-$k$ constrained SP of $f$ on $\calM$ and $P^\star$ is an orthogonal projector onto the lowest $k$-dimensional invariant subspace of $\hess{\calM}f(\vecx^\star)$. 
\end{theorem}
\begin{proof}
	The proof revolves around the spectrum of $\upD\vech(\vecw^\star)$. 
	
	\medskip
	
	\par\textbf{Necessity.} Suppose that $\vecw^\star$ is a linearly steady state of dynamic\revise{\sout{s}} \eqref{eqn:deterministic CSD compact}. Therefore, $\dd\vecw/\dd t$ vanishes at $\vecw^\star$, which implies that 
	\begin{itemize}
		\item $\grad{\calM}f(\vecx^\star)=0$, i.e., $\vecx^\star$ is a critical point of $f$ on $\calM$. This further yields $\vecm_1(\vecw^\star)=0$, $\vecm_2(\vecw^\star)=0$, and $\vecm_4(\vecw^\star)=0$ from their definitions in Lemma \ref{lem:differential of solution map};
		
		\item $[P^\star,[P^\star,\hess{\calM}f(\vecx^\star)]]=0$, i.e., $P^\star$ is an orthogonal projector onto a $k$-dimensional invariant subspace of $\hess{\calM}f(\vecx^\star)$. 
	\end{itemize}
	The differential $\upD\vech(\vecw^\star)[D]$ in Lemma \ref{lem:differential of solution map} thus simplifies to
	\begin{equation}
		\upD\vech(\vecw^\star)[D]=\begin{pmatrix}
			-R_{\vecw^\star}\big(\hess{\calM}f(\vecx^\star)[\bdelta]\big)\\[0.1cm]
			-[P^\star,[\Delta,\hess{\calM}f(\vecx^\star)]]+\vecm_3(\vecw^\star)[\bdelta]
		\end{pmatrix}.
		\label{eqn:differential of solution map}
	\end{equation}
	\revise{Since $\vecw^\star$ is a steady state, $\upD\vech(\vecw^\star)[D]$ coincides with the covariant derivative of $\vech$ at $\vecw^\star$ along $D$ and thus belongs to $\upT_{\vecw^\star}\Gr_k(\upT\calM)$.} 
	
	\par By the \revise{Hurwitz} linearization theorem \cite{arrowsmith1992dynamical}, the linear stability of $\vecw^\star$ implies that all of the eigenvalues of $\upD\vech(\vecw^\star)$ have negative real parts. By Lemma \ref{lem:tangent space of Grassmann bundle}, these eigenvalues coincide with \revise{the generalized eigenvalues of the matrix pair $(A,B)$,} where $A\in\R^{\hat d\times\hat d}$ and $B\in\R_{\sym}^{\hat d\times\hat d}$ (with $\hat d:=d+k(d-k)$) are \revise{respectively the representation of $\upD\vech(\vecw^\star)$ under the basis 
    \begin{multline*}
        \lrbrace{D_q^\star:=(\bdelta_q^\star,\Delta_q^\star):=\lrbracket{\vecv_q^\star,\widehat\two_{\vecx^\star}(\vecv_q^\star,P^\star)}}_{1\le q\le d}\\
        \bigcup\lrbrace{D_{ij}^\star:=(\bdelta_{ij}^\star,\Delta_{ij}^\star):=\lrbracket{0,\frac{1}{\sqrt{2}}\big(\vecv_i^\star(\vecv_j^\star)^\top+\vecv_j^\star(\vecv_i^\star)^\top\big)}}_{1\le i\le k\atop k+1\le j\le d}
    \end{multline*}
    and the Gram matrix} using the Sasaki-type metric in Definition \ref{def:sasaki metric}: 
	\begin{align*}
		&A_{q,q'}:=\inner{D_{q'}^\star,\upD\vech(\vecw^\star)[D_q^\star]}_{\vecw^\star},& &B_{q,q'}:=\inner{D_{q'}^\star,D_q^\star}_{\vecw^\star},\\
        &\revise{A_{q,ij}:=\inner{D_{ij}^\star,\upD\vech(\vecw^\star)[D_q^\star]}_{\vecw^\star},}& &\revise{B_{q,ij}:=\inner{D_{ij}^\star,D_q^\star}_{\vecw^\star},}\\
        &\revise{A_{ij,q}:=\inner{D_q^\star,\upD\vech(\vecw^\star)[D_{ij}^\star]}_{\vecw^\star},}& &\revise{B_{ij,q}:=\inner{D_q^\star,D_{ij}^\star}_{\vecw^\star},}\\
		&A_{ij,i'j'}:=\inner{D_{i'j'}^\star,\upD\vech(\vecw^\star)[D_{ij}^\star]}_{\vecw^\star},& &B_{ij,i'j'}:=\inner{D_{i'j'}^\star,D_{ij}^\star}_{\vecw^\star}.
	\end{align*}
	Here, $\{(\mu_q^\star,\vecv_q^\star)\}_{q=1}^d$ are the eigenpairs of $\hess{\calM}f(\vecx^\star)$ such that $P^\star=\sum_{\ell=1}^k\vecv_\ell^\star(\vecv_\ell^\star)^\top$. 
    \revise{W}e could work out the entries in $A$ and $B$ explicitly:
	\begin{align*}
		A_{q,q'}&=\left\{\begin{array}{ll}
			\mu_q^\star\delta_{q,q'}, & q,q'=1,\ldots,k, \\
			-\mu_q^\star\delta_{q,q'}, & q,q'=k+1,\ldots,d,\\
			0, & \text{otherwise},
		\end{array}\right.\\
		\revise{A_{ij,q}}&=0,~q=1,\ldots,d,~i=1,\ldots,k,~j=k+1,\ldots,d,\\
		\revise{A_{q,ij}}&=\inner{\Delta_{ij}^\star,\vecm_3(\vecw^\star)[\vecv_q^\star]},~q=1,\ldots,d,~i=1,\ldots,k,~j=k+1,\ldots,d,\\
		A_{ij,i'j'}&=(\mu_i^\star-\mu_j^\star)\delta_{i,i'}\delta_{j,j'},~i,i'=1,\ldots,k,~j,j'=k+1,\ldots,d,
	\end{align*}
	and 
    \revise{\begin{align*}
        B_{q,q'}&=\delta_{q,q'},~q,q'=1,\ldots,d,\\
        B_{q,ij}&=B_{ij,q}=0,~q=1,\ldots,d,~i=1,\ldots,k,~j=k+1,\ldots,d,\\
        B_{ij,i'j'}&=\delta_{i,i'}\delta_{j,j'},~i,i'=1,\ldots,k,~j,j'=k+1,\ldots,d.
    \end{align*}}
	Consequently, \revise{$B$ is an identity matrix, $A$ is lower (or upper) triangular, the generalized eigenvalues of $(A,B)$ are identical to the eigenvalues of $A$, which} are exactly
	$$\{\mu_q^\star\}_{q=1}^k,~~\{-\mu_q^\star\}_{q=k+1}^d,~~\text{and}~~\{\mu_i^\star-\mu_j^\star\}_{i=1,\ldots,k\atop j=k+1,\ldots,d}.$$
	Since they are all negative, we get $\max\{\mu_1^\star,\ldots,\mu_k^\star\}<0<\min\{\mu_{k+1}^\star,\ldots,\mu_d^\star\}$. This indicates that $\vecx^\star$ is an index-$k$ constrained SP of $f$ on $\calM$, $P^\star$ is an orthogonal projector onto the lowest $k$-dimensional invariant subspace of $\hess{\calM}f(\vecx^\star)$, and $\{\mu_1^\star,\ldots,\mu_k^\star\}=\{\lambda_1^\star,\ldots,\lambda_k^\star\}$, as desired.
	
	\medskip
    
	\par\textbf{Sufficiency.} Suppose that $\vecx^\star$ is an index-$k$ constrained SP of $f$ on $\calM$ and $P^\star$ is an orthogonal projector onto the lowest $k$-dimensional invariant subspace of $\hess{\calM}f(\vecx^\star)$. Therefore, $\grad{\calM}f(\vecx^\star)=0$ and $\lambda_1^\star\le\cdots\le\lambda_k^\star<0<\lambda_{k+1}^\star\le\cdots\lambda_d^\star$. By similar arguments as above, we could obtain that the eigenvalues of $\upD\vech(\vecw^\star)$ are exactly
	$$\{\lambda_q^\star\}_{q=1}^k,~~\{-\lambda_q^\star\}_{q=k+1}^d,~~\text{and}~~\{\lambda_i^\star-\lambda_j^\star\}_{i=1,\ldots,k\atop j=k+1,\ldots,d},$$
	which are all negative, implying that $\vecw^\star=(\vecx^\star,P^\star)$ is linearly steady. 
\end{proof}

\begin{remark}\label{rem:weaker assumptions}
	In comparison with the linear stability results in the existing work \cite[Theorem 1]{yin2022constrained}, Theorem \ref{thm:linear stability of CSD} not only \revise{avoids the use of explicit global/local level-set representations of $\calM$}, 
    but also requires weaker assumptions on the eigenvalues of Riemannian Hessian. More specifically, we only ask for a positive gap between $\lambda_k^\star$ and $\lambda_{k+1}^\star$, whereas the existing one assumes additionally that $\{\lambda_i^\star\}_{i=1}^k$ are all simple (cf. \cite[Eqs. (34) and (40)]{yin2022constrained}). This advantage is due to the fact that we respect the intrinsic quotient structure \eqref{eqn:Grassmann bundle quotient} by introducing the Grassmann bundle $\Gr_k(\upT\calM)$.
\end{remark}

\subsection{Local convergence analysis of the discretized algorithm}

\par In what follows, we show the local convergence properties of Algorithm \ref{alg:discretized deterministic constrained saddle dynamics}. \revise{We first recall an auxiliary lemma, whose proof can be largely covered by the literature.}

\begin{lemma}\label{lem:diffeomorphism nbd}
	Suppose that $\calN$ is \revise{an embedded} Riemannian submanifold of a Euclidean space, $\upR^{\calN}:\upT\calN\to\calN$ is a retraction on $\calN$, and $\vecy\in\calN$. Then there exist constants $c_1>0$ and $r_1>0$ such that the following two statements hold at the same time:
	\begin{itemize}
		\item $\Exp_{\vecy}^{\calN}:\upT_{\vecy}\calN\to\calN$ is a diffeomorphism between $\calB_{\vecy}(0,r_1):=\{\vecu\in\upT_{\vecy}\calN\mid\snorm{\vecu}<r_1\}$ and $\calU:=\Exp_{\vecy}^{\calN}(\calB_{\vecy}(0,r_1))$.
		
		\item The inequality $\dist_{\calN}(\vecy',\upR^{\calN}_{\vecy'}(\vecu'))\le c_1\snorm{\vecu'}$ holds for any $(\vecy',\vecu')$ such that $\vecy'\in\Exp_{\vecy}(\mathrm{cl}(\calB_{\vecy}(0,r_1)))$ and $\vecu'\in\upT_{\vecy'}\calN$ with $\snorm{\vecu'}\le r_1$. In particular, if $\upR^{\calN}=\Exp^{\calN}$, then the equality holds with $c_1=1$. 
	\end{itemize}
\end{lemma}
\begin{proof}
	For the first statement, it suffices to note $\upD\Exp_{\vecy}^{\calN}(0)=\Id_{\upT_{\vecy}\calN}$ and then use the inverse function theorem. For the second one, please refer to \cite[Lemma 6.32 and Proposition 10.22]{boumal2023introduction}.
\end{proof}

\begin{lemma}[Residual after single iteration]\label{lem:residual mapping}
	Suppose that $f$ is $C^3$, $\calM$ is \revise{an embedded} Riemannian submanifold of a Euclidean space $\calE$ with $\dim(\calM)=d$, \revise{$1\le k<d$}, $\vecx^\star\in\calM$ is an index-$k$ constrained SP of $f$ on $\calM$, $\hess{\calM}f(\vecx^\star)$ is nondegenerate, whose eigenvalues are $\lambda_1^\star\le\cdots\le\lambda_k^\star<\lambda_{k+1}^\star\le\cdots\le\lambda_d^\star$, and $P^\star\in\Gr_k(\upT_{\vecx^\star}\calM)$ is an orthogonal projector onto the lowest $k$-dimensional invariant subspace of $\hess{\calM}f(\vecx^\star)$. Let $\vecw^\star:=(\vecx^\star,P^\star)\in\Gr_k(\upT\calM)$. Consider the following iterative mapping $F:\upT_{\vecw^\star}(\Gr_k(\upT\calM))\to\upT_{\vecw^\star}(\Gr_k(\upT\calM))$ defined as
	\begin{equation}
		F(D):=\lrbracket{\Exp_{\vecw^\star}^{\Gr_k(\upT\calM)}}^{-1}\lrbracket{\upR^{\Gr_k(\upT\calM)}_{\tilde\vecw(D)}\big(\eta\cdot\vech(\tilde\vecw(D))\big)}=:(F_1(D),F_2(D)),
		\label{eqn:iterative mapping}
	\end{equation}
	for any $D:=(\bdelta,\Delta)\in\upT_{\vecw^\star}(\Gr_k(\upT\calM))$, where $\tilde\vecw(D):=\Exp_{\vecw^\star}^{\Gr_k(\upT\calM)}(D)\in\Gr_k(\upT\calM)$, and $\eta$ is a step size satisfying
	\begin{equation}
		0<\eta<\frac{\min\{(r_1-r_2)/c_1,r_1\}}{\sup_{\snorm{D}_{\vecw^\star}\le r_1}\snorm{\vech(\tilde\vecw(D))}_{\revise{\tilde\vecw(D)}}},
		\label{eqn:step size ub}
	\end{equation}
	with $c_1$ and $r_1$ the constants in Lemma \ref{lem:diffeomorphism nbd} (associated to $\calN=\Gr_k(\upT\calM)$, $\upR^{\calN}=\upR^{\Gr_k(\upT\calM)}$ or $\Exp^{\Gr_k(\upT\calM)}$, and $\vecy=\vecw^\star$) and $r_2\in(0,r_1)$. Then $F$ is well-defined and smooth over the subset $\{D\in\upT_{\vecw^\star}(\Gr_k(\upT\calM)):\snorm{D}_{\vecw^\star}<r_2\}$. Moreover, there exist constants $r_3\in(0,r_2]$ and $c_3\ge0$ such that
	\begin{equation}
		\begin{pmatrix}
			\snorm{F_1(D)}\\
			\snorm{F_2(D)_{\vecx^\star}}
		\end{pmatrix}\le A(\eta)\begin{pmatrix}
			\snorm{\bdelta}\\
			\snorm{\Delta_{\vecx^\star}}
		\end{pmatrix}
		\label{eqn:residual recursion}
	\end{equation}
	holds over the subset $\{D\in\upT_{\vecw^\star}(\Gr_k(\upT\calM)):\snorm{D}_{\vecw^\star}<r_3\}$, where
	$$A(\eta):=\begin{pmatrix}
		q_1(\eta)+c_3r_3 & c_3r_3\\
		\eta M+c_3r_3 & q_2(\eta)+c_3r_3
	\end{pmatrix},$$
	$M:=\snorm{\vecm_3(\vecw^\star)_{\vecx^\star}}$, $F_2(D)_{\vecx^\star}$ and $\vecm_3(\vecw^\star)_{\vecx^\star}$ are defined similarly as in Eq. \eqref{eqn:projected direction},
    \begin{align*}
        \revise{F_2(D)_{\vecx^\star}:=}&~\revise{\Proj_{\upT_{\vecx^\star}\calM}\circ F_2(D)\circ\Proj_{\upT_{\vecx^\star}\calM},}\\
        \revise{\vecm_3(\vecw^\star)_{\vecx^\star}:=}&~\revise{\Proj_{\upT_{\vecx^\star}\calM}\circ\vecm_3(\vecw^\star)\circ\Proj_{\upT_{\vecx^\star}\calM},}
    \end{align*}
    and
	\begin{align*}
		&q_1(\eta):=\max\lrbrace{\abs{1-\eta\lambda_{\min}^\star},\abs{1-\eta\lambda_{\max}^\star}},\\
		&q_2(\eta):=\max\lrbrace{\abs{1-\eta\Delta\lambda_{k+1,k}^\star},\abs{1-\eta\Delta\lambda_{d1}^\star}},\\
		&\lambda_{\min}^\star:=\min_{i=1}^d\abs{\lambda_i^\star},\quad\lambda_{\max}^\star:=\max_{i=1}^d\abs{\lambda_i^\star},\\
		&\Delta\lambda_{ji}^\star:=\lambda_j^\star-\lambda_i^\star,\quad i=1,\ldots,k,~j=k+1,\ldots,d.
	\end{align*}
\end{lemma}
\begin{proof}
	We first show the well-definedness of $F$, i.e., the Riemannian distance between $\vecw^\star$ and $\upR^{\Gr_k(\upT\calM)}_{\tilde\vecw(D)}\big(\eta\cdot\vech(\tilde\vecw(D))\big)$ falls below $r_1$. By the triangle inequality of the Riemannian distance, 
	\begin{align*}
		&~\dist_{\Gr_k(\upT\calM)}\lrbracket{\vecw^\star,\upR^{\Gr_k(\upT\calM)}_{\tilde\vecw(D)}\big(\eta\cdot\vech(\tilde\vecw(D))\big)}\\
		\le&~\dist_{\Gr_k(\upT\calM)}(\vecw^\star,\tilde\vecw(D))+\dist_{\Gr_k(\upT\calM)}\lrbracket{\tilde\vecw(D),\upR^{\Gr_k(\upT\calM)}_{\tilde\vecw(D)}\big(\eta\cdot\vech(\tilde\vecw(D))\big)}\\
		\le&~\snorm{D}_{\vecw^\star}+c_1\eta\snorm{\vech(\tilde\vecw(D))}_{\revise{\tilde\vecw(D)}}\le r_2+c_1\eta\sup_{\snorm{D}_{\vecw^\star}\le r_1}\snorm{\vech(\tilde\vecw(D))}_{\revise{\tilde\vecw(D)}}<r_1,
	\end{align*}
	where the second inequality leverages Lemma \ref{lem:diffeomorphism nbd}, the fact that $\snorm{D}_{\vecw^\star}<r_2<r_1$, and the assumption \eqref{eqn:step size ub} on $\eta$, in that
	$$\eta\snorm{\vech(\tilde\vecw(D))}_{\revise{\tilde\vecw(D)}}\le\eta\sup_{\snorm{D}_{\vecw^\star}\le r_1}\snorm{\vech(\tilde\vecw(D))}_{\revise{\tilde\vecw(D)}}<r_1,$$
	the last one is again due to the assumption \eqref{eqn:step size ub} on $\eta$. The smoothness of $F$ then follows. 
	
	\par Regarding the estimate for $\snorm{F(D)}_{\vecw^\star}$, we first expand $F$ around \revise{the} origin up to first order: there exist constants $r_3\in(0,r_2]$ and $c_3\ge0$ such that
	\begin{align*}
		F(D)&=F(D)-F(0)=\upD F(0)[D]+r(D),\quad\forall~D\in\upT_{\vecw^\star}(\Gr_k(\upT\calM)):\snorm{D}_{\vecw^\star}<r_3,
	\end{align*}
	where \revise{we have used $F(0)=0$ due to 
    $$\vech(\vecw^\star)=0,\quad\upR_{\vecw^\star}^{\Gr_k(\upT\calM)}(0)=\vecw^\star,\quad\big(\Exp_{\vecw^\star}^{\Gr_k(\upT\calM)}\big)^{-1}(\vecw^\star)=0,$$} 
    $r(D):=(r_1(D),r_2(D))$ collects the higher-order terms and satisfies 
	$$\snorm{r_1(D)},~\snorm{r_2(D)_{\vecx^\star}}<c_3\lrbracket{\snorm{\bdelta}^2+\snorm{\Delta_{\vecx^\star}}^2},$$
	with $r_2(D)_{\vecx^\star}$ defined in a similar way as in Eq. \eqref{eqn:projected direction} \revise{via two compositions with projectors}. More calculations on the first-order term yield that
	\begin{align*}
		\upD F(0)[D]&=\upD\upR^{\Gr_k(\upT\calM)}(\vecw^\star,0)\lrsquare{\lrbracket{D,\eta\cdot\upD\vech(\vecw^\star)[D]}}=D+\eta\cdot\upD\vech(\vecw^\star)[D]\\
		&=\begin{pmatrix}
			\bdelta-\eta R_{\vecw^\star}\big(\hess{\calM}f(\vecx^\star)[\bdelta]\big)\\[0.1cm]
			\Delta-\eta[P^\star,[\Delta,\hess{\calM}f(\vecx^\star)]]+\eta\vecm_3(\vecw^\star)[\bdelta]
		\end{pmatrix},
	\end{align*}
	where the second equality follows from \cite[Lemma 4.21]{boumal2023introduction} and the last one uses Eq. \eqref{eqn:differential of solution map}. \revise{All} in all, we have
	\begin{align*}
		&~\snorm{F_1(D)}=\norm{\bdelta-\eta R_{\vecw^\star}\big(\hess{\calM}f(\vecx^\star)[\bdelta]\big)+r_1(D)}\\
		\le&~\norm{\Id_{\upT_{\vecx^\star}\calM}-\eta R_{\vecw^\star}\circ\hess{\calM}f(\vecx^\star)}\snorm{\bdelta}+c_3\lrbracket{\snorm{\bdelta}^2+\snorm{\Delta_{\vecx^\star}}^2}\\
		\le&~ (q_1(\eta)+c_3r_3)\snorm{\bdelta}+c_3r_3\snorm{\Delta_{\vecx^\star}},
	\end{align*}
	and
	\begin{align*}
		&~\snorm{F_2(D)_{\vecx^\star}}=\snorm{\Proj_{\upT_{\vecx^\star}\calM}\circ F_2(D)\circ\Proj_{\upT_{\vecx^\star}\calM}}\\
		=&~\norm{\Delta_{\vecx^\star}-\eta[P^\star,[\Delta_{\vecx^\star},\hess{\calM}f(\vecx^\star)]]+\revise{\eta}(\vecm_3(\vecw^\star)[\bdelta])_{\vecx^\star}+r_2(D)_{\vecx^\star}}\\
		\le&~\norm{\Id_{\upT_{P^\star}\Gr_k(\upT_{\vecx^\star}\calM)}-\eta[P^\star,[\Id_{\upT_{P^\star}\Gr_k(\upT_{\vecx^\star}\calM)},\hess{\calM}f(\vecx^\star)]]}\snorm{\Delta_{\vecx^\star}}\\
		&~+\eta\norm{\vecm_3(\vecw^\star)_{\vecx^\star}}\snorm{\bdelta}+c_3r_3\lrbracket{\norm{\bdelta}+\snorm{\Delta_{\vecx^\star}}}\\
		\le&~(q_2(\eta)+c_3r_3)\snorm{\Delta_{\vecx^\star}}+(\eta M+c_3r_3)\snorm{\bdelta}.
	\end{align*}
	The proof is complete.
\end{proof}

\begin{lemma}\label{lem:Cayley-Hamilton}
	Let $X$ be a $2\times2$ matrix with distinct eigenvalues $\mu_1$ and $\mu_2$. Then for any $\revise{m}\in\N$, it holds that
	$$X^{\revise{m}}=\frac{\mu_1^{\revise{m}}-\mu_2^{\revise{m}}}{\mu_1-\mu_2}X-\frac{\mu_2\mu_1^{\revise{m}}-\mu_1\mu_2^{\revise{m}}}{\mu_1-\mu_2}I_2.$$
\end{lemma}
\begin{proof}
	Since $X$ is $2\times2$, the Cayley-Hamilton theorem tells us that $X$ satisfies its own characteristic equation:
	\begin{equation}
		X^2-(\mu_1+\mu_2)X+\mu_1\mu_2I_2=0.
		\label{eqn:2*2 matrix characteristic equation}
	\end{equation}
	In what follows, we prove the desired result by induction. The conclusion holds obviously for the cases of $\revise{m}=0,1,2$. Now suppose that the conclusion holds for the case of $\revise{m}=\ell\in\N$ and consider the case of $\revise{m}=\ell+1$. Direct calculations yield
	\begin{align*}
		X^{\ell+1}&=X^\ell\cdot X=\lrbracket{\frac{\mu_1^\ell-\mu_2^\ell}{\mu_1-\mu_2}X-\frac{\mu_2\mu_1^\ell-\mu_1\mu_2^\ell}{\mu_1-\mu_2}I_2}X\\
		&=\frac{\mu_1^\ell-\mu_2^\ell}{\mu_1-\mu_2}\lrbracket{(\mu_1+\mu_2)X-\mu_1\mu_2I_2}-\frac{\mu_2\mu_1^\ell-\mu_1\mu_2^\ell}{\mu_1-\mu_2}X\\
		&=\frac{\mu_1^{\ell+1}-\mu_2^{\ell+1}}{\mu_1-\mu_2}X-\frac{\mu_2\mu_1^{\ell+1}-\mu_1\mu_2^{\ell+1}}{\mu_1-\mu_2}I_2.
	\end{align*}
	Here the second equality follows from mathematical induction and the third one uses Eq. \eqref{eqn:2*2 matrix characteristic equation}. As a result, the conclusion holds for the case of $\revise{m}=\ell+1$ as well. The proof is complete.
\end{proof}

\begin{theorem}[Local convergence of Algorithm \ref{alg:discretized deterministic constrained saddle dynamics}]\label{thm:local convergence of discretized algorithm}
	Suppose that $f$ is $C^3$, $\calM$ is \revise{an embedded} Riemannian submanifold of a Euclidean space $\calE$ with $\dim(\calM)=d$, \revise{$1\le k<d$}, $\vecx^\star\in\calM$ is an index-$k$ constrained SP of $f$ on $\calM$, $\hess{\calM}f(\vecx^\star)$ is nondegenerate, whose eigenvalues are $\lambda_1^\star\le\cdots\le\lambda_k^\star<0<\lambda_{k+1}^\star\le\cdots\le\lambda_d^\star$, and $P^\star\in\Gr_k(\upT_{\vecx^\star}\calM)$ is an orthogonal projector onto the lowest $k$-dimensional invariant subspace of $\hess{\calM}f(\vecx^\star)$. Let $\vecw^\star:=(\vecx^\star,P^\star)\in\Gr_k(\upT\calM)$ and $\{\vecw^{(t)}:=(\vecx^{(t)},P^{(t)})\}$ be the sequence generated by Algorithm \ref{alg:discretized deterministic constrained saddle dynamics}. If the following two assumptions hold:
	\begin{itemize}
		\item \textup{(Smallness of the step size)} the step size fulfills
		\begin{equation}
			0<\eta<\min\lrbrace{\frac{2}{\Delta\lambda_{d1}^\star},\frac{q_1(\eta)+q_2(\eta)}{M},\frac{\min\{(r_1-r_2)/c_1,r_1\}}{\sup_{\snorm{D}_{\vecw^\star}\le r_1}\snorm{\vech(\tilde\vecw(D))}_{\revise{\tilde\vecw(D)}}}}
			\label{eqn:step size ub 2}
		\end{equation}
		with $c_1$ and $r_1$ the constants in Lemma \ref{lem:diffeomorphism nbd} (associated to $\calN=\Gr_k(\upT\calM)$, $\upR^{\calN}=\upR^{\Gr_k(\upT\calM)}$ or $\Exp^{\Gr_k(\upT\calM)}$, and $\vecy=\vecw^\star$) and $r_2\in(0,r_1)$;
		
		\item \textup{(Proximity of the initial point to the desired constrained SP)} the initial point satisfies
		\begin{equation}
			\dist_{\Gr_k(\upT\calM)}(\vecw^\star,\vecw^{(0)})<\frac{r_3}{\max\lrbrace{\max\lrbrace{\frac{\mu_1(\eta)^2}{\mu_1(\eta)-\mu_2(\eta)},1}\snorm{A(\eta)},1}},
			\label{eqn:initial point}
		\end{equation}
		where the constants $c_3$ and $r_3$ are defined in Lemma \ref{lem:residual mapping} and are chosen such that
		\begin{equation}
			c_3r_3<\frac{(1-q_{\max}(\eta))^2}{\eta M+2(1-q_{\max}(\eta))}\quad\text{with}\quad q_{\max}(\eta):=\max\lrbrace{q_1(\eta),q_2(\eta)},
			\label{eqn:higher order radius}
		\end{equation}
		and $\mu_1(\eta)$, $\mu_2(\eta)$ are respectively
		\begin{equation}
			\begin{aligned}
				\mu_1(\eta)&:=\frac{q_1(\eta)+q_2(\eta)+2c_3r_3+\sqrt{(q_1(\eta)-q_2(\eta))^2+4c_3r_3(\eta M+c_3r_3)}}{2},\\
				\mu_2(\eta)&:=\frac{q_1(\eta)+q_2(\eta)+2c_3r_3-\sqrt{(q_1(\eta)-q_2(\eta))^2+4c_3r_3(\eta M+c_3r_3)}}{2},
			\end{aligned}
			\label{eqn:2*2 matrix eigenvalue}
		\end{equation}
	\end{itemize}
	then $\{\vecw^{(t)}\}$ converges linearly to $\vecw^\star$.
\end{theorem}
\begin{proof}
	We first show by induction that $\dist_{\Gr_k(\upT\calM)}(\vecw^\star,\vecw^{(t)})<r_3$ holds for any $t$, so that $\Exp_{\vecw^\star}^{\Gr_k(\upT\calM)}$ remains a diffeomorphism and Eq. \eqref{eqn:residual recursion} is applicable for all the iterates. The case of $t=0$ is evident due to the assumption \eqref{eqn:initial point} on the initial point. Now suppose that the statement holds for the case of $t=\ell\in\N$ and consider the case of $t=\ell+1$. By induction and Eq. \eqref{eqn:residual recursion}, there exists a unique $D^{(t+1)}\in\upT_{\vecw^\star}(\Gr_k(\upT\calM))$ such that $\snorm{D^{(t+1)}}_{\vecw^\star}<r_3$ and $D^{(t+1)}=F(D^{(t)})$ hold for $t=0,\ldots,\ell-1$, and that
	\begin{equation}
		\begin{pmatrix}
			\snorm{F_1(D^{(\ell)})}\\
			\snorm{F_2(D^{(\ell)})_{\vecx^\star}}
		\end{pmatrix}\le A(\eta)^{\ell+1}\begin{pmatrix}
			\snorm{\bdelta^{(0)}}\\
			\snorm{\Delta_{\vecx^\star}^{(0)}}
		\end{pmatrix}.
		\label{eqn:induction residual}
	\end{equation}
	It is not difficult to verify that $\mu_1(\eta)$ and $\mu_2(\eta)$ in Eq. \eqref{eqn:2*2 matrix eigenvalue} are the two distinct real eigenvalues of $A(\eta)$, and both of them belong to $(0,1)$:
	\begin{align*}
		2\mu_2(\eta)&= q_1(\eta)+q_2(\eta)+2c_3r_3-\sqrt{(q_1(\eta)-q_2(\eta))^2+4c_3r_3(\eta M+c_3r_3)}\\
		&> q_1(\eta)+q_2(\eta)+2c_3r_3-\sqrt{(q_1(\eta)-q_2(\eta))^2+4c_3r_3(q_1(\eta)+q_2(\eta)+c_3r_3)}\\
		&=q_1(\eta)+q_2(\eta)+2c_3r_3-\sqrt{(q_1(\eta)+q_2(\eta)+2c_3r_3)^2-4q_1(\eta)q_2(\eta)}\ge0,
	\end{align*}
	where the strict inequality is due to the assumption \eqref{eqn:step size ub 2}, and
	\begin{align*}
		&~2\mu_1(\eta)\\
        =&~q_1(\eta)+q_2(\eta)+2c_3r_3+\sqrt{(q_1(\eta)-q_2(\eta))^2+4c_3r_3(\eta M+c_3r_3)}\\
		\le&~ q_1(\eta)+q_2(\eta)+2c_3r_3+\abs{q_1(\eta)-q_2(\eta)}+2\sqrt{c_3r_3(\eta M+c_3r_3)}\\
		=&~2\lrbracket{c_3r_3+q_{\max}(\eta)+\sqrt{c_3r_3(\eta M+c_3r_3)}}\\
		<&~2\lrbracket{\frac{(1-q_{\max}(\eta))^2}{\eta M+2(1-q_{\max}(\eta))}+q_{\max}(\eta)+\sqrt{\frac{(1-q_{\max}(\eta))^2(\eta M+1-q_{\max}(\eta))^2}{(\eta M+2(1-q_{\max}(\eta)))^2}}}\\
		=&~2\lrbracket{\frac{(1-q_{\max}(\eta))^2}{\eta M+2(1-q_{\max}(\eta))}+q_{\max}(\eta)+\frac{(1-q_{\max}(\eta))(\eta M+1-q_{\max}(\eta))}{\eta M+2(1-q_{\max}(\eta))}}=2,
	\end{align*}
	where the strict inequality uses the assumption \eqref{eqn:higher order radius} and the third equality uses $q_1(\eta),q_2(\eta)<1$ from the assumption \eqref{eqn:step size ub 2}. By Lemma \ref{lem:Cayley-Hamilton}, 
	\begin{align}
		&~A(\eta)^{\ell+1}\begin{pmatrix}
			\snorm{\bdelta^{(0)}}\\
			\snorm{\Delta_{\vecx^\star}^{(0)}}
		\end{pmatrix}\label{eqn:residual upper}\\
		=&~\lrbracket{\frac{\mu_1(\eta)^{\ell+1}-\mu_2(\eta)^{\ell+1}}{\mu_1(\eta)-\mu_2(\eta)}A(\eta)-\frac{\mu_2(\eta)\mu_1(\eta)^{\ell+1}-\mu_1(\eta)\mu_2(\eta)^{\ell+1}}{\mu_1(\eta)-\mu_2(\eta)}I_2}\begin{pmatrix}
			\snorm{\bdelta^{(0)}}\\
			\snorm{\Delta_{\vecx^\star}^{(0)}}
		\end{pmatrix}\nonumber\\
		\le&~\max\lrbrace{\frac{\mu_1(\eta)^2}{\mu_1(\eta)-\mu_2(\eta)},1}A(\eta)\begin{pmatrix}
			\snorm{\bdelta^{(0)}}\\
			\snorm{\Delta_{\vecx^\star}^{(0)}}
		\end{pmatrix},\nonumber
	\end{align}
	where the inequality is due to $0<\mu_2(\eta)<\mu_1(\eta)<1$. Eqs. \eqref{eqn:induction residual}, \eqref{eqn:residual upper}, and the assumption \eqref{eqn:initial point} on the initial point then together imply 
	\begin{align*}
		\snorm{F(D^{(\ell)})}_{\vecw^\star}^2&=\snorm{F_1(D^{(\ell)})}^2+\snorm{F_2(D^{(\ell)})_{\vecx^\star}}^2\revise{\le}\norm{A(\eta)^{\ell+1}\begin{pmatrix}
				\snorm{\bdelta^{(0)}}\\
				\snorm{\Delta_{\vecx^\star}^{(0)}}
		\end{pmatrix}}^2\\
		&\le\max\lrbrace{\frac{\mu_1(\eta)^2}{\mu_1(\eta)-\mu_2(\eta)},1}^2\norm{A(\eta)\begin{pmatrix}
				\snorm{\bdelta^{(0)}}\\
				\snorm{\Delta_{\vecx^\star}^{(0)}}
		\end{pmatrix}}^2\\
		&\le\max\lrbrace{\frac{\mu_1(\eta)^2}{\mu_1(\eta)-\mu_2(\eta)},1}^2\snorm{A(\eta)}^2\snorm{D^{(0)}}_{\vecw^\star}^2<r_3^2.
	\end{align*}
	Consequently, there exists a unique $D^{(\ell+1)}\in\upT_{\vecw^\star}(\Gr_k(\upT\calM))$ such that $\snorm{D^{(\ell+1)}}<r_3$ and $D^{(\ell+1)}=F(D^{(\ell)})$. The statement thus holds for any $t\in\N$. 
	
	\par Now by Eqs. \eqref{eqn:induction residual} and \eqref{eqn:residual upper}, we have for any $t\ge1$ that
	$$\begin{pmatrix}
		\snorm{\bdelta^{(t)}}\\
		\snorm{\Delta_{\vecx^\star}^{(t)}}
	\end{pmatrix}\le\lrbracket{\frac{\mu_1(\eta)^{t}-\mu_2(\eta)^{t}}{\mu_1(\eta)-\mu_2(\eta)}A(\eta)-\frac{\mu_2(\eta)\mu_1(\eta)^{t}-\mu_1(\eta)\mu_2(\eta)^{t}}{\mu_1(\eta)-\mu_2(\eta)}I_2}\begin{pmatrix}
		\snorm{\bdelta^{(0)}}\\
		\snorm{\Delta_{\vecx^\star}^{(0)}}
	\end{pmatrix},$$
	which combined with $0\le\mu_2(\eta)<\mu_1(\eta)<1$ and $\snorm{D^{(t)}}_{\vecw^\star}^2=\snorm{\bdelta^{(t)}}^2+\snorm{\Delta_{\vecx^\star}^{(t)}}^2$ implies that $D^{(t)}$ converges linearly to 0, or equivalently, $\{\vecw^{(t)}\}$ converges linearly to $\vecw^\star$ as $t\to\infty$. The proof is completed.
\end{proof}

\begin{remark}[Satisfiability of the assumptions]\label{rem:satisfiability of assumptions}
	It might seem to be difficult at first glance to check the satisfiability of the assumptions \eqref{eqn:step size ub 2} and \eqref{eqn:higher order radius}. Indeed, after noting that $q_1(\eta)$, $q_2(\eta)\to1-$ as $\eta\to0+$, the assumption \eqref{eqn:step size ub 2} can always be met with a sufficiently small $\eta$. Once $\eta$ is fixed, one could choose an $r_3$ small enough to meet the assumption \eqref{eqn:higher order radius}, since $c_3$, as a function of $r_3$, is uniformly bounded in a small neighborhood of \revise{the} origin. 

    \par \revise{It is worth mentioning that the constants showing up in Eqs. \eqref{eqn:step size ub 2} and \eqref{eqn:higher order radius} are not intended as directly computable quantities. They depend on, e.g., local smoothness of the objective function, Riemannian manifold, and retraction, and are introduced mainly for theoretical convenience. In addition, conditions \eqref{eqn:step size ub 2} and \eqref{eqn:initial point} are sufficient but not definitely necessary for the local result, and the established error bound is not guaranteed to be sharp. Nonetheless, we shall point out that Theorem \ref{thm:local convergence of discretized algorithm} should better be understood as the first local iterate convergence result for a discretized saddle-search algorithm in the embedded-submanifold setting. Completely closing the gap between theory and practice is important and will be left for future investigations.}
\end{remark}

\begin{remark}[Local convergence rates and condition numbers]\label{rem:local convergence rates and condition numbers}
	Lemma \ref{lem:residual mapping} provides some information about the local convergence rates of $\vecx$- and $P$-residuals. Indeed, as long as $r_3$ is sufficiently small, the contraction matrix $A(\eta)$ in Eq. \eqref{eqn:residual recursion} becomes close to diagonal, with entries dominated by $q_1(\eta)$ and $q_2(\eta)$. They can then be viewed as the contractive factors of $\vecx$- and $P$-residuals, respectively. 
	
	\par Moreover, if two different step sizes $\eta_{\vecx}$ and $\eta_P$ are used, it is possible for us to approach the best linear convergence rates. By the definitions of $q_1(\eta)$ and $q_2(\eta)$ in Lemma \ref{lem:residual mapping}, it is easy to show that the step sizes minimizing these two terms are
	\begin{align}
		\eta_{\vecx}^\star&:=\frac{2}{\lambda_{\min}^\star+\lambda_{\max}^\star}=\frac{2}{\min_{i=1}^d|\lambda_i^\star|+\max_{i=1}^d|\lambda_i^\star|},\label{eqn:theoretical step size x}\\
		\eta_{P}^\star&:=\frac{2}{\Delta\lambda_{k+1,k}^\star+\Delta\lambda_{d1}^\star}=\frac{2}{\lambda_d^\star+\lambda_{k+1}^\star-\lambda_k^\star-\lambda_1^\star}.\label{eqn:theoretical step size v}
	\end{align}
	The condition numbers are thus estimated by
	\begin{equation}
		\kappa_{\vecx}:=\frac{\lambda_{\max}^\star}{\lambda_{\min}^\star}=\frac{\max_{i=1}^d|\lambda_i^\star|}{\min_{i=1}^d|\lambda_i^\star|},\quad\kappa_P:=\frac{\Delta\lambda_{d1}^\star}{\Delta\lambda_{k+1,k}^\star}=\frac{\lambda_d^\star-\lambda_1^\star}{\lambda_{k+1}^\star-\lambda_k^\star}.
		\label{eqn:condition numbers}
	\end{equation}
	We shall note that the best step sizes in practice are not necessarily given by Eqs. \eqref{eqn:theoretical step size x} and \eqref{eqn:theoretical step size v} exactly, since Eq. \eqref{eqn:residual recursion} provides only an estimate from above and the terms other than $q_1(\eta)$ and $q_2(\eta)$ in the contraction matrix $A(\eta)$ might not be negligible. For numerical illustrations, see Section \ref{subsec:linear eigenvalue problem}.
\end{remark}

\begin{remark}[Comparison with the existing results for discretized algorithms]\label{rem:removal of unnecessary assumptions}
	All the existing results \cite{gao2015iterative,gould2016dimer,levitt2017convergence,luo2025accelerated,luo2022convergence,zhang2012shrinking} are set in \revise{\sout{the}} unconstrained settings\footnote{\revise{In \cite{zhang2023discretization}, for a fixed time interval, the authors derives error estimates between the solution of the continuous constrained high-index saddle dynamic\revise{\sout{s}} and its numerical approximation as the time step tends to zero. This is of a different nature from the iterate-to-saddle convergence result established in this work.}}. Therefore, it makes no sense to conduct a direct comparison. Nevertheless, we shall point out that the local convergence rate and condition number for the $\vecx$-residual in Remark \ref{rem:local convergence rates and condition numbers} recover the ones in \cite{luo2022convergence} if we let $\calM=\calE$, yet with weaker assumptions on the eigenvalues of the Hessian (see Remark \ref{rem:weaker assumptions}). Moreover, since we treat $\vecx$ and $P$ equally on the Grassmann bundle $\Gr_k(\upT\calM)$, by leveraging \revise{\sout{the}} geometrical tools, we manage to establish the explicit local convergence rates and condition numbers for both parts in the manifold constrained settings (see Remark \ref{rem:local convergence rates and condition numbers}), and do not require unnecessary (or even uncheckable) assumptions; for example, the angle assumption in \cite{luo2022convergence}, which reads
	$$\exists~\alpha\in[0,1),~\st~\norm{V^{(t)}(V^{(t)})^\top-V(\vecx^{(t)})V(\vecx^{(t)})^\top}\le\alpha,~\forall~t.$$
	Here, $V(\vecx^{(t)})\in\St_k(\upT_{\vecx^{(t)}}\calM)$ spans exactly the lowest $k$-dimensional invariant subspace of $\hess{\calM}f(\vecx^{(t)})$. 
\end{remark}

\begin{example}[Condition numbers at index-1 constrained SPs]\label{exm:condition number at index 1 Grassmann}
	Consider a linear objective function on $\calM=\Gr_p(\R^n)$ ($p>1$) defined by $f(P)=\tr(PA)$ for any $P\in\calM$, where $A\in\R_{\sym}^{n\times n}$. It is not difficult to obtain \revise{a} closed-form expression \revise{for the} Riemannian Hessian: 
	$$\hess{\calM}f(P)[\Gamma]=[P,[\Gamma,A]],\quad\forall~P\in\calM,~\Gamma\in\upT_P\calM.$$ 
	For simplicity, suppose that $A$ is diagonal, i.e., $A=\diag(\eps_1,\ldots,\eps_n)$ with $\eps_1\le\cdots\le\eps_n$, and that all the sums of $p$ eigenvalues of $A$ are different. Under this assumption, one could verify that the unique global minimizer (GM) and unique index-1 constrained SP of $f$ on $\calM$ are respectively 
	$$P_{\rm GM}:=\sum_{i=1}^p\vece_i\vece_i^\top\quad\text{and}\quad P_{\rm SP}:=\sum_{i=1}^{p-1}\vece_i\vece_i^\top+\vece_{p+1}\vece_{p+1}^\top,$$
	where $\vece_i$ is the $i$-th unit vector in $\R^n$ ($i=1,\ldots,n$). The eigenpairs of $\hess{\calM}f(P_{\rm GM})$ are given explicitly by 
	$$\lrbrace{\left.\lrbracket{\eps_a-\eps_i,\frac{1}{\sqrt2}\lrbracket{\vece_i\vece_a^\top+\vece_a\vece_i^\top}}\right|i=1,\ldots,p,~a=p+1,\ldots,n},$$
	and for $\hess{\calM}f(P_{\rm SP})$, 
	$$\lrbrace{\left.\lrbracket{\eps_a-\eps_i,\frac{1}{\sqrt2}\lrbracket{\vece_i\vece_a^\top+\vece_a\vece_i^\top}}\right|i=1,\ldots,p-1,p+1,~a=p,p+2,\ldots,n}.$$
	The condition number at $P_{\rm GM}$ is 
	$$\kappa_{P_{\rm GM}}=\frac{\eps_n-\eps_1}{\eps_{p+1}-\eps_p},$$
	which is well known in the literature. The denominator is also called the ``eigengap'' in some applications. Regarding the condition number at $P_{\rm SP}$, with reference to Eq. \eqref{eqn:condition numbers}, we could identify $\lambda_1^\star=\eps_p-\eps_{p+1}<0$, $\lambda_2^\star=\min\{\eps_{p+2}-\eps_{p+1},\eps_p-\eps_{p-1}\}>0$, and $\lambda_d^\star=\eps_n-\eps_1>0$ (with $d=\dim(\calM)=p(n-p)$), and therefore, 
	$$\lambda_{\min}^\star=\min\lrbrace{\eps_p-\eps_{p-1},\eps_{p+1}-\eps_p,\eps_{p+2}-\eps_{p+1}},\quad\lambda_{\max}^\star=\eps_n-\eps_1,$$
	which implies 
	$$\kappa_{P_{\rm SP}}=\frac{\eps_n-\eps_1}{\min\lrbrace{\eps_p-\eps_{p-1},\eps_{p+1}-\eps_p,\eps_{p+2}-\eps_{p+1}}}.$$
	Comparing $\kappa_{P_{\rm SP}}$ with $\kappa_{P_{\rm GM}}$, we observe that finding the index-1 constrained SPs can be worse conditioned. Moreover, it asks for not only a positive eigengap, but also nondegeneracy above $\eps_{p+1}$ and below $\eps_p$. Another useful message is that, if the problem is reformulated on the Stiefel manifold $\St_p(\R^n)$ and is treated by the saddle search algorithms in the existing works \cite{li2015gentlest,liu2023constrained,yin2022constrained,zhang2012constrained}, we could anticipate their poor performance due to the inherent degeneracy. For illustrations, see Section \ref{subsec:linear eigenvalue problem}.
\end{example}

\section{Numerical experiments}\label{sec:numerical experiments}

In this part, we report numerical results on the linear eigenvalue problem and electronic excited-state calculations. Both tasks sit on the Grassmann manifold. In particular, with the former one, we elucidate the importance of using nonredundant parametrizations in finding SPs (cf. Example \ref{exm:condition number at index 1 Grassmann}) and show the influence of problem and algorithm settings on the convergence rates (cf. Remark \ref{rem:local convergence rates and condition numbers}). \revise{We also conduct comparisons between algorithms implementing different retractions (see Remark \ref{rem:example retractions}) on the linear eigenvalue problem, and the results can be found in Appendix \ref{appsec:numerical comparisons between retractions}. The code for producing all the numerical results is available at \url{https://github.com/huyukuan/SadMan}.}

\subsection{Linear eigenvalue problem}\label{subsec:linear eigenvalue problem}

Given a real symmetric matrix $A\in\R_{\sym}^{n\times n}$, the linear eigenvalue problem amounts to finding a critical point $X^\star$ of the quadratic function $f_{\St}(X):=\frac12\trace(X^\top AX)$ on the Stiefel manifold $\St_p(\R^n)$. It is well known that the columns of any such point span a $p$-dimensional invariant subspace of $A$. Moreover, the function value $f_{\St}(X^\star)$ is exactly half the sum of $p$ eigenvalues of $A$. If all the sums ($\binom{n}{p}$ in total) are different from each other, the indices of the saddle points can be indicated by their function values, in that $f_{\St}(X_{{\rm SP},k}^\star)<f_{\St}(X_{{\rm SP},\ell}^\star)$ holds for the index-$k$ and index-$\ell$ constrained SPs whenever $k<\ell$ (the reverse might not be true). Note that $f_{\St}$ is invariant under the transformation $X\mapsto XQ$ for any $Q\in\calO(p)$. Therefore, one could instead consider a linear function $f_{\Gr}(P):=\frac12\trace(PA)$ on the Grassmann manifold $\Gr_p(\R^n)$ \revise{where any critical point} $P^\star$ is the orthogonal projector onto a $p$-dimensional invariant subspace of $A$. The statements regarding the function values and indices of constrained SPs remain valid similarly. 

\medskip

\par\noindent\textit{Implementation details.} In what follows, we construct the test matrix $A$ as $A=U\Sigma U^\top$, where $U\in\calO(n)$ is obtained from the orthonormalization of a matrix whose entries are random numbers drawn independently and identically from the standard normal distribution (random seed = 0), and $\Sigma=\diag(\sigma_1,\ldots,\sigma_n)\in\R^{n\times n}$ with $\sigma_i:=\xi^{i-n}$ ($i=1,\ldots,n$) for a parameter $\xi>1$. The problem data $n$, $p$, and $\xi$ are specified later. We consider Algorithm \ref{alg:discretized deterministic constrained saddle dynamics lifted} of Stiefel and Grassmann versions ($\calM=\St_p(\R^n)$ and $\calM=\Gr_p(\R^n)$). \revise{The retraction consists of a retraction on $\calM$, direct addition and metric projection (see Lemma \ref{lem:general retraction on Grassmann bundle}, Remark \ref{rem:example retractions} for their definitions and Appendix \ref{appsec:numerical comparisons between retractions} for a numerical justification). For the retraction on $\calM$, we adopt the QR decomposition-based second-order retraction on $\St_p(\R^n)$ and the exponential mapping on $\Gr_p(\R^n)$, respectively.} 
The step sizes will \revise{be} specified later\footnote{\revise{We do not claim that these step sizes either necessarily satisfy the assumptions of Theorem \ref{thm:local convergence of discretized algorithm} or need to be optimal; see also Remark \ref{rem:satisfiability of assumptions}. The same also applies to Section \ref{subsec:electronic excited-state calcs}.}}. The maximum iteration number and convergence tolerance are respectively set as $maxit=\infty$ and $tol=10^{-8}$. If not stated, the initial feasible points $(X^{(0)},V^{(0)})$ and $(P^{(0)},\Gamma^{(0)})$ are randomly generated as follows: with a given random seed,
\begin{verbatim}
from jax import random, numpy
key = random.PRNGKey(randseed) # random seed
key1, key2 = random.split(key)
X, _ = numpy.linalg.qr(random.normal(key1, (n, p)))
P = X @ X.T
Gamma = random.normal(key2, (n, n))
Gamma = (Gamma + Gamma.T) / 2.0
Gamma = proj_tangent(P, Gamma) # project onto the tangent space
V = Gamma @ X # horizontal lift
\end{verbatim}

\par
\begin{figure}[!t]
	\centering
	\includegraphics[width=\linewidth]{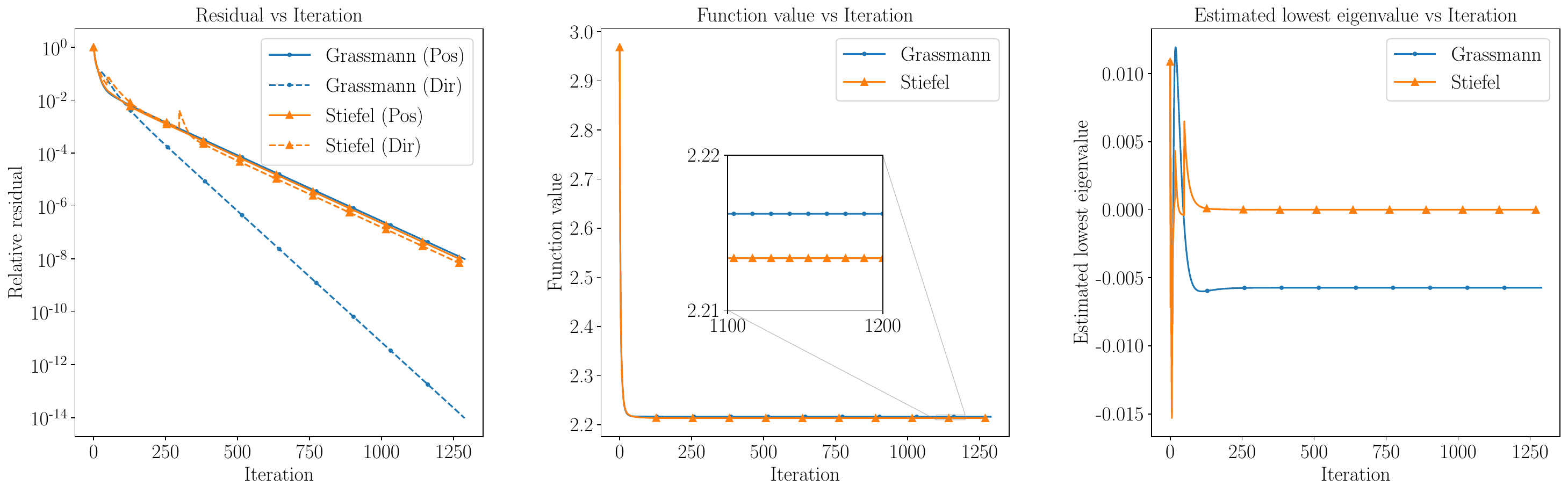}
	\caption{Convergence curves of Algorithm \ref{alg:discretized deterministic constrained saddle dynamics lifted} of Stiefel (orange lines) and Grassmann (blue lines) versions on the linear eigenvalue problem ($n=64$, $p=8$, and $\xi=1.01$). Left: position (solid lines) and direction (dashed lines) relative residuals vs iteration. Middle: function value vs iteration. Right: estimated lowest eigenvalue of Riemannian Hessian vs iteration.}
	\label{fig:LEP compare Gr St}
\end{figure}

\noindent\textit{Importance of nonredundant parametrizations.} Let $n=64$, $p=8$, and $\xi=1.01$. Following the above problem description, we run both the Stiefel and Grassmann versions of Algorithm \ref{alg:discretized deterministic constrained saddle dynamics lifted} to find the index-1 constrained SPs of $f_{\St}$ and $f_{\Gr}$, respectively. The initial feasible points are generated with random seed $=1$. The step sizes are specified as \revise{$\eta_{\St}=\eta_{\Gr}=2$}. 
The convergence curves are depicted in Figure \ref{fig:LEP compare Gr St}. 

\par We first observe from the left panel of Figure \ref{fig:LEP compare Gr St} the convergence of both versions of the algorithm. The correctness of the results can be checked by comparing the function values; see the middle panel of Figure \ref{fig:LEP compare Gr St}. Due to the construction of $A$, the function value at the index-1 constrained SP is exactly $\frac12(\sum_{i=1}^{p-1}\sigma_i+\sigma_{p+1})\approx2.216$. Initialized from the same point (more precisely, the Stiefel version is fed with the representative of the initial point for the Grassmann version), the Grassmann version converges well to the desired constrained SP, while the Stiefel one gets trapped to the global minimizer, whose function value is $\frac12\sum_{i=1}^p\sigma_i\approx2.213$. This comparison can also be seen in the right panel of Figure \ref{fig:LEP compare Gr St}, where the Riemannian Hessian at the point given by the Stiefel version is found to be positive semidefinite. 

\par We remark that the failure of the Stiefel version should be ascribed to the parametrization redundancy. The formulation on the Stiefel manifold does not take into account the quotient structure, as explained in the first paragraph of this subsection. The redundancy brings degeneracies above and/or below the eigengap and leads to an infinite condition number at the desired index-1 constrained SP (cf. Example \ref{exm:condition number at index 1 Grassmann}). The deficiency cannot be neglected as in the task of finding the global minimizer because the update direction for the position is not a horizontal lift of that in the Grassmann version. Roughly, suppose that the initial point lies in a region where the Riemannian Hessian is positive semidefinite (say, around the global minimizer) and the direction step converges quickly to the lowest eigenvector of the Riemannian Hessian. Due to the parametrization redundancy, the lowest eigenvalue of the Riemannian Hessian is zero and the direction variable $V(t)$ is almost vertical. This implies that climbing up along $V(t)$ makes little difference in the first order and the position variable mainly follows the gradient flow down to the minimizer. 

\par The above arguments are supported with numerical results. We run the Grassmann and Stiefel algorithms with the same problem and algorithm settings as before and randomly sampled initial points perturbed from the global minimizer and index-1 constrained SP. Specifically, we consider the perturbation level $\beta\in\{10^{-3},10^{-2},\ldots,10^1\}$ and each of them is tested with 100 independent random samples (random seed = $0\sim99$):
\begin{verbatim}
key = random.PRNGKey(randseed) # random seed
key1, key2 = random.split(key)
X, _ = numpy.linalg.qr(X_ref + beta * random.normal(key1, (n, p)))
# repeat previous procedures to create P, Gamma, and V
\end{verbatim}
Here $X_{\rm ref}$ is either the global minimizer or index-1 constrained SP. The empirical success rates of the two algorithms at different perturbation levels are recorded in Table \ref{tab:LEP perturbation}. 
\begin{table}[!t]
	\centering
	\caption{Empirical success rates in finding the index-1 constrained SP of Algorithm \ref{alg:discretized deterministic constrained saddle dynamics lifted} of Grassmann and Stiefel versions starting from the neighborhoods of the global minimizer and index-1 constrained SP of the linear eigenvalue problem ($n=64$, $p=8$, and $\xi=1.01$).}
	\label{tab:LEP perturbation}
	\resizebox{\linewidth}{!}{\begin{tabular}{c|rrrrr|rrrrr}
			\toprule
			\multirow{2}{*}{Algorithms} & \multicolumn{5}{c}{Perturbation from global minimizer} & \multicolumn{5}{|c}{Perturbation from constrained SP} \\\cmidrule{2-11}
			& $10^{-3}$ & $10^{-2}$ & $10^{-1}$ & $10^0$ & $10^1$ & $10^{-3}$ & $10^{-2}$ & $10^{-1}$ & $10^0$ & $10^1$\\\midrule
			Grassmann & 100\% & 100\% & 100\% & 100\% & 100\% & 100\% & 100\% & 100\% & 100\% & 100\%\\
			Stiefel & 2\% & 0\% & 1\% & 38\% & 36\% & 100\% & 99\% & 82\% & 34\% & 34\%\\\bottomrule
	\end{tabular}}
\end{table}
It turns out that the Grassmann version converges surprisingly well across the board to the index-1 constrained SP, whereas the Stiefel one often gets attracted by the global minimizer, with an empirical success rate of approximately 30\% if randomly initialized. The rate can be tremendously improved if the initial point is selected around the desired index-1 constrained SP. We only use the Grassmann version for the subsequent experiments.

\medskip


\par
\begin{figure}[!t]
    \centering
    \includegraphics[width=\linewidth]{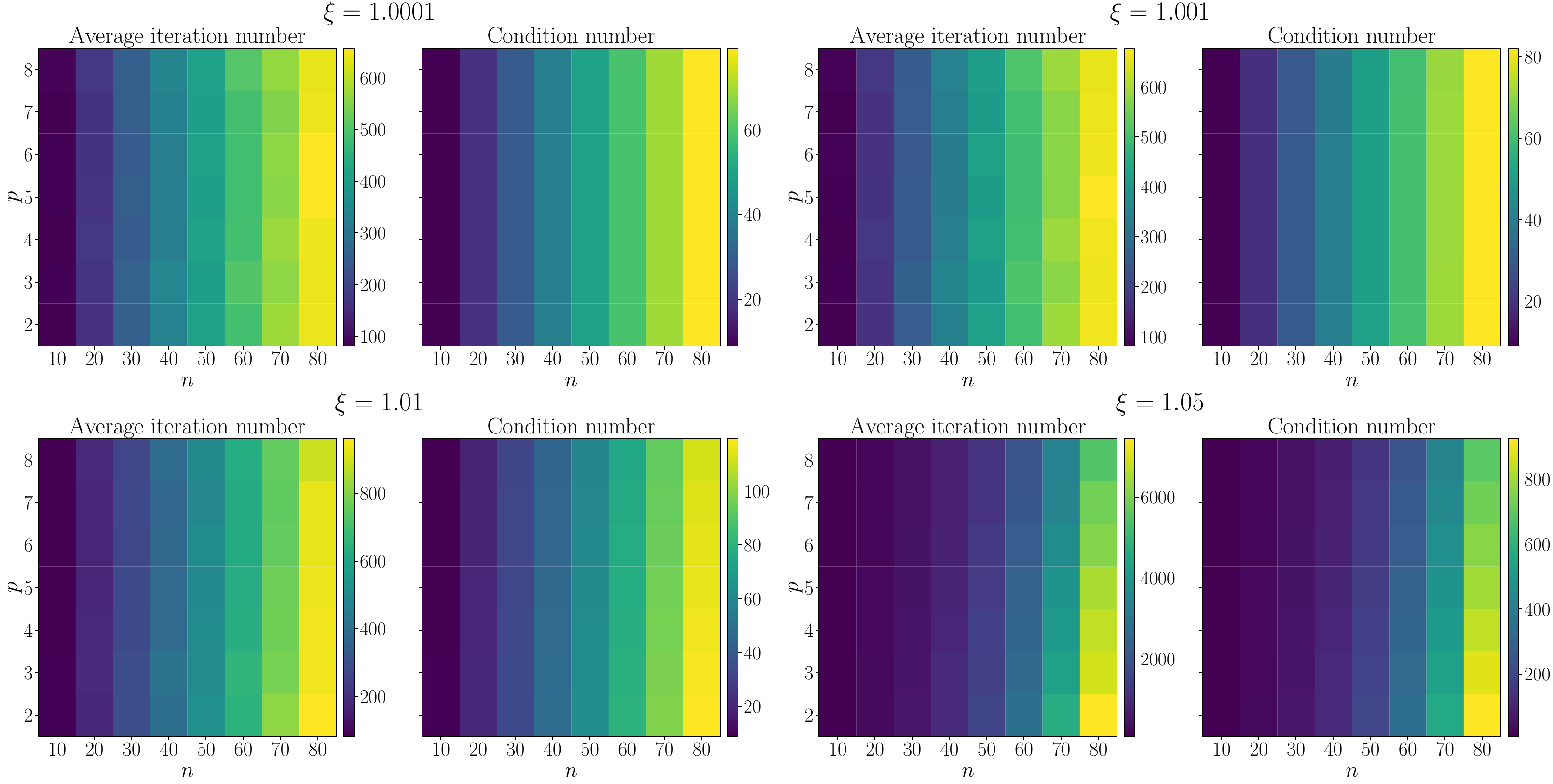}
    \caption{\revise{Average iteration number and estimated condition number vs $(n,p)$ on the linear eigenvalue problem class. Upper left two panels: $\xi=1.0001$. Upper right two panels: $\xi=1.001$. Lower left two panels: $\xi=1.01$. Lower right two panels: $\xi=1.05$.}}
    \label{fig:problem data xi}
\end{figure}

 \noindent\textit{Influence of problem and algorithm settings on the convergence rates.} We investigate the performance of Algorithm \ref{alg:discretized deterministic constrained saddle dynamics lifted} in finding the index-1 constrained SP with different problem data ($n$, $p$, and $\xi$) and algorithm settings (initialization and step sizes). First consider varying the problem data in $n\in\{10,20,\ldots,80\}$, $p\in\{2,3,\ldots,8\}$, $\xi\in\{1.0001, 1.001, 1.01, 1.05\}$. The step sizes are determined by the estimates \eqref{eqn:theoretical step size x} and \eqref{eqn:theoretical step size v}. For each problem setting, we perform ten independent runs of the algorithm with the initial points randomly selected at the perturbation level of $10^{-3}$ from the index-1 constrained SP (random seed = $0\sim9$). The average iteration numbers and the condition numbers estimated through Eq. \eqref{eqn:condition numbers}, $\kappa=\frac{\xi^{n-1}-1}{\xi^{p-2}(\xi-1)}$, are shown with \revise{heatmaps in Figure \ref{fig:problem data xi}}.
The trends of the iteration numbers are found to align qualitatively well with those of the estimated condition numbers, showcasing the validity of our theoretical analysis. 

\par
\begin{figure}[htb]
	\centering
	\includegraphics[width=\linewidth]{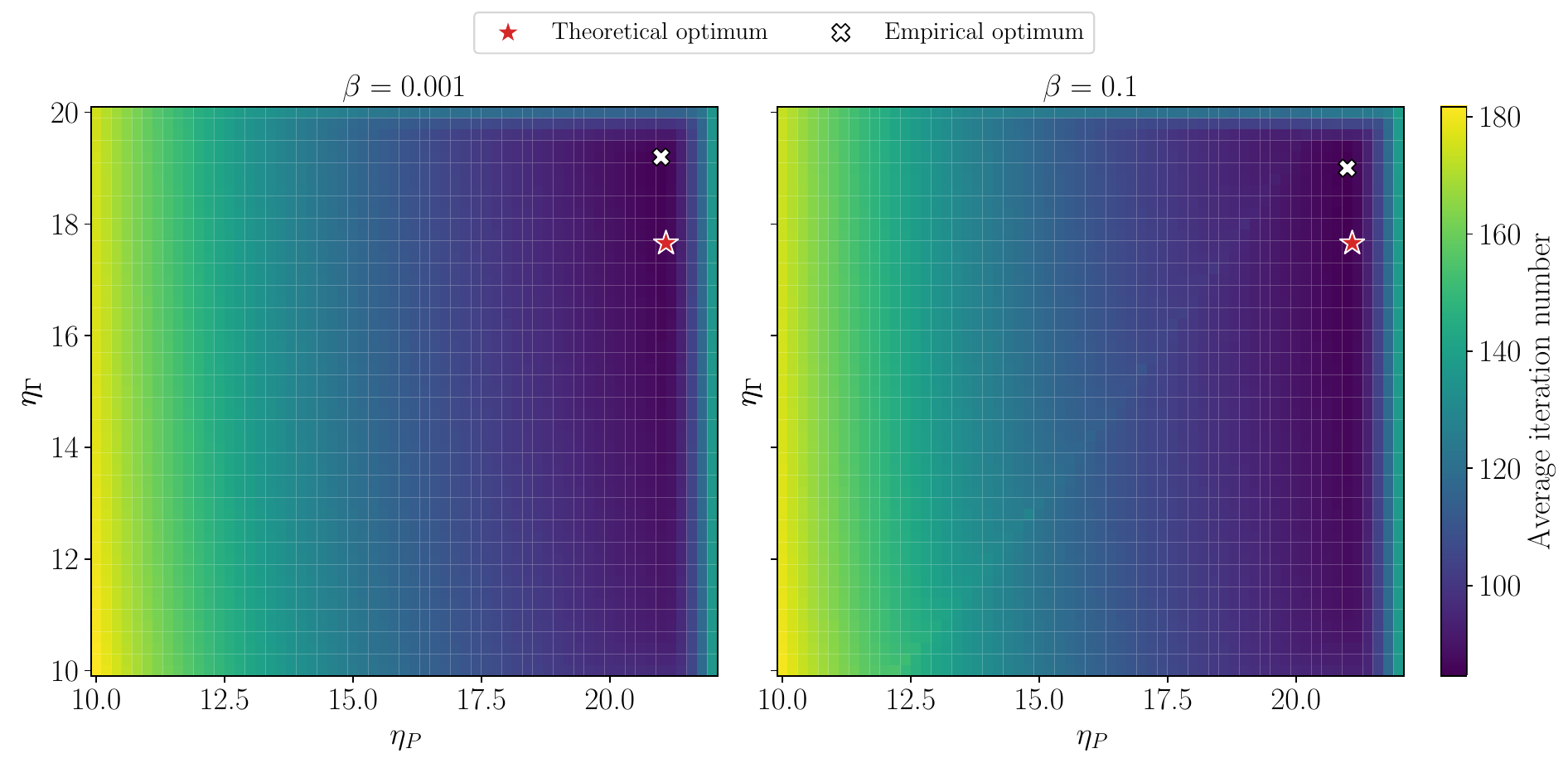}
	\caption{Average iteration number vs $(\eta_P,\eta_\Gamma)$ on the linear eigenvalue problem ($n=10$, $p=2$, and $\xi=1.01$). \revise{The theoretical and empirical best pairs of step sizes are marked out by red stars and white crosses, respectively.} Left: perturbation level of $10^{-3}$. Right: perturbation level of $10^{-1}$.}
	\label{fig:algorithm setting}
\end{figure}

 Next we consider varying the algorithm settings. We fix the problem data to be $n=10$, $p=2$, and $\xi=1.01$. Since the best step sizes are estimated by Eqs. \eqref{eqn:theoretical step size x} and \eqref{eqn:theoretical step size v} to be $(\eta_P^\star,\eta_\Gamma^\star)\approx(21.096,17.656)$ for this instance, the step sizes $\eta_P$ and $\eta_\Gamma$ are varied \revise{respectively within $\eta_P\in\{10,10.2,10.4,\ldots,21.8,22\}$ and $\eta_{\Gamma}\in\{10,10.2,10.4,\ldots,19.8,20\}$}\footnote{\revise{In fact, we have also tested step sizes outside these two boxes, but observed non-convergence of the algorithm.}}. For each pair of $(\eta_P,\eta_\Gamma)$, we perform ten independent runs of the algorithm with the initial points randomly selected at the perturbation level of $10^{-3}$ or $10^{-1}$ from the index-1 constrained SP (random seed $=0\sim9$). Similarly, we visualize the results with \revise{heatmaps} in Figure \ref{fig:algorithm setting}.
For this test instance, the least iteration numbers are achieved when the pair \revise{$(\eta_P,\eta_\Gamma)=(21,19.2)$ given $\beta=0.001$ and $(\eta_P,\eta_\Gamma)=(21,19)$ given $\beta=0.1$}. 
We also find that the step size $\eta_P$ for the position variable is far more pivotal than $\eta_\Gamma$ for the direction variable. Incidentally, \revise{the empirically best pairs of step sizes are close to the estimated pair $(\eta_P^\star,\eta_\Gamma^\star)\approx(21.096,17.656)$. The modest shift in the direction step size is consistent with} the residual reduction inequality \eqref{eqn:residual recursion} \revise{being an upper bound together with the influence of other coupling and higher-order terms. Performance at the estimated pair is already near-optimal over the tested grid.}  

\medskip

\par \noindent\textit{Computation of higher-index constrained SPs.} We proceed to compute the constrained SPs of all indices for an instance. We fix the problem data to be $n=10$, $p=2$, and $\xi=1.01$ and the step sizes to be \revise{$\eta_P=\eta_\Gamma=12$}. The possible index of the constrained SP on this instance is at most 15. For each index in $\{0,\ldots,15\}$, we perform \revise{400} runs of Algorithm \ref{alg:discretized deterministic constrained saddle dynamics lifted} from randomly generated initial points (\revise{random seed = $0\sim399$}). The indices and the function values at the obtained constrained SPs as well as the required iterations on average are listed in Table \ref{tab:all index SPs}. The configurations of eigenvalues of $A$ corresponding to the function values are also included.
\begin{table}[!t]
	\centering
	\caption{Indices, function values, the corresponding configurations of eigenvalues of $A$, and the needed iterations on average for finding constrained SPs on the linear eigenvalue problem ($n=10$, $p=2$, and $\xi=1.01$). The configurations are indicated by doublets; e.g., (1, 3) means that the function value equals $\frac12(\sigma_1+\sigma_3)$. \revise{For each configuration, we also provide its occurrence in 400 random trials.}}
	\label{tab:all index SPs}
	\resizebox{.9\linewidth}{!}{\begin{tabular}{|c|cccr||c|cccr|}
			\hline
			Indices & Func. vals. & Configs. & Iters. & Occs. & Indices & Func. vals. & Configs. & Iters. & Occs. \\\hline
			0 & 0.918912 & (1, 2) & 148.3 & 400 & 8 & \tabincell{c}{0.956223\\0.956318\\0.956507\\0.956791\\0.957170} & \tabincell{c}{(5, 6)\\(4, 7)\\(3, 8)\\(2, 9)\\(1, 10)} & \tabincell{c}{156.4\\162.2\\164.0\\162.5\\169.0} & \tabincell{r}{71\\198\\113\\17\\1}\\\hline
			
			1 & 0.923529 & (1, 3) & 161.6 & 400 & 9 & \tabincell{c}{0.961028\\0.961171\\0.961409\\0.961742} & \tabincell{c}{(5, 7)\\(4, 8)\\(3, 9)\\(2, 10)} & \tabincell{c}{161.0\\160.6\\158.9\\155.8} & \tabincell{r}{181\\163\\50\\6}\\\hline
			
			2 & \tabincell{c}{0.928101\\0.928193} & \tabincell{c}{(2, 3)\\(1, 4)} & \tabincell{c}{161.9\\162.0} & \tabincell{r}{126\\274} & 10 & \tabincell{c}{0.965785\\0.965881\\0.966072\\0.966359} & \tabincell{c}{(6, 7)\\(5, 8)\\(4, 9)\\(3, 10)} & \tabincell{c}{156.3\\160.5\\161.2\\159.6} & \tabincell{r}{79\\212\\97\\12}\\\hline
			
			3 & \tabincell{c}{0.932764\\0.932903} & \tabincell{c}{(2, 4)\\(1, 5)} & \tabincell{c}{165.3\\162.4} & \tabincell{r}{221\\179} & 11 & \tabincell{c}{0.970638\\0.970782\\0.971023} & \tabincell{c}{(6, 8)\\(5, 9)\\(4, 10)} & \tabincell{c}{159.7\\160.5\\160.8} & \tabincell{r}{201\\166\\33}\\\hline
			
			4 & \tabincell{c}{0.937382\\0.937474\\0.937660} & \tabincell{c}{(3, 4)\\(2, 5)\\(1, 6)} & \tabincell{c}{161.1\\163.8\\163.1} & \tabincell{r}{87\\221\\92} & 12 & \tabincell{c}{0.975443\\0.975540\\0.975733} & \tabincell{c}{(7, 8)\\(6, 9)\\(5, 10)} & \tabincell{c}{153.7\\157.7\\158.4} & \tabincell{r}{97\\229\\74}\\\hline
			
			5 & \tabincell{c}{0.942092\\0.942232\\0.942465} & \tabincell{c}{(3, 5)\\(2, 6)\\(1, 7)} & \tabincell{c}{163.7\\163.0\\158.8}& \tabincell{r}{196\\162\\42} & 13 & \tabincell{c}{0.980345\\0.980490} & \tabincell{c}{(7, 9)\\(6, 10)} & \tabincell{c}{157.4\\154.5} & \tabincell{r}{266\\134}\\\hline
			
			6 & \tabincell{c}{0.946755\\0.946849\\0.947037\\0.947318} & \tabincell{c}{(4, 5)\\(3, 6)\\(2, 7)\\(1, 8)} & \tabincell{c}{159.9\\163.9\\164.5\\163.6} & \tabincell{r}{56\\225\\105\\14} & 14 & \tabincell{c}{0.985198\\0.985295} & \tabincell{c}{(8, 9)\\(7, 10)} & \tabincell{c}{150.9\\154.1} & \tabincell{r}{136\\264}\\\hline
			
			7 & \tabincell{c}{0.951513\\0.951654\\0.951890\\0.952219} & \tabincell{c}{(4, 6)\\(3, 7)\\(2, 8)\\(1, 9)} & \tabincell{c}{162.3\\162.3\\163.3\\159.3} & \tabincell{r}{169\\172\\50\\9} & 15 & 0.990148 & (8, 10) & 152.3 & 400\\\hline
	\end{tabular}}
\end{table}
\revise{We observe that every run converged to a valid critical point, and every admissible configuration has been observed at least once (cf. Example \ref{exm:condition number at index 1 Grassmann}), although the sampling frequencies are highly nonuniform.} 
Note that the ``index-0'' constrained SP for this example \revise{is} exactly the global minimizer. 

\subsection{Electronic excited-state calculations}\label{subsec:electronic excited-state calcs}

\par The core of electronic calculations for molecular systems is the electronic Schr\"odinger equation (ESE) \cite{schrodinger1926undulatory}, which is in fact a linear eigenvalue problem. Nevertheless, the ESE is intractable in general due to the curse of dimensionality. For numerical purpose\revise{s}, various approximations have been proposed in some atomic-orbital basis, such as the full configuration interaction (FCI), Hartree-Fock (HF) methods, and post-HF methods \cite{helgaker2000molecular}, among others. Excited states define the optical and reaction properties of atoms and molecules \cite{balzani2014photochemistry,mai2020molecular,turro2009principles}. Characterizing the excited states is challenging due to electron correlation effects \cite{gonzalez2012progress}. Finding the constrained SPs of the quantum chemical approximated methods in use arises as a natural methodology \cite{burton2022energy,burton2020energy,cances2006computing,lewin2004solutions,marie2023excited,saade2024excited}. Moreover, these approximations usually come together with manifold structures. In the following, we briefly introduce the HF methods and report the numerical results of finding constrained SPs as candidates of excited states, with FCI calculations (performed by \texttt{PySCF} \cite{sun2020recent}) as reference. 
More advanced levels of theory, such as the complete active space self-consistent field method \cite{roos1980complete}, involve complicated manifolds \revise{and \textit{ad-hoc} treatments,} and \revise{are investigated in a parallel work; see a recent preprint \cite{grazioli2026critical2} by the authors and other co-workers}. Before proceeding, we remark that our methodology falls into the class of state-specific methods in quantum chemistry for electronic excited states; other popular ones include linear response theory \cite{cances1999time,casida1995time,grazioli2026critical,olsen1985linear} and state-average methods for multiconfigurational approximations \cite{werner1981quadratically}. 

\medskip

\par\noindent\textit{The HF methods.} The HF approximation restricts the electronic wavefunction to a single Slater determinant parametrized by orthonormal molecular orbitals \cite{fock1930naherungsmethode,hartree1928wave}. In this work, we consider the restricted HF (in short, RHF) method. The RHF method assumes all molecular orbitals to be doubly occupied, by one spin-up and one spin-down electron. As a result, the spatial orbitals are considered to be the same for both spin-up and spin-down electrons. The spatial orbitals are expressed as \revise{\sout{the}} linear combinations of \revise{the} chosen atomic-orbitals, with coefficients to be determined. The atomic-orbitals are assumed to be real hereafter. 

\par For a closed-shell system with $N_{\rm elec}\in\N$ electrons, the RHF approximation gives rise to the following energy functional on the Grassmann manifold:
$$E^{\RHF}(\gamma):=2\trace(h\gamma)+\trace((2J(\gamma)-K(\gamma))\gamma)\revise{+E_{\rm nuc}}\quad\text{with}\quad\gamma\in\Gr_{N_{\rm o}}(\R^{N_{\rmb}}),$$
where $N_{\rmb}\in\N$ is the size of real atomic-orbital basis $\{\phi_i\}_{i=1}^{N_{\rmb}}$, $N_{\rm o}\in\N$ the number of occupied molecular orbitals ($2N_{\rm o}=N_{\rm elec}$), $h\in\R_{\sym}^{N_{\rmb}\times N_{\rmb}}$ the discretized one-body Hamiltonian, \revise{$E_{\rm nuc}\in\R$ the nuclear-repulsion energy,} and $\gamma$ the discretized one-body reduced density matrix. Here, $J$, $K:\R_{\sym}^{N_{\rmb}\times N_{\rmb}}\to\R_{\sym}^{N_{\rmb}\times N_{\rmb}}$ are respectively the Coulomb and exchange functionals, defined as
$$[J(\gamma)]_{pq}:=\sum_{r,s=1}^{N_{\rmb}}g_{pqrs}\gamma_{sr},\quad[K(\gamma)]_{pq}:=\sum_{r,s=1}^{N_{\rmb}}g_{psrq}\gamma_{sr},\quad p,q=1,\ldots,N_{\rmb},$$
with 
$$g_{pqrs}:=\int_{\R^3}\int_{\R^3}\frac{\phi_p(\vecr)\phi_q(\vecr)\phi_r(\vecr')\phi_s(\vecr')}{\abs{\vecr-\vecr'}}\dd\vecr\dd\vecr',\quad p,q,r,s=1,\ldots,N_{\rmb}$$
the two-body integrals. 

\medskip

\par\noindent\textit{Implementation details.} We consider the H$_2$ and LiH molecules. The bond length of the H$_2$ molecule is varied from 0.1 bohr\footnote{1 bohr $\approx$ $5.29\times×10^{-11}$ m.} to 4.0 bohr with a spacing of 0.1 bohr, while that of the LiH molecule is from 2.0 bohr to 6.0 bohr with a spacing of 0.2 bohr\revise{.} The H$_2$ and LiH molecules are described by RHF with the 6-31G basis set ($N_{\rm elec}=2$, $N_{\rmb}=4$, $N_{\rmo}=1$) \cite{ditchfield1971self,hehre1972self} and the STO-3G basis set ($N_{\rm elec}=4$, $N_{\rmb}=6$, $N_{\rm o}=2$) \cite{hehre1969self}, respectively. The numbers of total degrees of freedom (DOFs) are thus $N_{\rmo}(N_{\rmb}-N_{\rmo})=3$ for the H$_2$ molecule and 8 for the LiH molecule. We investigate the energy landscape of RHF for both molecules with the varying bond lengths by searching for the constrained SPs of allowed indices. They are found by running Algorithm \ref{alg:discretized deterministic constrained saddle dynamics lifted} from 1,000 random initial points (random seed $=0\sim999$); see the pseudocodes in Section \ref{subsec:linear eigenvalue problem} for \revise{the} initialization. \revise{Regarding} retraction, Algorithm \ref{alg:discretized deterministic constrained saddle dynamics lifted} is equipped with the exponential mapping for the position part; the treatment for the direction part is similar to that in the previous subsection. The step size is specified as $\eta=10^{-1}$ and $\eta=10^{-2}$ for the H$_2$ and LiH molecules, respectively, which are not necessarily optimal. The maximum iteration number and convergence tolerance are respectively set as \revise{$maxit=20000$} and $tol=10^{-6}$. The results are compared with those obtained by solving FCI based on the RHF ground-state calculations, which is exact under the basis in use. In particular, FCI gives 16 states for the H$_2$ molecule and 225 states for the LiH molecule. 

\medskip

\par
\begin{figure}[!t]
	\centering
	\includegraphics[width=0.55\linewidth]{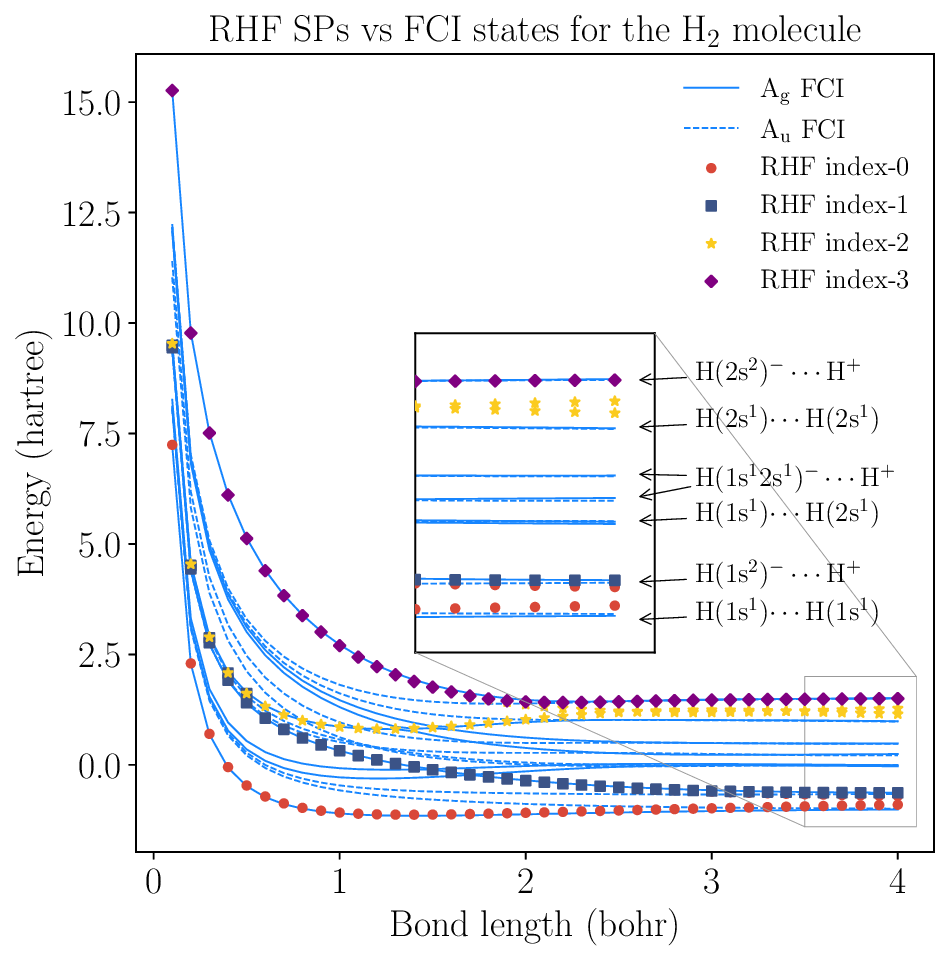}
	\caption{An overview of the constrained SPs of RHF found across the bond length interval on the H$_2$ molecule. The blue solid and dashed lines represent the energies of FCI states belonging to the irreps ${\rm A_g}$ and ${\rm A_u}$, respectively. The red dots, deepblue squares, yellow stars, and purple diamonds stand for the energies of the RHF SPs of indices 0, 1, 2, and 3, respectively. The zoom-in inset describes the FCI states by atomic fragments when the molecule is close to its dissociation limit.}
	\label{fig:H2 RHF FCI all}
\end{figure}

 \noindent\textit{Results on the H$_2$ molecule.} An overview of the constrained SPs of RHF identified across the considered bond length interval is shown in Figure \ref{fig:H2 RHF FCI all}, together with the ground-/excited-state energies of FCI as reference. 
It can be seen that the RHF energy landscape varies smoothly with the bond length. Note that FCI states are classified by their irreducible representations (irreps); \revise{for a complete reference to point group symmetry theory, please refer to \cite{mcweeny2012symmetry}}. For the H$_2$ molecule in the computational point group ${\rm D_{2h}}$, the relevant irreps are ${\rm A_{g}}$ and ${\rm A_{u}}$, with eight FCI states belonging to each. It is observed that the RHF SPs are only able to describe a small subset (e.g., four near the equilibrium geometry, three close to the dissociation limit) of FCI states in terms of their energies. This behavior is consistent with the fact that RHF neglects electronic correlation \cite{barca2014communication}. The deficiency can be mitigated to some extent by resorting to post-HF methods, which is investigated in a parallel work. 

\par
\begin{figure}[!t]
	\centering
	\includegraphics[width=.9\linewidth]{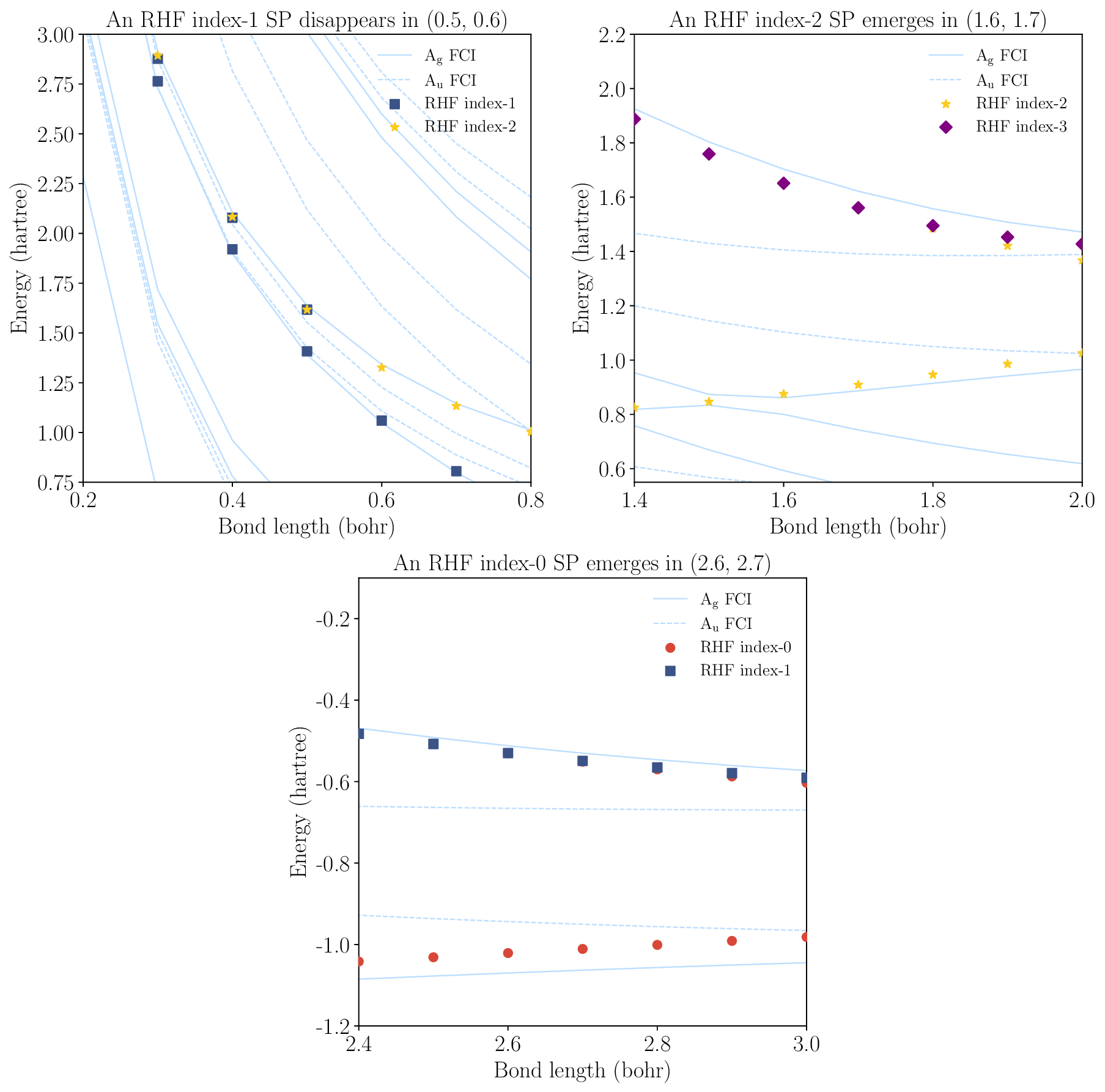}
	
	\caption{Zoom-in views of the constrained SPs of RHF. Upper left: an RHF index-1 SP disappears in (0.5 bohr, 0.6 bohr). Upper right: an RHF index-2 SP emerges in (1.6 bohr, 1.7 bohr). Lower: an RHF index-0 SP emerges in (2.6 bohr, 2.7 bohr).}
	\label{fig:H2 RHF FCI zoom in}
\end{figure}

 Our results also reveal numerically that the varying bond length, as an external parameter, can lead to disappearance or emergence of constrained SPs. Zoom-in views are given in Figure \ref{fig:H2 RHF FCI zoom in}. 
Concretely, in the interval (0.5 bohr, 0.6 bohr), an index-1 SP gets close to the index-2 SP, in terms of energy, and disappears; in the interval (1.6 bohr, 1.7 bohr), an index-2 SP close to the index-3 SP emerges; and in the interval (2.6 bohr, 2.7 bohr), an index-0 SP (i.e., a local minimizer) close to the index-1 SP emerges. 

\medskip

\par
\begin{figure}[!t]
	\centering
	\includegraphics[width=0.9\linewidth]{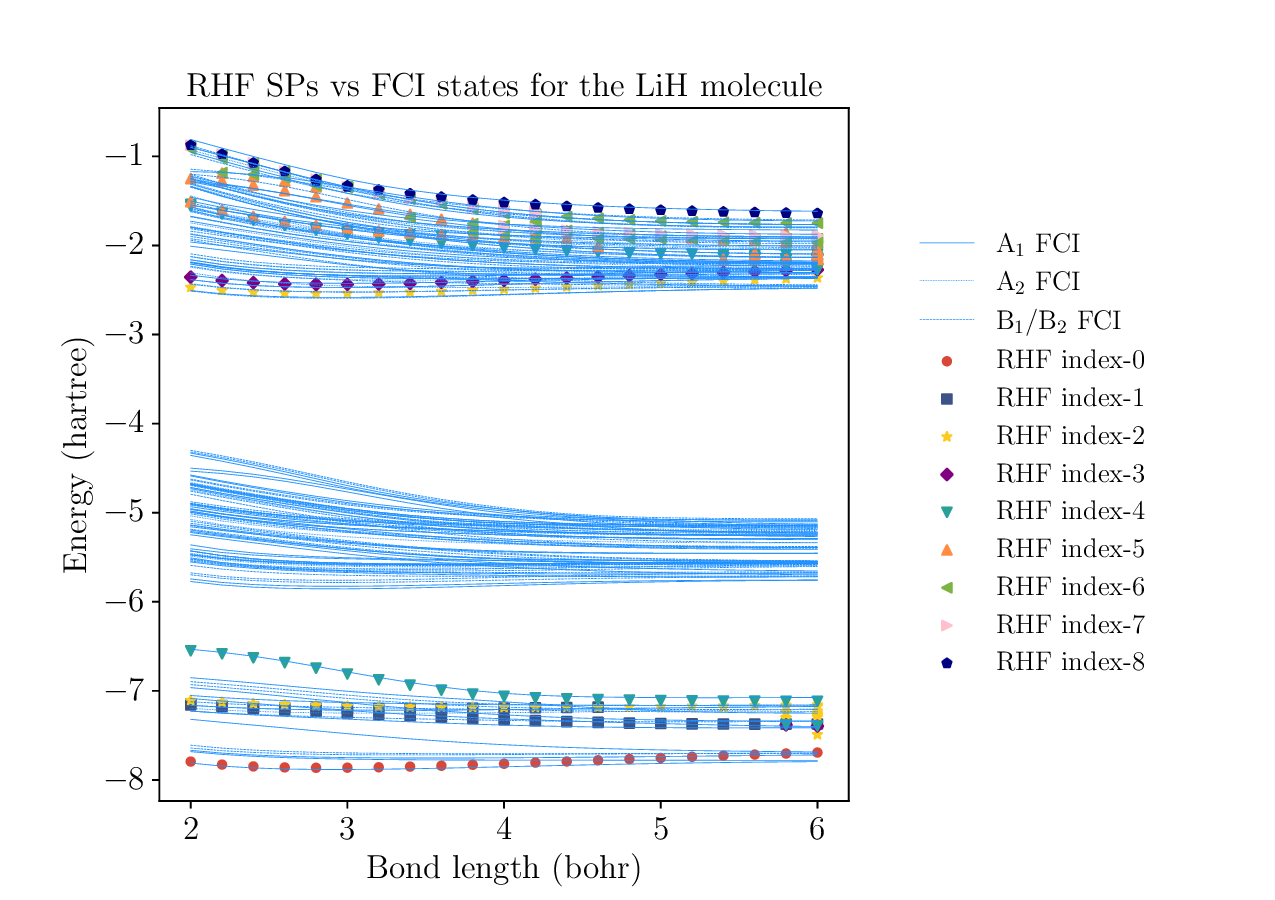}
	\caption{An overview of the constrained SPs of RHF found across the bond length interval on the LiH molecule. The blue solid, dotted, and dashed lines represent the energies of FCI states belonging to the irreps A$_1$, A$_2$, B$_1$ (or B$_2$), respectively. The red dots, deepblue squares, yellow stars, purple diamonds, down-pointing cyan triangles, up-pointing orange triangles, left-pointing green triangles, right-pointing pink triangles, and navyblue pentagons stand for the energies of the RHF SPs of indices 0 to 8, respectively.}
	\label{fig:LiH RHF FCI all}
\end{figure}

\noindent\textit{Results on the LiH molecule.} The identified constrained SPs of RHF are shown in a similar way in Figure \ref{fig:LiH RHF FCI all}, again with the FCI-state energies as reference. 
For the LiH molecule in the computational group $\rm C_{2v}$, the relevant irreps are A$_1$, A$_2$, B$_1$, and B$_2$. Note that the B$_1$ and B$_2$ FCI states are two-fold degenerate across the considered bond length interval. In this case, the RHF SPs appear to well describe the low-lying FCI states with energies lower than $-6$ hartree\footnote{1 hartree $\approx4.36\times10^{-18}$ J}. But the FCI states with energies in ($-6$ hartree, $-4$ hartree) are completely missed by the RHF SPs, probably due to inherent shortcomings of the RHF approximation. In addition, the potential energy surfaces of higher-index RHF SPs are found to be more complicated than those of low-index SPs, especially when the molecule is driven to the dissociation limit. 

\medskip

\par The above results provide proof-of-concept evidence for the effectiveness of our algorithms in excited-state calculations. Nonetheless, we shall point out that a comprehensive quantum chemical analysis of the obtained SPs, though beyond the scope of present work, is essential for practical applications. As a nonlinear approximation to the exact theory, RHF may produce critical points that lack physical meanings; e.g., a spurious non-global local minimizer emerges when the bond length exceeds 2.7 bohr in the right panel of Figure \ref{fig:H2 RHF FCI zoom in}. Moreover, since RHF is a low-dimensional approximation, a one-to-one correspondence between RHF SPs and excited states in the same energetic order no longer holds. Molecular geometries must also be taken into account. In the exact theory, the ESE is parameterized by atomic configurations. Variations in molecular geometry may induce degeneracies and crossings in states, which would lead to unfavorable features of potential energy surfaces within approximate theories (cf. Figure \ref{fig:H2 RHF FCI zoom in}). From the computational aspect, it is of great importance to develop schemes capable of efficiently navigating the nonconvex landscape \cite{xu2025general}, instead of relying on random multi-start. 

\section{Conclusions}\label{sec:conclusions}

\par We have developed a constrained saddle dynamic\revise{\sout{s}} for finding SPs on \revise{any embedded Riemannian submanifold of Euclidean space}. The dynamic\revise{\sout{s}} is formulated compactly on the Grassmann bundle of the tangent bundle, and achieves broad\revise{er} applicability \revise{than those in the literature exploiting explicit global/local level-set manifold representations}. 
By investigating the Grassmann bundle geometry, we have rigorously established the theoretical properties of both the dynamic\revise{\sout{s}} and the resulting discretized algorithms. Remarkably, our analysis provides the \revise{first linear iterate convergence results of the discretized algorithms in embedded-submanifold settings}. Moreover, compared with existing results, we eliminate unnecessary nondegeneracy assumptions on the eigenvalues of the Riemannian Hessian by adopting a single orthogonal projector as the direction variable, thereby respecting the underlying quotient structure. We have also characterized how the spectrum of the Riemannian Hessian affects the local convergence rates and highlighted the importance of using nonredundant parametrizations. Both of these two points have been validated \revise{with} numerical results on linear eigenvalue problems. Finally, we have applied the proposed algorithms to electronic excited-state calculations. 

\par There remain\revise{\sout{s}} lots of directions to be explored. \revise{The constants showing up in the assumptions of convergence results are still not computable, which brings difficulties to choosing step sizes. An adaptively tuned yet theoretically guaranteed step size strategy is much more preferable. Moreover, t}he numerical performance of the discretized algorithms is highly sensitive to condition numbers, as evidenced by their local convergence rates. It is thus desirable to incorporate higher-order contributions without sacrificing local convergence properties. In addition, a globally convergent method for locating SPs on Riemannian manifolds is still lacking, due to the absence of a global merit function. One possible avenue is to extend the analysis in \cite{lelievre2024using} and develop stochastic methods on manifolds. \revise{An intrinsic extension of the formulation in this work to abstract Riemannian manifolds should be possible, but would require replacing the ambient constructions by intrinsic bundle-geometric counterparts, and is left for future work.} From the perspective of quantum chemistry, it would also be valuable to investigate the manifold geometry underlying more complicated levels of theory and to devise efficient yet physically meaningful strategies for navigating the associated nonconvex energy landscapes.

\section*{Acknowledgements}

This work has received funding from the European Research Council (ERC) under the
European Union's Horizon 2020 research and innovation program (grant agreement EMC2
No 810367). The authors are grateful to \'Eric Canc\`es, Tony Leli\`evre, and Panos Parpas for useful discussions. 

\appendix

\section{Fundamental concepts of embedded Riemannian submanifolds of Euclidean space}\label{appsec:preliminaries on Riemannian manifolds}

\par We recall briefly some fundamental concepts of \revise{embedded Riemannian submanifolds of Euclidean space}. For interested readers, we refer to the monographs \cite{absil2008optimization,boumal2023introduction,lee2012introduction,lee2018introduction}. \revise{In the following
, let $\calM$ be embedded in a Euclidean space $\calE$ with $\dim(\calM)=d$. }

\medskip

\par\noindent {\it Tangent space and tangent bundle.} For each $\vecx\in\calM$, the tangent space to $\calM$ at $\vecx$ is referred to as $\upT_{\vecx}\calM$, which is defined as 
$$\upT_{\vecx}\calM:=\lrbrace{c'(0)\mid c:\R\supseteq\calI\to\calM~\text{smooth},~c(0)=\vecx}.$$
The vectors in $\upT_{\vecx}\calM$ are called tangent vectors to $\calM$ at $\vecx$. The tangent space $\upT_{\vecx}\calM$ is endowed with the Riemannian metric induced from the inner product of the ambient Euclidean space. For any $\vecx\in\calM$, the orthogonal projection operator $\Proj_{\upT_{\vecx}\calM}$ from the ambient space onto the tangent space $\upT_{\vecx}\calM$ is defined as
$$\Proj_{\upT_{\vecx}\calM}(\vecv):=\argmin_{\vecu}\norm{\vecu-\vecv},\quad\forall~\vecv\in\calE,$$
where $\snorm{\bullet}:\calE\to\R_+$ refers to the norm induced by the inner product. 

\par The tangent bundle of $\calM$ is denoted by $\upT\calM:=\{(\vecx,\vecv)\mid\vecx\in\calM,~\vecv\in\upT_{\vecx}\calM\}$, i.e., the disjoint union of the tangent spaces to $\calM$. From \cite[Theorem 3.43]{boumal2023introduction}, $\upT\calM$ is a $2d$-dimensional submanifold embedded in $\calE\times\calE$. 

\medskip

\par\noindent {\it Riemannian gradient and Hessian.} For a smooth function $f$, its Riemannian gradient at $\vecx\in\calM$, denoted by $\grad{\calM}f(\vecx)$, is defined as the unique element of $\upT_{\vecx}\calM$ satisfying
$$\inner{\grad{\calM}f(\vecx),\vecv}=\upD f(\vecx)[\vecv],\quad\forall~\vecv\in\upT_{\vecx}\calM,$$
where $\upD f(\vecx)[\vecv]$ stands for the directional derivative of $f$ at $\vecx$ along the tangent vector $\vecv$. Since $\calM$ is a Riemannian submanifold embedded in a Euclidean space, $\grad{\calM}f(\vecx)$ can be readily computed via
$$\grad{\calM}f(\vecx)=\Proj_{\upT_{\vecx}\calM}\big(\grad{}f(\vecx)\big).$$

\par The definition of Riemannian Hessian in the general cases necessitates the concept of Riemannian connection \cite[Section 5.4]{boumal2023introduction}. Again due to the fact that $\calM$ is an embedded Riemannian submanifold in $\calE$, we recall for simplicity the following characterization:
\begin{equation}
	\hess{\calM}f(\vecx)[\vecv]=\Proj_{\upT_{\vecx}\calM}\big(\upD\bar G(\vecx)[\vecv]\big),
	\label{eqn:definition of Riemannian Hessian}
\end{equation}
where $\bar G$ is any smooth extension of $\grad{\calM}f$ to a neighborhood of $\calM$ in $\calE$. 

\medskip

\par\noindent {\it Retraction.} A retraction on $\calM$ is a smooth mapping $\upR:\upT\calM\to\calM$, $\upT\calM\owns(\vecx,\vecv)\mapsto\upR_{\vecx}(\vecv)\in\calM$, satisfying $\upR_{\vecx}(0)=\vecx$ and that $\upD\upR_{\vecx}(0)$ is the identity mapping on $\upT_{\vecx}\calM$ for any $\vecx\in\calM$. By leveraging the retraction, we can obtain a point by moving away from $\vecx\in\calM$ along some $\vecv\in\upT_{\vecx}\calM$, while remaining on $\calM$. \revise{It} follows that \revise{retractions define update rules that preserve feasibility}. One typical example of retraction is the exponential mapping, denoted specially by $\Exp:\upT\calM\to\calM$, which is determined by a set of second-order ordinary differential equations and yields geodesics on $\calM$. If $\calM$ is complete, then it holds that \cite[Proposition 10.22]{boumal2023introduction}
\begin{equation*}
	\dist_{\calM}(\vecx,\Exp_{\vecx}(\vecv))=\snorm{\vecv},
\end{equation*}
\revise{for any $\vecx\in\calM$ and $\vecv$ in a neighborhood of the origin of $\upT_{\vecx}\calM$,} where $\dist_{\calM}(\bullet,\bullet):\calM\times\calM\to\R_+$ represents the Riemannian distance on $\calM$:
$$\dist_{\calM}(\vecx,\vecy):=\inf_c\lrbrace{\left.\int_a^b\snorm{c'(t)}\dd t~\right|c~\text{piecewise smooth on}~\calM,~c(a)=\vecx,~c(b)=\vecy}.$$

\medskip

\par\noindent \revise{{\it Parallel transport and vector transport.} 
Let \(c:[0,1]\to\calM\) be a smooth curve. For any \(\vecv\in \upT_{c(0)}\calM\), there exists a unique vector field \(W(t)\in\upT_{c(t)}\calM\) along \(c\) satisfying
\[
\nabla_{\dot c(t)}W(t)=0,
\qquad
W(0)=v.
\]
The parallel transport of \(\vecv\) along \(c\) from \(c(0)\) to \(c(t)\) is defined by
\[
\PT^{c}_{t\leftarrow0}(\vecv):=W(t).
\]
Parallel transport is a linear isometry between the corresponding tangent spaces. Since \(\calM\) is embedded in a Euclidean space, the ambient derivative of a vector field $W$ along \(c\) satisfies the Gauss formula
\[
\frac{\dd W}{\dd t}
=
\nabla_{\dot c}W+\two_{c(t)}(\dot c,W).
\]
Consequently, if \(W\) is parallel along \(c\), then
\[
\frac{\dd W}{\dd t}
=
\two_{c(t)}(\dot c,W).
\]

\par A vector transport associated with a retraction \(\upR\) on $\calM$ is a smooth mapping that assigns, to \(\vecx\in\calM\), \(\vecu,\vecv\in\upT_{\vecx}\calM\), a vector
\[
\VT_{(\vecx,\vecu)}(\vecv)
\in\upT_{\upR_{\vecx}(\vecu)}\calM.
\]
It is required to be linear in \(\vecv\) and to satisfy
\[
\VT_{(\vecx,0)}(\vecv)=\vecv.
\]
Unlike parallel transport, a general vector transport need not be isometric or globally invertible. A common choice is the metric-projection transport
\[
\VT_{(x,\vecu)}(\vecv)
:=
\Proj_{\upT_{\upR_{\vecx}(\vecu)}\calM}(\vecv).
\]}

\medskip

\par\noindent {\it Second fundamental form.} For a given $\vecx\in\calM$, the normal space $\upN_{\vecx}\calM$ is the orthogonal complement of $\upT_{\vecx}\calM$ in $\calE$. The second fundamental form at $\vecx$ can be identified as the mapping $\two_{\vecx}:\upT_{\vecx}\calM\times\upT_{\vecx}\calM\to\upN_{\vecx}\calM$ defined as $\two_{\vecx}(\vecu,\vecv):=\upP_{\vecx,\vecu}(\vecv)$ for any $\vecu$, $\vecv\in\upT_{\vecx}\calM$, where
\begin{equation}
	\upP_{\vecx,\vecu}:=\upD(\vecy\mapsto\Proj_{\upT_{\vecy}\calM})(\vecx)[\vecu],
	\label{eqn:differential of orthogonal projector}
\end{equation}
namely, the directional derivative of $\vecy\mapsto\Proj_{\upT_{\vecy}\calM}$ at $\vecx$ along $\vecu$. That the range of $\two_{\vecx}$ is $\upN_{\vecx}\calM$ can be shown by taking the directional derivative on the both sides of the identity $\Proj_{\upT_{\vecx}\calM}\circ\Proj_{\upT_{\vecx}\calM}=\Proj_{\upT_{\vecx}\calM}$, which yields
\begin{equation}
	\upP_{\vecx,\vecu}\circ\Proj_{\upN_{\vecx}\calM}=\Proj_{\upT_{\vecx}\calM}\circ\upP_{\vecx,\vecu},\quad\upP_{\vecx,\vecu}\circ\Proj_{\upT_{\vecx}\calM}=\Proj_{\upN_{\vecx}\calM}\circ\upP_{\vecx,\vecu}.
	\label{eqn:exchange property of projection differential}
\end{equation}
By the above definition, the second fundamental form accounts for the changes in the way the tangent spaces sit inside the ambient $\calE$. 

\revise{
\section{Numerical comparisons between different retractions}\label{appsec:numerical comparisons between retractions}

We test Algorithm \ref{alg:discretized deterministic constrained saddle dynamics lifted} equipped with different choices of retractions on the linear eigenvalue problem. For the problem description and most implementation details, please refer to Section \ref{subsec:linear eigenvalue problem}. We fix $n=80$, $p=8$, and $\xi=1.05$, and search for the index-1 constrained SP. The initial points are generated either with random seed $=1$ or from perturbation of the index-1 constrained SP (with the perturbation level $\beta=10^{-1}$) with random seed $=1$. We adopt the estimated best step sizes in Eqs. \eqref{eqn:theoretical step size x} and \eqref{eqn:theoretical step size v}.

\par A general construction of retraction on $\Gr_k(\upT\calM)$ is given in Lemma \ref{lem:general retraction on Grassmann bundle}. For the linear eigenvalue problem, we have $\calM=\Gr_p(\R^n)$ and take $\upR^{\calM}=\Exp^{\calM}$, whose explicit form can be found in, e.g., \cite{bendokat2024grassmann}. For the other two mappings showing up in (ii) and (iii) of Lemma \ref{lem:general retraction on Grassmann bundle}, we consider the concrete examples provided in Remark \ref{rem:example retractions}, which amount to four combinations in Table \ref{tab:retractions to be tested}. 
\begin{table}[htb]
    \centering
    \caption{\revise{Retractions in comparison.}}
    \label{tab:retractions to be tested}
    \begin{tabular}{|l|cc|}
        \hline
        Combinations & $\widehat\upR^{\Gr_k(\upT_{\vecx}\calM)}$ & $\VT^{\calM}$ \\\hline
        Exp + PT & Exponential mapping & Parallel transport \\
        Exp + Proj & Exponential mapping & Metric projection\\
        Add + PT & Direct addition & Parallel transport \\
        Add + Proj & Direct addition & Metric projection\\
        \hline
    \end{tabular}
\end{table}

\par In experiments, all these combinations converge to the same solution from two different initializations, with the converged objective value 0.101898. As noted at the end of Remark \ref{rem:example retractions}, these combinations differ from the second order and can yield distinct numerical performance. Their convergence curves in wall-clock time on the linear eigenvalue instance are collected in Figure \ref{fig:retraction comparison}. 
\begin{figure}[!t]
    \centering
    \includegraphics[width=.9\linewidth]{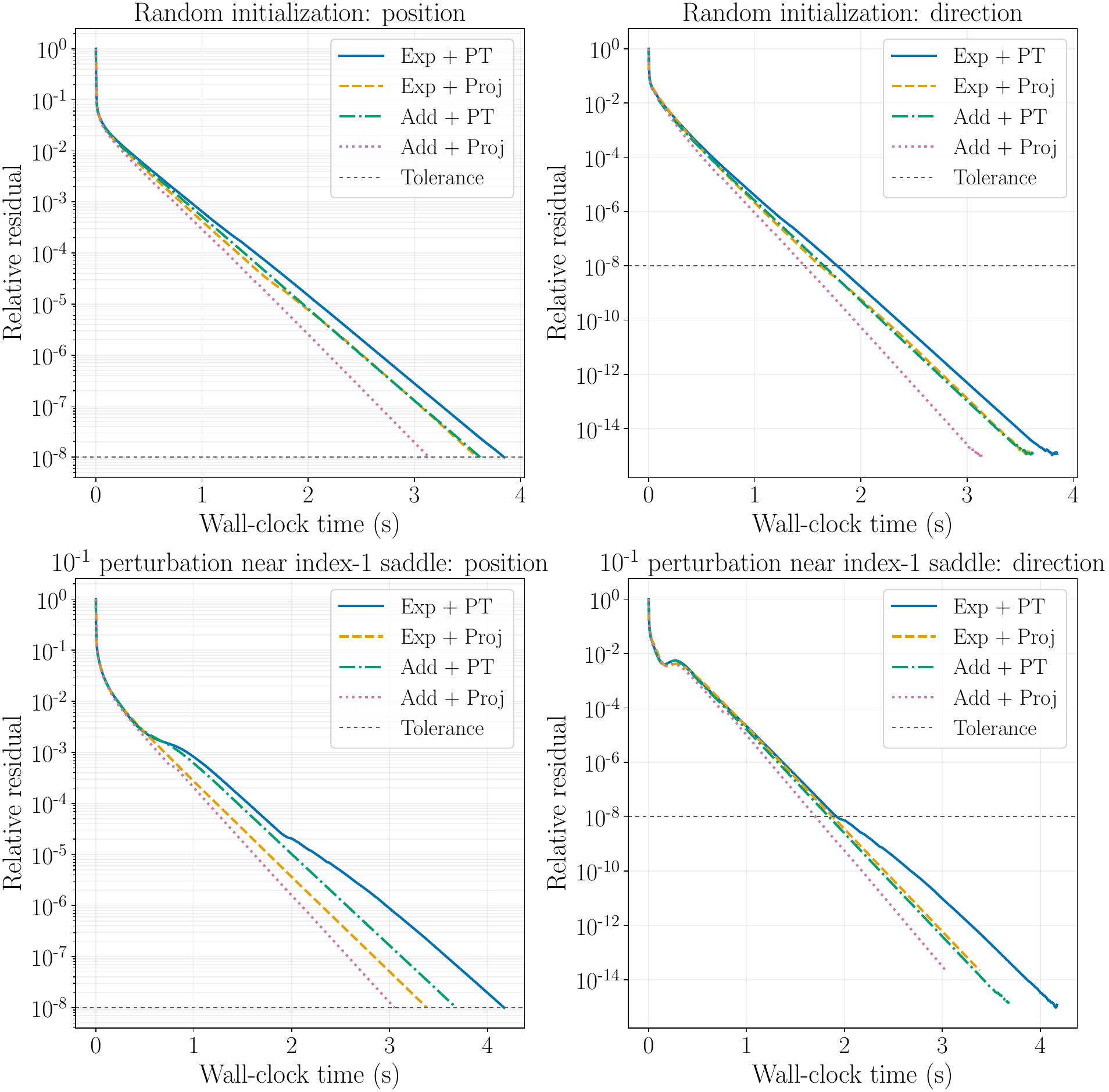}
    \caption{\revise{Convergence curves of Algorithm \ref{alg:discretized deterministic constrained saddle dynamics lifted} implementing different retractions on the linear eigenvalue problem ($n=80$, $p=8$, and $\xi=1.05$). The blue solid, orange dashed, green dash-dotted, and red dotted lines denote the results by the combinations ``Exp + PT'', ``Exp + Proj'', ``Add + PT'', ``Add + Proj'', respectively. The black dashed line refers to the tolerance $10^{-8}$. Top: random initialization. Bottom: initialization from $10^{-1}$ perturbation near the index-1 SP. Left: position relative residual. Right: direction relative residual.}}
    \label{fig:retraction comparison}
\end{figure}
The results showcase that the combination ``Add + Proj'' in Table \ref{tab:retractions to be tested} is the most efficient in time. This can be ascribed to the fact that the operations of direction addition and metric projection are cheaper than those of exponential mapping and parallel transport. We note that the implementations in Section \ref{sec:numerical experiments} all adopt ``Add + Proj''. In terms of iteration number, the combinations ``Exp + PT'', ``Exp + Proj'', ``Add + PT'', ``Add + Proj'' are only slightly different from each other, requiring respectively 5148, 5160, 5146, 5157 iterations with random initialization, and respectively 5288, 4943, 5214, 4961 iterations with perturbed initialization. 
}


\end{document}